\documentclass{amsart}	

\usepackage{amsmath}
\usepackage{amssymb}
\usepackage{amsthm}
\usepackage{mathrsfs}
\usepackage[all]{xy} 
\usepackage{bm}
\usepackage{comment}
\usepackage{color}
\usepackage{manfnt}
\usepackage{indentfirst}% Indent the first paragraph

\usepackage{autobreak}
\usepackage{tikz}
    \usetikzlibrary{arrows.meta} 

\usepackage{mathtools}
\usepackage{array,mathdots}
\usepackage[T1]{fontenc}

\usepackage{cleveref}
\usepackage{autonum} 

\usepackage{aliascnt}

\newtheorem{Theorem}{Theorem}[section]

\newaliascnt{Lemma}{Theorem}
\newtheorem{Lemma}[Lemma]{Lemma}
\aliascntresetthe{Lemma}

\newaliascnt{Proposition}{Theorem}
\newtheorem{Proposition}[Proposition]{Proposition}
\aliascntresetthe{Proposition}

\newaliascnt{Corollary}{Theorem}
\newtheorem{Corollary}[Corollary]{Corollary}
\aliascntresetthe{Corollary}

\crefname{Theorem}{Theorem}{Theorems}
\crefname{Lemma}{Lemma}{Lemmas}
\crefname{Proposition}{Proposition}{Propositions}
\crefname{Corollary}{Corollary}{Corollaries}

\theoremstyle{definition}
\newaliascnt{definition}{Theorem}
\newtheorem{definition}[definition]{Definition}
\aliascntresetthe{definition}

\newaliascnt{Remark}{Theorem}
\newtheorem{Remark}[Remark]{Remark}
\aliascntresetthe{Remark}

\newaliascnt{Example}{Theorem}
\newtheorem{Example}[Example]{Example}
\aliascntresetthe{Example}

\newaliascnt{problem}{Theorem}

\aliascntresetthe{problem}

\crefname{definition}{Definition}{Definitions}
\Crefname{definition}{Definition}{Definitions}

\crefname{Remark}{Remark}{Remarks}
\Crefname{Remark}{Remark}{Remarks}

\crefname{Example}{Example}{Examples}
\Crefname{Example}{Example}{Examples}

\crefname{problem}{Problem}{Problems}
\Crefname{problem}{Problem}{Problems}

\crefname{equation}{Equation}{Equations}
\Crefname{equation}{Equation}{Equations}

\numberwithin{equation}{section}

\DeclareMathOperator{\trdeg}{\mathrm{tr\text{.}deg}}

\crefname{Theorem}{Theorem}{Theorems}
\crefname{Lemma}{Lemma}{Lemmas}
\crefname{Proposition}{Proposition}{Propositions}
\crefname{Corollary}{Corollary}{Corollaries}

\newcommand{\Aloc}{\ensuremath A_{(\zeta)}}
\newcommand{\Acomp}{\ensuremath A_{v}}
\newcommand{\Kcomp}{\ensuremath K_{v}}
\newcommand{\Kalg}{\ensuremath {K^{\mathrm{alg}}}}
\newcommand{\quotfield}{\ensuremath \mathbb{F}_v}
\newcommand{\GL}{\ensuremath \operatorname{GL}}
\newcommand{\Mat}{\ensuremath \operatorname{Mat}}
\newcommand{\Omegav}{\ensuremath \Omega_v}

\newcommand{\prelog}{\ensuremath \mathfrak{L}}
\newcommand{\omegav}{\ensuremath {{\omega_{v}}}}
\newcommand{\rank}{\ensuremath {\operatorname{rank}}}
\newcommand{\Entire}{\ensuremath {\mathbb{E}_\infty}}
\newcommand{\Entirev}{\ensuremath {\mathbb{E}_v}}
\newcommand{\Frac}{\ensuremath {\operatorname{Frac}}}

\newlength\PowW
\newcolumntype{P}{>{\mathmakebox[\PowW][c]}c}

\allowdisplaybreaks[4]

\title[Algebraic independence of $v$-adic Carlitz multiple polylogarithms]{Algebraic independence of $v$-adic Carlitz multiple polylogarithms via period interpretations}

\author{Daichi Matsuzuki}
\address{Department of Mathematics, Texas A\&M University, College Station, TX 77843, U.S.A}
\email{daichi\_matsuzuki@tamu.edu}
\date{\today}

\begin{document}
\begin{abstract}

   Based on Hartl's de Rham--crystalline comparison isomorphism, we formulate a notion of Andr\'{e}'s $v$-adic periods for a class of Anderson $\mathbb{F}_q[z]$-motives with respect to a degree-one finite place $v$ of a rational function field over a finite field, following Ancona--Fr\u{a}\c{t}il\u{a}'s formulation of $p$-adic periods. 
    The formulation naturally produces deformation series of the resulting periods and Frobenius difference equations satisfied by these series. 
    We then explain how the formulation brings the resulting $v$-adic periods within the scope of Chen's $v$-adic ABP criterion and of difference Galois theory, and gives us a strategy to study the transcendence degrees of $v$-adic periods through their deformation series.
   We show that certain special values of $v$-adic Carlitz multiple polylogarithms can be expressed in terms of Andr\'e's $v$-adic periods.
   Applying the strategy to these examples, we obtain algebraic independence results for certain special values of $v$-adic Carlitz multiple polylogarithms;
   more precisely, the transcendence degree of the field generated by a suitable collection of special values of $v$-adic Carlitz multiple polylogarithms is one less than the transcendence degree of the field generated by the corresponding special values of $\infty$-adic Carlitz multiple polylogarithms together with the Carlitz period.

\end{abstract}
\maketitle

\tableofcontents

%%%%%%%%%%%%%%%
\section{Introduction}
\subsection{Background and motivation}

The aim of the present paper is to formulate a positive-characteristic analogue of Andr\'e's $p$-adic periods, following the framework developed in \cite{Ancona2025}, to develop an approach to studying transcendence properties of the resulting periods, and to apply the strategy to the study of examples given by special values of $v$-adic Carlitz multiple polylogarithms.
Ancona and Fr\u{a}\c{t}il\u{a} developed a formulation of Andr\'e's $p$-adic periods and observed that certain special values of the $p$-adic gamma function and $p$-adic logarithms arise as basic examples in their theory;
furthermore, they related algebraic independence properties of $p$-adic periods to the structure of motivic Galois groups and proposed a $p$-adic analogue of Grothendieck's period conjecture.
This point of view has also inspired subsequent work; for instance, we refer the reader to \cite{DanCohen2025} for Andr\'e's $p$-adic periods of mixed Tate motives.
The present paper focuses on a finite place $v$ of degree $1$ of the rational function field over a finite field. 
In this setting, we shall introduce a notion of Andr\'e's $v$-adic periods and show that certain special values of $v$-adic Carlitz multiple polylogarithms (\cite{Chang2019a}) admit expressions in terms of such periods; further, we obtain algebraic independence results for those values based on the formulation of Andr\'e's $v$-adic periods and the interpretations of the special values.
We hope that this provides a first step toward a positive-characteristic counterpart of the theory of Ancona and Fr\u{a}\c{t}il\u{a}.

%\subsubsection{Key idea from the theory of $p$-adic periods}
Andr\'{e}'s $p$-adic periods come from Berthelot's comparison isomorphisms between de Rham and crystalline realizations of motives, whereas the usual periods (see \cite{Kontsevich2001} for example) arise from the comparison isomorphisms between de Rham and Betti realizations. 
Unlike de Rham and Betti realizations, crystalline realizations of motives are not naturally equipped with a canonical $\overline{\mathbb{Q}}$-structure, at least to the best of the author's knowledge; this absence makes it harder to formulate $p$-adic periods and study their transcendence properties.
One classical way to address this difficulty is to use the Frobenius endomorphism on the crystalline realization, as in the approach used in \cite{Brown2017} and \cite{Furusho2007}.
Ancona and Fr\u{a}\c{t}il\u{a} proceed differently; they use cycles on the reduction of the motive modulo $p$ to define $p$-adic periods.
In the present paper, we adapt the latter method to the function-field setting and formulate a corresponding notion of $v$-adic periods for degree-one finite places.

\subsubsection{The function-field setting}
Fixing a prime $p$, a power $q$ of $p$, and an indeterminate $\theta$, we consider the one-variable polynomial ring $A:=\mathbb{F}_q[\theta]$ over the finite field with $q$ elements to be an analogue of the ring $\mathbb{Z}$ of integers. 
Then, the field of fractions $\mathbb{F}_q(\theta)$ of $A$, which is denoted by $K$, can be regarded as an analogue of the field $\mathbb{Q}$; we write $\Kalg$ for an algebraic closure of $K$.
The field of Laurent series $\mathbb{F}_q((1/\theta))$, which is denoted by $K_\infty$, is considered to be an analogue of the field $\mathbb{R}$ of real numbers and the completion $\mathbb{C}_\infty:=\widehat{K_\infty^{\textrm{alg}}}$ of the algebraic closure of $K_\infty$ is regarded as an analogue of the field $\mathbb{C}$ of complex numbers.
Let us take a degree-one prime ideal $v$ of $A$ with a monic generator $\zeta$ and write $\Aloc$ for the localization of $A$ at $v$.
We write $\Acomp$ and $\Kcomp$ for the $v$-adic completions of $A$ and $K$, respectively. Furthermore, we put $\mathbb{C}_v:=\widehat{K_v^{\textrm{alg}}}$.
Let $t$ be an indeterminate independent of $\theta$, and put $z:=\zeta|_{\theta=t}\in\mathbb{F}_q[t]$.
Since $\deg v=1$, we have $\zeta=\theta-\alpha$ and $z=t-\alpha$ for some $\alpha\in\mathbb{F}_q$,  thus 
$\mathbb{F}_q[z]=\mathbb{F}_q[t]$ holds.
We shall also regard $z$ as an indeterminate over $K$, obtained from $t$ by this change of variable.

In the line of research on function-field arithmetic that this paper follows, the role of periods in the sense of \cite{Kontsevich2001} is played by periods of Drinfeld modules, $t$-modules, and $t$-motives (\cite{Anderson1986,Papanikolas2008}).
Drinfeld (\cite{Drinfelprimed1974}) introduced objects called Drinfeld modules as positive characteristic analogues of elliptic curves; the $t$-modules, which were introduced by Anderson (\cite{Anderson1986}), are higher dimensional versions of Drinfeld modules and are analogues of Abelian varieties.
Anderson also introduced a notion of $t$-motives inspired by Drinfeld's idea of shtukas. 
The theory of $t$-motives can be seen as a sort of theory of motives for $t$-modules, in the sense that $t$-motives form a linear category and cohomology theories of $t$-modules factor through the category of $t$-motives, for example; see \cite{Hartl2020} for details.

Periods of these objects and their arithmetic properties have been studied extensively.
A fundamental example is the Carlitz period $\tilde{\pi}$, which is a period of the Carlitz module, the most basic example of a Drinfeld module; this quantity is regarded as an analogue of the complex period $2\pi\sqrt{-1}$.
Wade proved that $\tilde{\pi}$ is transcendental over $K$ (\cite{Wade1941}).
Yu's sub-$t$-module theorem \cite[Theorem 0.1]{Yu1997} has been a fundamental tool in the study of $\Kalg$-linear relations among periods of Drinfeld modules and $t$-modules, as well as their logarithms at algebraic points; see \cite{Brownawell2002} for example.
On the other hand, the ABP criterion \cite[Theorem 1.3.2]{Anderson2004} and Papanikolas' theory \cite{Papanikolas2008} form the foundation of the transcendence theory of periods of $t$-motives; see \cite{Chang2010b,Chang2010c} for example.

The study of a $v$-adic counterpart of the theory of periods appears to be emerging.
Although Chen (\cite{Chen2025}) established a $v$-adic version of the ABP criterion, it has been unclear how his criterion is related to the theory of $t$-motives, while the original ABP criterion plays a pivotal role in the study of $\infty$-adic periods of $t$-motives.
Chang, Wei, and Yu (\cite{Chang2025}) studied $v$-adic periods of the $t$-motives associated with higher Carlitz modules.
Related algebraic independence results for special values of the $v$-adic arithmetic gamma function were also obtained by Adam \cite{Adam2023}; see \cite[Remark 1.4.5]{Chang2025} for a comparison with these results.
More recently, $v$-adic versions of geometric gamma values have also been studied (\cite{Wei2026}).
At present, however, there does not seem to be a fully developed $v$-adic counterpart of the existing theory of $\infty$-adic periods.

\subsubsection{Carlitz multiple polylogarithms}
Let us review ($\infty$-adic) Carlitz multiple polylogarithms (\cite[Definition 5.1.1]{Chang2014}); special values of these functions are known to be examples of periods of $t$-motives \cite[\S 5.3]{Chang2014}.
For each tuple $\mathbf{s}:=(s_1,\,\dots,\,s_r)$ of positive integers, we define
\begin{equation}\label{eq:Dinfty}
    \mathbb{D}_{\mathbf{s},\,\infty}:=\left\{(x_1,\,\dots,\,x_r) \in \mathbb{C}_\infty^r\,\middle| \, |x_1|_\infty <| \theta|_\infty^{\frac{s_1q}{q-1}},\,|x_i|_\infty \leq | \theta|_\infty^{\frac{s_iq}{q-1}} \text{ for } 2 \leq i \leq r \right\};
\end{equation}
in this domain, the series
\begin{equation}\label{eq:defCMPL}
    \sum_{i_1>i_2>\cdots>i_r\geq 0}\frac{x_1^{q^{i_1}}\cdots x_r^{q^{i_r}}}{\Big(\left(\theta-\theta^{q}\right)\cdots\left(\theta-\theta^{q^{i_1}}\right)\Big)^{s_1}\cdots\Big(\left(\theta-\theta^{q}\right)\cdots\left(\theta-\theta^{q^{i_r}}\right)\Big)^{s_r}} 
\end{equation}
converges in $\mathbb{C}_\infty$ and the \textit{($\infty$-adic) Carlitz multiple polylogarithm} $\mathrm{Li}_{s_1,\,\dots,\,s_r}(x_1,$ $\dots,$ $x_r)$ is defined to be the function given by this series; we put $\left(\theta-\theta^{q}\right)\cdots\left(\theta-\theta^{q^{i_r}}\right):=1$ if $i_r=0$ by convention.
Special values of $\infty$-adic Carlitz multiple polylogarithms provide an important class of examples to which the aforementioned transcendence theories apply; see \cite{Chang2007,Chang2014,Mishiba2017} among others. 
These special values admit interpretations in terms of both logarithms of suitable $t$-modules and periods of suitable (dual) $t$-motives. 
The former interpretation enables us to apply Yu's sub-$t$-module theorem, whereas the latter brings these values within the framework of the ABP criterion and Papanikolas' theory. 
The period interpretation, which is more closely related to the content of the present paper, gives rise to deformation series satisfying Frobenius difference equations whose specializations recover the special values of Carlitz multiple polylogarithms; the series enable us to apply the ABP criterion and Papanikolas' theory to study the transcendence properties of the special values.

On the $v$-adic side, we can show that the series \eqref{eq:defCMPL} converges in $\mathbb{C}_v$ within the domain
\begin{equation}\label{eq:Dv}
    \mathbb{D}_{\mathbf{s},\,v}^{\textrm{Conv}}:=\left\{(x_1,\,\dots,\,x_r) \in \mathbb{C}_v^r\,\,\middle| \, \,|x_1|_v<1,\,|x_i|_v \leq 1 \text{ for } 2 \leq i \leq r \right\},
\end{equation}
see \cite[\S 2.2]{Chang2019a}; following this reference, we write $\mathrm{Li}_{s_1,\,\dots,\,s_r}(x_1,\,\dots,\,x_r)_v$ for the resulting function and call it the \textit{$v$-adic Carlitz multiple polylogarithm}.
It turns out in the present paper that certain special values of these $v$-adic functions can be expressed using some examples of Andr\'e's $v$-adic periods.
Chang and Mishiba extend the domain of definition of $\mathrm{Li}_{s_1,\,\dots,\,s_r}(x_1,\,\dots,\,x_r)_v$ to a wider domain; we refer to \cite[Definition 4.2.3]{Chang2019a}.
%\begin{equation}
%    \mathbb{D}_{\mathbf{s},\,v}^{\textrm{Def}}:=\left\{(x_1,\,\dots,\,x_r) \in \mathbb{C}_v^r\,\,\middle| \,\, |x_i|_v \leq 1 \text{ for } 1 \leq i \leq r \right\}.
%\end{equation}
On the $v$-adic side, special values of $v$-adic Carlitz multiple polylogarithms have also been related to the theory of $t$-modules and their logarithms (\cite[Theorem 3.3.3 and Lemma 4.2.1]{Chang2019a}).  
Moreover, Chen constructed deformation series satisfying Frobenius difference equations and used them, together with his $v$-adic analogue of the ABP criterion, to obtain linear independence results for special values of $v$-adic Carlitz multiple polylogarithms.
However, in contrast with the $\infty$-adic case, it has been unclear how special values of $v$-adic Carlitz multiple polylogarithms and Chen's deformations relate to the theory of $t$-motives in a form similar to that of the theory of usual periods.

The main contributions of this paper are as follows.  
First, we formulate a notion of Andr\'e's $v$-adic periods for a certain class of Anderson $\mathbb{F}_q[z]$-motives and for a degree-one place $v$.
A key feature of this formulation is that it naturally produces deformation series satisfying the Frobenius difference equations needed to apply Chen's $v$-adic ABP criterion.
Second, we show that the special values of $v$-adic Carlitz multiple polylogarithms at points all of whose coordinates lie in $v \Aloc$ can be expressed in terms of such periods.
This result explains how, for the special values considered here, Chen's deformation series essentially arise from the theory of $t$-motives.
Third, we relate the transcendence properties of special values of $v$-adic Carlitz multiple polylogarithms and those of special values of $\infty$-adic Carlitz multiple polylogarithms.
As applications, we obtain some algebraic independence results for special values of $v$-adic Carlitz multiple polylogarithms.

\subsection{The formulation of $v$-adic periods}

\subsubsection{Difficulty in studying $v$-adic periods and their transcendence}
Before discussing a formulation of $v$-adic periods, let us briefly recall the notion of the usual periods in the arithmetic of positive characteristic function fields. 
In the literature, the transcendence theory of periods relevant to this paper is usually based on the theory of dual $t$-motives and pre-$t$-motives, not on the theory of $t$-motives,
whereas the objects treated in the present paper are what we call \textit{Anderson $\mathbb{F}_q[t]$-motives} (\cref{Def:t-motive}), which are variants of the $t$-motives introduced by Anderson \cite[\S 1.2]{Anderson1986}.
For the definition of Anderson $\mathbb{F}_q[t]$-motives, we follow \cite[Definition 2.3]{Gazda2024}; see also \cite[Definition 2.1]{Hartl2019}.
We recall the notion of periods in terms of Anderson $\mathbb{F}_q[t]$-motives by interpreting the formulation via the tensor anti-equivalence between $\mathbb{F}_q[t]$-motives and dual $\mathbb{F}_q[t]$-motives, variants of dual $t$-motives (\cite[Proposition 2.4.3]{Hartl2020}).

Taking a uniformizable (also referred to as rigid analytically trivial) Anderson $\mathbb{F}_q[t]$-motive $M$ over $\Kalg$ satisfying certain conditions, we can say that the definition of periods of $M$ hinges on the following three ingredients (\cite[\S 2.3]{Hartl2020} and \cite[1.1.4]{Papanikolas2008}).
The first one is the Betti realization 
$
\operatorname{H}^1_{\mathrm{Betti}}(M,\mathbb{F}_q[t])
$,
which is an $\mathbb{F}_q[t]$-module attached to $M$;
the second is the $\Kalg$-vector space
$
\operatorname{H}^1_{\mathrm{dR}}(M,\Kalg)
$
called the de Rham realization associated with $M$ (\cite[\S 2.3.5]{Hartl2020}).
Third, we have a comparison isomorphism
\begin{equation}
\operatorname{H}^1_{\mathrm{Betti}}(M,\mathbb{F}_q[t])\otimes_{\mathbb{F}_q[t],\,t\mapsto \theta}\mathbb{C}_\infty\simeq \operatorname{H}^1_{\mathrm{dR}}(M,\Kalg)\otimes_{\Kalg}\mathbb{C}_\infty
\end{equation}
under certain conditions; see \cite[Theorem 2.3.37]{Hartl2020} for details.
Choosing bases of $\operatorname{H}^1_{\mathrm{Betti}}(M,\mathbb{F}_q[t])$ and $\operatorname{H}^1_{\mathrm{dR}}(M,\Kalg)$, one defines the periods as the entries of the matrix representing the comparison isomorphism. 
Because of the $\mathbb{F}_q[t]$-structure and $\Kalg$-structure of the realizations, we can study transcendence properties of the periods over $\Kalg$. 
%\textcolor{red}{[Need to be confirm]}

In this paper, we define $v$-adic periods using de Rham and crystalline realizations. 
Since $\deg v=1$, we may equivalently work with Anderson $\mathbb{F}_q[t]$-motives or Anderson $\mathbb{F}_q[z]$-motives via the change of variable $z=t-\alpha$.
We use the latter notation for our $v$-adic construction.
Given an Anderson $\mathbb{F}_q[z]$-motive $M$ over $K$ satisfying certain conditions (see \Cref{Section:formulation_of_v-adic_periods} for details), we write $\overline{M}$ for the reduction of $M$ modulo $v$; then a natural comparison isomorphism
\begin{equation}
\operatorname{H}^1_{\mathrm{dR}}(M,\Kalg)\otimes_{\Kalg}\mathbb{C}_v\simeq \operatorname{H}^1_{\mathrm{cris}}(\overline{M},\Kcomp)\otimes_{\Kcomp}\mathbb{C}_v
\end{equation}
is known (\cite[Theorem 3.5.18]{Hartl2020a}).
However, in the standard formulation, the crystalline realization is a vector space over a local field, and, to the best of our knowledge, there exists no generally accepted choice of a $\Kalg$-structure on its scalar extension to $\mathbb{C}_v$ for use in defining periods;
by a $\Kalg$-structure, we mean a $\Kalg$-vector space $V$ together with a $\mathbb{C}_v$-linear
isomorphism
$
    V\otimes_{\Kalg}\mathbb{C}_v
    \xrightarrow{\sim}
    \operatorname{H}^1_{\mathrm{cris}}(\overline{M},\Kcomp)
    \otimes_{\Kcomp}\mathbb{C}_v.
$
Therefore, it is not straightforward to define $v$-adic periods and formulate the problem of their transcendence.
The notion of Andr\'e's $p$-adic periods resolves a similar issue in characteristic $0$; we adapt the approach of Ancona and Fr\u{a}\c{t}il\u{a} to formulate $v$-adic periods and study their transcendence.

\subsubsection{The idea coming from the work of Ancona and Fr\u{a}\c{t}il\u{a}}\label{subsec:intro_idea}

Following the idea of Ancona and Fr\u{a}\c{t}il\u{a} \cite{Ancona2025}, we define Andr\'e's $v$-adic periods for a certain class of Anderson $\mathbb{F}_q[z]$-motives. 
The guiding principle is that morphisms from the unit object to the object modulo $v$ play the role of the cycles in the formulation of $p$-adic periods in \cite[\S 2 (4)]{Ancona2025} and provide distinguished elements of the crystalline realization, thereby allowing us to define $v$-adic periods without requiring a canonical $\Kalg$-structure on the entire crystalline realization.

Roughly speaking, given an Anderson $\mathbb{F}_q[z]$-motive $M$ satisfying certain conditions (\Cref{Section:formulation_of_v-adic_periods}) with reduction modulo $v$ denoted by $\overline{M}$, we consider a morphism from the unit object to $\overline{M}$, in analogy with the formulation in \cite[\S 2 (4)]{Ancona2025}.
Such a morphism gives rise to an element of the crystalline realization of the original Anderson $\mathbb{F}_q[z]$-motive; via the de Rham--crystalline comparison isomorphism, this element can be expressed with respect to a basis coming from the natural $\Kalg$-structure on the de Rham realization of $M$.  We regard the resulting coefficients as Andr\'e's $v$-adic periods.

A central feature of the present formulation is that the construction naturally produces the Frobenius difference equations required for applying Chen's $v$-adic ABP criterion, together with suitable deformation series of the $v$-adic periods satisfying these equations.
This feature also places the resulting deformations of periods in a difference-Galois-theoretic framework.  
To establish and exploit this connection with the $v$-adic ABP criterion and difference Galois theory in a form suitable for transcendence arguments, we realize the idea above somewhat indirectly in the body of the paper: instead of defining the periods directly as the coordinates described above, we first construct the relevant deformation series and then obtain the periods by specialization.  
For the particular objects used in this paper to study special values of $v$-adic Carlitz multiple polylogarithms, this construction essentially recovers Chen's deformation series and the Frobenius difference equations satisfied by them.

An approach based on Frobenius actions on crystalline realizations has already appeared in the work of Chang, Wei, and Yu \cite[\S 3.7]{Chang2025}; we leave a general comparison with our formulation to future work.
We discuss a special case of this relationship in \cref{Rem:CWY}.
Mishiba (\cite{Mishiba2012}) also introduced a notion of $v$-adic periods of $t$-motives;
his formulation is based on comparison isomorphisms between \'{e}tale and Betti realizations of $t$-motives, while the work in \cite{Chang2025} and our formulation use comparison isomorphisms between crystalline and de Rham realizations.

\subsection{The main result}

The main result of the present paper is the following theorem relating transcendence properties of special values of $v$-adic Carlitz multiple polylogarithms to those of special values of $\infty$-adic Carlitz multiple polylogarithms.
Writing $\mathbb{N}$ for the set of positive integers, we define the subset $\tilde{\mathbb{D}}_\infty \subset \mathbb{N}\times \Kalg$ to be the set of pairs $(s,\,a)$ satisfying the condition $|a|_\infty<| \theta|_\infty^{\frac{sq}{q-1}}$ and define the subset $\tilde{\mathbb{D}}_v \subset \mathbb{N}\times \Aloc$ to be the set of pairs $(s,\,a)$ such that $|a|_v<1$.
For $\bullet \in \{ \infty, v\}$, we shall write
\begin{equation}\label{eq:polylogabb}
    \mathrm{Li}(\tilde{\mathbf{a}})_\bullet:=
    \begin{cases}
        \mathrm{Li}_{s_1,\,\dots,\,s_r}(a_1,\,\dots,\,a_r) & \bullet=\infty\\
        \mathrm{Li}_{s_1,\,\dots,\,s_r}(a_1,\,\dots,\,a_r)_v & \bullet=v
    \end{cases}
\end{equation}
for $\tilde{\mathbf{a}}=\big((s_1,\,a_1),\,\dots,\,(s_r,\,a_r)\big) \in \tilde{\mathbb{D}}_\bullet^r$ in what follows; we simplify the notation by this abuse of language.
\begin{Theorem}\label{Thm:main_CMPLs}
    Take a prime ideal $v\subset A$ of degree $1$ and a finite subset $I$ of $\bigcup_{r=1}^\infty(\tilde{\mathbb{D}}_\infty \cap \tilde{\mathbb{D}}_v)^r$. Assume that $\tilde{\mathbf{a}}=\big((s_1,\,a_1),\,\dots,\,(s_r,\,a_r)\big)\in I$ implies $\big((s_i,\,a_i),$ $\dots,$ $(s_j,\,a_j)\big)\in I$ for any $1 \leq i \leq j \leq r$.
    Then it follows that
    \begin{equation}\label{eq:mainthm}
        \trdeg_{\Kalg}\Kalg \left(\mathrm{Li}(\tilde{\mathbf{a}})_{v} \,\middle|\,\tilde{\mathbf{a}} \in I \right)=\trdeg_{\Kalg}\Kalg \left(\tilde{\pi},\,\mathrm{Li}(\tilde{\mathbf{a}})_{\infty}\,\middle|\,\tilde{\mathbf{a}} \in I \right)-1.
    \end{equation}
\end{Theorem}

For calculating the transcendence degree on the right-hand side of \cref{eq:mainthm}, powerful machinery is already available in \cite[Theorem 5.2.2]{Papanikolas2008}.
Hence, it seems that \cref{Thm:main_CMPLs} provides us with an approach to studying algebraic independence of special values of $v$-adic Carlitz multiple polylogarithms.
Chen studied linear independence of special values of $v$-adic multiple polylogarithms in \cite{Chen2025}; one aim of this paper is to develop a method for applying his idea to the problem of algebraic independence of those special values.
The results of Chang, Chen, and Mishiba in \cite{Chang2019a}, \cite{Chang2021}, and \cite{Chang2022} imply that the right-hand side of \cref{eq:mainthm} is greater than or equal to the left-hand side; see \Cref{SS:opposite_ineq}. 
The contribution of the present paper is the opposite inequality.

\subsection{Organization of the paper}

The paper is organized as follows. 
In \Cref{Section:infty-side}, we recall the part of Papanikolas' theory (\cite{Papanikolas2008}) of $\infty$-adic periods and their transcendence that will be used in the sequel, focusing on special values of Carlitz multiple polylogarithms.  
In particular, we review the relevant deformation series, their Frobenius difference equations, and the application of the ABP criterion.  
Consequently, we can relate the transcendence properties of special values of $\infty$-adic Carlitz multiple polylogarithms to those of their deformations.

In \Cref{Section:formulation_of_v-adic_periods}, we formulate Andr\'{e}'s $v$-adic periods for a class of Anderson $\mathbb{F}_q[z]$-motives with respect to a degree-one prime $v$.  
We first recall Anderson $\mathbb{F}_q[z]$-motives, following the formulation in \cite[Definition 2.3]{Gazda2024} and \cite[Definition 2.1]{Hartl2019}, as well as their reductions modulo $v$ and the de Rham--crystalline comparison isomorphism (\cite[Lemma 3.5.8 and Theorem 3.5.18]{Hartl2020a}) used in the construction.  
We then define the series in the variable $z$ called $v$-adic pre-periods; the $v$-adic periods are defined as their specializations at $z=\zeta$.  
We also give several examples of the construction.  
In particular, we show that special values of $v$-adic Carlitz multiple polylogarithms at points all of whose coordinates lie in $v\Aloc$ can be written in terms of such periods and observe that the pre-periods are essentially the same as the deformation series introduced by Chen in his study of linear independence of $v$-adic Carlitz multiple polylogarithms.

In \Cref{Section:v-adic_ABP}, we develop a method for studying the transcendence of $v$-adic periods using Chen's ABP criterion (\cite[Theorem 2.2.2]{Chen2025}).
We start by recalling Chen's $v$-adic ABP criterion in a form used in this paper. 
Following the argument of Papanikolas, we derive a comparison between the transcendence degrees of suitable deformation series satisfying Frobenius difference equations and those of their specializations.  
We then observe that the construction of $v$-adic periods naturally implies that the pre-periods satisfy Frobenius difference equations.  
Applying this result to the case of Carlitz multiple polylogarithms, we obtain a relation between the transcendence properties of the $v$-adic special values and those of their corresponding pre-periods.

In \Cref{sec:difference galois}, we develop a difference-Galois-theoretic method that serves as a bridge between the deformation series of $\infty$-adic periods and the corresponding $v$-adic pre-periods.  
We consider the difference domain generated by the deformation series of the $\infty$-adic periods and the difference field obtained by adjoining the corresponding $v$-adic pre-periods; the difference structures of these objects are given by the $(-1)$-st twists.
After replacing these twists by their iterates, we use standard arguments from difference Galois theory to construct a difference domain containing both objects as difference subrings.  
Within this larger difference domain, we directly compare the transcendence properties of the $\infty$-adic deformation series and the corresponding $v$-adic pre-periods.  
Applying this technique to Carlitz multiple polylogarithms, we obtain an inequality between the transcendence degree of the $v$-adic pre-periods and that of the $\infty$-adic deformation series.  

Finally, \Cref{Section:conclusiongeneral} is devoted to the proof of \cref{Thm:main_CMPLs}. 
We first recall the results of Chang, Chen, and Mishiba (\cite{Chang2022}, \cite{Chang2019a}, and \cite{Chang2021}) that give one of the inequalities needed for \cref{eq:mainthm}. 
We then obtain the reverse inequality by combining results in the previous sections, which completes the proof of \cref{Thm:main_CMPLs}.
We conclude the paper by stating some algebraic independence results for special values of $v$-adic Carlitz multiple polylogarithms that follow from the theorem and by discussing possible extensions of the present framework.

\begingroup
\renewcommand{\addcontentsline}[3]{}
\subsection*{Acknowledgements}

The author is grateful to Prof. Chang, Prof. Furusho, and Prof. Papanikolas for their mentorship, encouragement, and continued support.
He is also thankful to Prof. Mishiba for his support and valuable discussions.

This work was partially supported by the National Science and Technology Council of Taiwan and the NTHU Postdoctoral Researcher Mini-Project Program when the author was a postdoctoral researcher at National Tsing Hua University. 
He is also supported by JSPS Overseas Research Fellowships.

\par
\endgroup

\begingroup
\renewcommand{\addcontentsline}[3]{}
\subsection*{Declaration of AI Use}

During the preparation of the present paper, the author used AI tools (ChatGPT and Claude) to help check and improve the wording and organization of the manuscript, as well as some of its arguments.
In addition, he used the tools to scrutinize the proof of Proposition 3.1.3 in \cite{Anderson2004}, and the results of this examination contributed to the proof of \cref{Prop:ABP_condition}.
The author reviewed and verified all output, and he takes full responsibility for the final manuscript.

\par
\endgroup

\section{A brief review of the transcendence theory of special values of $\infty$-adic Carlitz multiple polylogarithms}\label{Section:infty-side}

In this section, we recall the part of the transcendence theory of special values of $\infty$-adic Carlitz multiple polylogarithms that will be used later, especially their deformation series and the associated Frobenius difference equations.
Our purpose is to explain how the transcendence properties of these special values can be studied via their deformation series.
Although these topics admit interpretations in the theory of $\infty$-adic periods, we do not use that formalism here.
Instead, we focus on the explicit constructions needed for our later arguments.
We mainly follow \cite{Chang2014,Papanikolas2008}.

Let us fix the notation. We write $A$ for the polynomial ring $\mathbb{F}_q[\theta]$ and $K$ for the field $\mathbb{F}_q(\theta)$ of fractions.
We also consider the completion $K_\infty:=\mathbb{F}_q((1/\theta))$ of $K$ with respect to the $\infty$-adic absolute value given by
\begin{equation}
    |a/b|_\infty:=q^{\deg a -\deg b} \quad (a,\,b\in A,\,b \neq0),
\end{equation}
and the completion of a fixed algebraic closure $K_\infty^{\mathrm{alg}}$ of $K_\infty$ is denoted by $\mathbb{C}_\infty$; we fix an inclusion $\Kalg \hookrightarrow \mathbb{C}_\infty$. 
By abuse of notation, we still denote by $|\cdot|_{\infty}$ the extended $\infty$-adic absolute value on $\mathbb{C}_\infty$. 
We let $t$ be a variable independent of $\theta$ and consider the field $\mathbb{C}_\infty((t))$ of Laurent series. 
For a Laurent series
\begin{equation}
    f=\sum_{i=m}^\infty a_i t^i \in \mathbb{C}_\infty((t))
\end{equation}
and $n\in \mathbb{Z}$, we define its $n$-th twist $f^{(n)}$ by
\begin{equation}  
f^{(n)}=\sum_{i=m}^\infty a_i^{q^n} t^i \in \mathbb{C}_\infty((t)).
\end{equation}
For a matrix $M$ with entries in $\mathbb{C}_\infty((t))$, we define its twist $M^{(n)}$ by applying the $n$-th twist entry-wise.
We denote by the symbol $\mathbb{T}_\infty$ the \textit{Tate algebra} over $\mathbb{C}_\infty$ defined as follows:
\begin{equation}
    \mathbb{T}_\infty:=\left\{ \sum_{i=0}^\infty a_i t^i \, \middle|\, a_i \in \mathbb{C}_\infty, \,|a_i|_{\infty}\to 0 \text{ as $i \rightarrow \infty$}\right\}.
\end{equation}
%The field of fractions of $\mathbb{T}$ is denoted by $\mathbb{L}$. 
A power series 
\begin{equation}
    f=\sum_{i=0}^\infty a_i t^i \in \Kalg[[t]]
\end{equation}
is said to be \textit{entire} if we have
\begin{equation}
    \lim_{i \to \infty}\sqrt[i]{|a_i|_\infty}=0 \text{ and } [K_\infty(a_1,\,a_2,\,\dots):K_\infty]<\infty
\end{equation}  following \cite{Anderson2004}.
An entire power series converges for any $t\in \mathbb{C}_\infty$, and we write $\Entire$ for the ring of entire power series. 
\begin{Example}[{\cite[\S 2.5]{Anderson1990}}]
    The \textit{Anderson-Thakur function} $\omega_\infty$ is defined as follows:
    \begin{equation}
       \omega_\infty=\omega_\infty(t)=(-\theta)^{\frac{1}{q-1}} \prod_{i=0}^\infty \left(1-\frac{t}{\theta^{q^i}}\right)^{-1}
    \end{equation}
    with a $(q-1)$-th root of $-\theta$ chosen.
    It holds that
    \begin{equation}\label{eq:omegaeq}
        \omega_\infty^{(1)}=(t-\theta)\omega_\infty.
    \end{equation}
    The function $\Omega_\infty(t):=(t-\theta)^{-1}\omega_\infty^{-1}=(\omega_\infty^{(1)})^{-1}$ satisfies 
    \begin{equation}\label{eq:Omega_inftyFrobeniuseq}
    \Omega_{\infty}^{(-1)}=(t-\theta)\Omega_\infty.
    \end{equation}
    It is known that the series $\Omega_\infty$ is an entire function (\cite[\S 3.3.2]{Anderson2004}) and that its value at $t=\theta$ is related to the Carlitz period as follows:
    \begin{equation}\label{eq:def_pitilde}
        \tilde{\pi}=-\frac{1}{\Omega_{\infty}}\bigg|_{t=\theta}=\theta(-\theta)^{\frac{1}{q-1}} \prod_{i=1}^\infty \left(1-\frac{\theta}{\theta^{q^i}}\right)^{-1}.
    \end{equation}
\end{Example}

\begin{Example}[{\cite[\S 5.3]{Chang2014}}]
    Writing 
    \begin{equation}\label{eq:def_lambda}
        \lambda_i := \left(t-\theta^{(0)}\right)\cdots\left(t-\theta^{(i)}\right)
    \end{equation}
    for $i \geq 0$, we define $\prelog_{\tilde{\mathbf{a}},\,m} \in \Kalg(t)$ to be 
    \begin{equation}\label{eq:deftruncatedprepolylog}
       \sum_{m\geq i_1>i_2>\cdots>i_r\geq0}\frac{a_1^{(i_1)}\cdots a_r^{(i_r)}}{\lambda_{i_1}^{s_1}\cdots\lambda_{i_r}^{s_r}}      
    \end{equation}    
 for $m \geq 0$ and $\widetilde{\mathbf{a}}=\big((s_1,\,a_1),\,\dots,\,(s_r,\,a_r)\big) \in (\mathbb{N}\times \Kalg)^r$; the sum is empty if $r>m+1$. 
Then, we put 
    \begin{equation}\label{eq:def infty prepolylog}
        \prelog_{\infty,\,\tilde{\mathbf{a}}}:=\lim_{m\to \infty} \prelog_{\tilde{\mathbf{a}},\,m}     \in \mathbb{C}_\infty((t))
    \end{equation}
    if $\tilde{\mathbf{a}} \in \tilde{\mathbb{D}}_\infty^r$ (\cite[Lemma 5.3.1]{Chang2014}).   
     The special values of these series at $t=\theta$ are related to the special values of $\infty$-adic Carlitz multiple polylogarithms as follows: 
        \begin{equation}\label{eq:inftyadicprepolylogpolylog}
            \prelog_{\infty,\,\tilde{\mathbf{a}}}^{(1)}\big|_{t=\theta}=\mathrm{Li}(\tilde{\mathbf{a}})_{\infty}-a_r\prelog_{\infty,\,\tilde{\mathbf{a}}^\prime}^{(1)}\big|_{t=\theta}.
        \end{equation}
        Here, we put $\tilde{\mathbf{a}}^\prime=\big((s_1,\,a_1),\,\dots,\,(s_{r-1},\,a_{r-1})\big)$ if $r \geq 2$ and put $\prelog_{\infty,\,\tilde{\mathbf{a}}^\prime}:=1$ if $r=1$.
        As the following lemma implies, these series are elements of $\Frac \Entire \subset \Frac \mathbb{T}_\infty$ after applying the $1$-fold twist.
\end{Example}

\begin{Lemma}[{\cite[Lemma 5.3.1]{Chang2014}}]\label{Lemma:prelog_in_Tate}
    For $\tilde{\mathbf{a}}=\big((s_1,\,a_1),\,\dots,\,(s_r,\,a_r)\big)\in \tilde{\mathbb{D}}_\infty^r$, the series $\Omega_\infty^{s_1+\cdots+s_r}\prelog_{\infty,\,\tilde{\mathbf{a}}}^{(1)}$ is entire.
\end{Lemma}
\begin{proof}
    Lemma 5.3.1 of \cite{Chang2014} proves that the series
    \begin{align}
        &\Omega_\infty^{s_1+\cdots+s_r}\prelog_{\infty,\,\tilde{\mathbf{a}}}^{(1)}+a_r\Omega_\infty^{s_1+\cdots+s_r}\prelog_{\infty,\,\tilde{\mathbf{a}}^\prime}^{(1)}\\
        &=\Omega_\infty^{s_1+\cdots+s_r}\sum_{i_1>i_2>\cdots>i_r\geq 0}\frac{{a_1^{(i_1)}}\cdots a_r^{(i_r)}}{\big(\left(t-\theta^{q}\right)\cdots\left(t-\theta^{q^{i_1}}\right)\big)^{s_1}\cdots\big(\left(t-\theta^{q}\right)\cdots\left(t-\theta^{q^{i_r}}\right)\big)^{s_r}} 
    \end{align}
    is entire, where we put $\tilde{\mathbf{a}}^\prime=\big((s_1,\,a_1),\,\dots,\,(s_{r-1},\,a_{r-1})\big)$. Therefore, we can show that the series $\Omega_\infty^{s_1+\cdots+s_r}\prelog_{\infty,\,\tilde{\mathbf{a}}}^{(1)}$ is entire by induction on $r$.
\end{proof}

The following theorem is a part of Chang's refinement (\cite[(1.1)]{Chang2009}) of Papanikolas' theorem (\cite[Theorem 5.2.2]{Papanikolas2008}). 
In its original form, the statement is formulated in the language of pre-$t$-motives,
and both its formulation and proof depend essentially on Tannakian duality and difference Galois theory.
Since our arguments only require the resulting transcendence statement in terms of matrices satisfying a certain difference equation, we do not recall those general aspects here.
Instead, we state only the concrete form needed later.

\begin{Theorem}[{\cite[(1.1)]{Chang2009}, \cite[Theorem 5.2.2]{Papanikolas2008}}]
\label{ThmPapanikolas}
Take matrices $\Phi$ in $\operatorname{Mat}_{ r}(\Kalg[t])\cap\operatorname{GL}_r(\Kalg(t))$ and $\Psi=(\Psi_{ij})$ in $\operatorname{Mat}_{ r}(\Entire)\cap\operatorname{GL}_r(\Frac\mathbb{T}_\infty)$ such that $\Psi^{(-1)}=\Phi \Psi$. 
If $\det \Phi |_{t=0}\neq 0$ and $\det \Phi |_{t=\theta^{1/q^{i}}}\neq 0$ for all $i \geq1$, then we have 
\begin{equation}\label{Papanikolasthmtrdeg}
    \operatorname{trdeg}_{\Kalg}\Kalg\left( \Psi_{ij} |_{t=\theta}\right)= \operatorname{trdeg}_{\Kalg(t)}\Kalg\left( t,\, \Psi_{ij} \right), 
\end{equation}
 where the indices $i$ and $j$ run through the range $1 \leq i,\,j \leq r$%(refer also to \cite[Theorem 3.2]{Mishiba2017})
 .
\end{Theorem}

We note that the proof of \cref{Papanikolasthmtrdeg} depends on the refined version (\cite[Theorem 1.2]{Chang2009}) of the ABP criterion (\cite[Theorem 1.3.2]{Anderson2004}). 
The following proposition is used to verify the hypotheses of \cref{ThmPapanikolas}.

\begin{Proposition}[{\cite[Proposition 3.1.3]{Anderson2004}}] \label{propentire}
  Given
        $ \Phi \in \operatorname{Mat}_{r \times r}(\Kalg[t]),
    $
    suppose that there exists  
    $\psi \in \operatorname{Mat}_{r\times 1}(\mathbb{T}_\infty)$
    such that $ \psi^{(-1)}=\Phi \psi.$
    If  $\det \Phi|_{t=0} \neq 0$, then all entries of $\psi$ are entire.
\end{Proposition}

%\begin{Example}\label{Example:inftypolylogperiodinterpretation0}
The following construction will be used in the proof of \cref{Cor:inftyprepolylogpolylog} and is almost the same as the one in \cite[\S 5.3]{Chang2014}, although the setting is slightly different.
    We take $\tilde{\mathbf{a}}=\big((s_1,\,a_1),\,\dots,\,(s_r,\,a_r)\big)\in \tilde{\mathbb{D}}_\infty^r$ and let $\Phi_{\tilde{\mathbf{a}}}\in \GL_{r+1}(\Kalg(t))$ be the matrix given by
    \begin{equation}\label{eq:def matrix T}
    \begin{pmatrix}
  1
  & 0 & \cdots & \cdots & 0\\%[25pt]
  \frac{a_1}{(t-\theta)^{s_1}} & \frac{1}{(t-\theta)^{s_1}} & \ddots & & \vdots\\%[25pt]
  0
  & \ddots
  & \frac{1}{(t-\theta)^{s_1+s_2}}
  &\ddots& \vdots\\%[25pt]
  \vdots
  &\ddots&
  \ddots
  & \ddots
  &0\\%[25pt]
  0 & \cdots & 0&  \frac{a_r}{(t-\theta)^{s_1+\cdots +s_r}}
  & \frac{1}{(t-\theta)^{s_1+\cdots +s_r}}
\end{pmatrix}.
\end{equation}
For $1 \leq i < j \leq r+1$, put
\begin{equation}
\prelog_{\infty,\,ji}:=\prelog_{\infty,\,\tilde{\mathbf{a}}^\prime},
\qquad
\tilde{\mathbf{a}}^\prime:=\big((s_i,\,a_i),\,\dots,\,(s_{j-1},\,a_{j-1})\big),
\end{equation}
where $\prelog_{\infty,\,\tilde{\mathbf{a}}^\prime}$ is defined by \cref{eq:def infty prepolylog}.
Then we define 
    \begin{equation}\label{eq:def matrix Psi}
       \begin{aligned}       
       &\Psi_{\tilde{\mathbf{a}}}:= \begin{pmatrix}
1 & 0 & \cdots & \cdots & 0\\
\prelog_{\infty,\,21}^{(1)} &\frac{1}{\Omega_{\infty}^{s_1}} & 0 & \cdots & 0\\
\vdots & \frac{\prelog_{\infty,\,32}^{(1)}}{\Omega_{\infty}^{s_1}} & \frac{1}{\Omega_{\infty}^{s_1+s_2}} & \ddots & \vdots\\
\prelog_{\infty,\,r1}^{(1)} & \vdots & \ddots & \ddots & 0\\
\prelog_{\infty,\,(r+1)1}^{(1)} &\frac{\prelog_{\infty,\,(r+1)2}^{(1)}}{\Omega_{\infty}^{s_1}} & \cdots & \frac{\prelog_{\infty,\,(r+1)r}^{(1)}}{\Omega_{\infty}^{s_1+\cdots+s_{r-1}}} & \frac{1}{\Omega_{\infty}^{s_1+\cdots+s_r}}
\end{pmatrix}\\
&=\begin{pmatrix}
1 & 0 & \cdots & \cdots & 0\\
\prelog_{\infty,\,21} &\omega_{\infty}^{s_1} & 0 & \cdots & 0\\
\vdots & \omega_{\infty}^{s_1}\prelog_{\infty,\,32} & \omega_{\infty}^{s_1+s_2} & \ddots & \vdots\\
\prelog_{\infty,\,r1} & \vdots & \ddots & \ddots & 0\\
\prelog_{\infty,\,(r+1)1} &\omega_{\infty}^{s_1}\prelog_{\infty,\,(r+1)2} & \cdots & \omega_{\infty}^{s_1+\cdots+s_{r-1}}\prelog_{\infty,\,(r+1)r} & \omega_{\infty}^{s_1+\cdots+s_r}
\end{pmatrix}^{(1)}.
\end{aligned}
    \end{equation}
    By \cref{Lemma:prelog_in_Tate}, we have
$
\Psi_{\tilde{\mathbf{a}}}\in \GL_{r+1}(\Frac\mathbb{T}_\infty)
$;
moreover, it holds that \begin{equation}\label{eq:inftypolylogmatrixdifferenceeq}
    \Psi_{\tilde{\mathbf{a}}}^{(-1)}=\Phi_{\tilde{\mathbf{a}}}\Psi_{\tilde{\mathbf{a}}},
    \end{equation}
     since we have \cref{eq:omegaeq} and
    \begin{equation}\label{eq:preloginftydifferenceeq}
        \prelog_{\infty,\,ji}=\frac{\prelog_{\infty,\,ji}^{(1)}}{(t-\theta)^{s_i+\cdots+s_{j-1}}}+\frac{a_{j-1}\prelog_{\infty,\,(j-1)i}^{(1)}}{(t-\theta)^{s_i+\cdots+s_{j-1}}},
    \end{equation}
    where we put $\prelog_{\infty,\,ii}=1$ for $1\leq i \leq r+1$ by convention.
    %\end{Example}
    Using \cref{ThmPapanikolas}, we obtain the following corollary:
\begin{Corollary}\label{Cor:inftyprepolylogpolylog}
    Take a finite subset $I$ of $\bigcup_{r=1}^\infty\tilde{\mathbb{D}}_\infty^r$ and assume that $\tilde{\mathbf{a}}=\big((s_1,\,a_1),$ $\dots,\,(s_r,\,a_r)\big)\in I$ implies $\big((s_i,\,a_i),$ $\dots,$ $(s_j,\,a_j)\big)\in I$ for any $1 \leq i \leq j \leq r$.
    Then it follows that
    \begin{equation}
        \trdeg_{\Kalg(t)}\Kalg \left(t,\,\Omega_\infty,\,\prelog_{\infty,\,\tilde{\mathbf{a}}}^{(1)}\,\middle|\,\tilde{\mathbf{a}} \in I \right)=\trdeg_{\Kalg}\Kalg \left(\tilde{\pi},\,\mathrm{Li}(\tilde{\mathbf{a}})_{\infty}\,\middle|\,\tilde{\mathbf{a}} \in I \right).
    \end{equation}
\end{Corollary}

\begin{proof}
    The matrix 
    \begin{equation}
    \Phi:=(t-\theta)\oplus \bigoplus_{\tilde{\mathbf{a}}\in I}\Phi_{\tilde{\mathbf{a}}}
    \end{equation}
    satisfies $\Psi^{(-1)}=\Phi\Psi$ with
    \begin{equation}
        \Psi:=(\Omega_{\infty})\oplus \bigoplus_{\tilde{\mathbf{a}}\in I}\Psi_{\tilde{\mathbf{a}}},
    \end{equation}
    by Equations \eqref{eq:Omega_inftyFrobeniuseq} and \eqref{eq:inftypolylogmatrixdifferenceeq}. 
    For sufficiently large $w$, we have $(t-\theta)^w\Phi \in \Mat(\Kalg[t])$, $\det \big( (t-\theta)^w\Phi\big)|_{t=0}\neq0$, $\det \big((t-\theta)^w\Phi\big)|_{t=\theta^{1/q^i}}\neq0$ for any $i \geq 1$, $\Omega_{\infty}^w \Psi \in  \Mat(\mathbb{T}_\infty)$, and
    \begin{equation}\label{eq:normalized_Frobenius_eq}
    \left(\Omega_\infty^w\Psi\right)^{(-1)}=(t-\theta)^w\Phi\left(\Omega_\infty^w\Psi\right)
    \end{equation}
     by \cref{eq:Omega_inftyFrobeniuseq} and the definitions of the matrices $\Phi_{\tilde{\mathbf{a}}}$ and $\Psi_{\tilde{\mathbf{a}}}$; then \cref{propentire} shows that  $\Omega_{\infty}^w \Psi \in  \Mat(\Entire)$.
    Now we can apply \cref{ThmPapanikolas} to \cref{eq:normalized_Frobenius_eq} and obtain 
    \begin{equation}
        \trdeg_{\Kalg(t)}\Kalg(t) \big(\text{entries of } \Omega_{\infty}^w \Psi\big)=\trdeg_{\Kalg}\Kalg \big(\text{entries of } (\Omega_{\infty}^w \Psi)|_{t=\theta} \big),
    \end{equation}
    which, together with \cref{eq:def_pitilde,eq:inftyadicprepolylogpolylog}, yields the desired result.
\end{proof}

\section{Formulation and examples of $v$-adic periods for degree one $v$}\label{Section:formulation_of_v-adic_periods}

This section is devoted to the formulation of Andr\'e's $v$-adic periods.
We shall recall the notion of Anderson $\mathbb{F}_q[z]$-motives in \Cref{subsection_Anderson_motives}; then we review a comparison theorem (\cref{Thm:comparison}), which appears in \Cref{subsec:comparisonthm} and plays a pivotal role in the formulation of Andr\'e's $v$-adic periods.
In \Cref{subsec:v-adic_periods}, we formulate Andr\'e's $v$-adic periods for a certain class of Anderson $\mathbb{F}_q[z]$-motives.

We continue to write $A$ for the polynomial ring $\mathbb{F}_q[\theta]$ and $K$ for the fraction field $\mathbb{F}_q(\theta)$.
Let us take a prime ideal $v$ of $A$, write $\zeta \in A$ for the monic generator of $v$, and put $z:=\zeta|_{\theta=t} \in \mathbb{F}_q[t]$. Let $\Aloc$ denote the localization of $A$ at $v$. 
The symbols $\Acomp$ and $\Kcomp$ stand for the $v$-adic completions of $\Aloc$ and $K$, respectively. 
Further, we consider the completion $\mathbb{C}_v$ of the algebraic closure of $\Kcomp$ and fix an inclusion $\Kalg \hookrightarrow \mathbb{C}_v$. 
We write $\quotfield$ for the residue field $A/v=\Aloc/v=\Acomp/v$.

Regarding $z$ as an independent variable over $\mathbb{C}_v$, we consider the ring 
\begin{equation}
    \Acomp[[z,\,1/z\}:=\left\{ \sum_{i=-\infty}^\infty a_i z^{i}\, \middle| a_i\in \Acomp,\,\, \lim_{i \to -\infty}|a_{i}|_v|\zeta^{i}|_v^r = 0 \text{ for any } r>0 \right\}
\end{equation}
following \cite[\S 2.5]{Hartl2009} and \cite[(3.5.2)]{Hartl2020a}.
For $r>0$, we consider the Gauss norm given by
\begin{equation}
    \Bigg|\Bigg|\sum_{i=-\infty}^\infty a_i z^{i}\Bigg|\Bigg|_{r}:=\max_{i \in \mathbb{Z}}|a_i|_v |\zeta|_v^{ri}.
\end{equation}
Then the ring $\Acomp[[z,\,1/z\}$ is equipped with the topology given by the family of these norms;
explicitly, a fundamental system of neighborhoods of zero
is given by the sets
\begin{equation}
    \left\{ f \in \Acomp[[z,\,1/z\}
    \,\middle|\,
    ||f||_{r_i}<\epsilon_i \text{ for } 1\leq i \leq m
    \right\}
\end{equation}
with $m \geq 1$, $r_i > 0$ and $\epsilon_i> 0$.
%\textcolor{red}{[Check]}.
We also note that the inversion map on
$(\Acomp[[z,\,1/z\})^\times$, endowed with the subspace topology,
is continuous; indeed, for a given $r>0$ and invertible elements $f$ and $g$ in $\Acomp[[z,\,1/z\}$ such that $||f-g||_r<||f||_r$, the ultrametric inequality gives $||g||_r=||f||_r$ and hence we have
\begin{equation}
 \left|\middle|\frac{1}{g}-\frac{1}{f}\middle|\right|_r
 =\left|\middle|\frac{f-g}{fg}\middle|\right|_r
 = \frac{||f-g||_r}{||f||_r||g||_r}=\frac{||f-g||_r}{||f||_r^2}
\end{equation}
as the norm is multiplicative.
The symbol $\mathbb{T}_v $ stands for the Tate algebra in the variable $1/z$ given as follows:
\begin{equation}
    \mathbb{T}_v:=\left\{ \sum_{i=0}^\infty a_i z^{-i} \, \middle| \, a_i \in \mathbb{C}_v,\,|a_i|_v \to 0 \text{ for } i \to \infty \right\},
\end{equation}
and we have $\Acomp[[1/z]]\cap\Acomp[[z,\,1/z\}\subset \mathbb{T}_v$.
We say that a power series 
\begin{equation}
    f=\sum_{i=0}^\infty a_i z^{-i} \in \Kalg[[1/z]]
\end{equation}
is \textit{entire} if we have
\begin{equation}
    \lim_{i \to \infty}\sqrt[i]{|a_i|_v}=0 \text{ and } [\Kcomp(a_1,\,a_2,\,\dots):\Kcomp]<\infty
\end{equation}  following \cite[\S 2.4]{Anderson2004} and \cite[\S 2.1]{Chen2025}.
An entire power series converges for any $1/z \in \mathbb{C}_v$, and we write $\Entirev$ for the ring of entire power series.
\begin{comment}
We consider the ring $\mathcal{R}_v$ defined to be
\begin{equation}
    \left\{
        f=\sum_{i\in\mathbb{Z}}b_iz^{i}
        \,\middle|\,
        \begin{array}{l}
            b_i\in\mathbb{C}_v,\ \text{and there exists $r$ with $0<r<1$, depending on $f$,}\\
            \text{such that $|b_i|_v\rho^{i}\to0$ as $i\to\pm\infty$ for every $\rho$ with $r<\rho<1$}
        \end{array}
    \right\},
\end{equation}
which contains $\Acomp[[z,$ $1/z\}$ as a subring; see \cite[Section 1.1]{Hartl2011} and \cite[Definitions 8.4.2 and 15.14]{Kedlaya2022}.
\end{comment}
For a formal bidirectional power series
\begin{equation}
    f=\sum_{i \in \mathbb{Z}}a_i z^{i},\quad a_i \in \mathbb{C}_v
\end{equation}
and $n\in \mathbb{Z}$, we define its $n$-th twist $f^{(n)}$ by
\begin{equation}  \label{eq:def_v_twisting}
f^{(n)}=\sum_{i \in \mathbb{Z}}a_i^{q^n} z^{i}.
\end{equation}
For every $n\in\mathbb{Z}$, the $n$-th twist induces a ring automorphism of $\mathbb{T}_v$. For every $n\geq0$, it also induces a ring endomorphism of $\Acomp[[z,\,1/z\}$.
Given a matrix $M$ whose entries are such formal bidirectional power series, we define its twist $M^{(n)}$ by applying the $n$-th twist entry-wise.

\begin{Example}[{\cite[Definition 3.1.1]{Chen2025}}]
If we put
\begin{equation}\label{eq:def_omegav}
    \omegav(z):=\prod_{i=0}^{\infty}
    \left(
    1-\frac{\zeta^{q^i}}{z}
    \right)^{-1},
\end{equation}
then this series satisfies the equation 
$
    \omegav^{(1)}=\left(1-\zeta/z\right)\omegav.
$
We can show that this series lies in $\mathbb{T}_v$ using the expansion
\begin{equation}
    \omegav(z):=\prod_{i=0}^{\infty}\sum_{j=0}^\infty \left(\frac{\zeta^{q^i}}{z}
    \right)^j
    =\sum_{j=0}^\infty 
    \left(
    \sum_{\substack{i_0,\,i_1,\dots \geq 0\\i_0+i_1+\cdots=j}}\zeta^{i_0}\zeta^{qi_1}\cdots
    \right)z^{-j}.
\end{equation}
The series
\begin{equation}\label{eq:def_Omegav}
    \Omegav:=\bigl(\omegav^{(1)}\bigr)^{-1}=\prod_{i=1}^{\infty}
    \left(
    1-\frac{\zeta^{q^i}}{z}
    \right)
\end{equation}
satisfies $\Omegav^{(-1)}=\left(1-\zeta/z\right)\Omegav$ and lies in $\Entirev$ (\cite[Proposition 3.1.3]{Chen2025}).
\end{Example}

We note that the series $\omegav^{-1}$ is in the ring $\Acomp[[z,\,1/z\}$ and maps to $1$ via the coefficient-wise modulo $v$ reduction $\Acomp[[z,\,1/z\} \rightarrow \quotfield((z))$. Thus, the reduction can be extended to the localization $\Acomp[[z,\,1/z\}[\omegav]$.

For the remainder of this section, we assume that the degree of the finite place $v$ is $1$; then there exists $\alpha \in \mathbb{F}_q$ such that $\zeta=\theta-\alpha$ and $z=t-\alpha$; hence we shall treat $z$ as a variable obtained by a change of variable from $t$ and have
\begin{equation}\label{eq:deg1consequence}
    z-\zeta^{q^i}=(t-\alpha)-(\theta-\alpha)^{q^i}=t-\alpha-\theta^{q^i}+\alpha^{q^i}
    =t-\theta^{q^i} \qquad (i \geq 0)
\end{equation}
as $\alpha^q=\alpha$.

\subsection{Anderson motives}\label{subsection_Anderson_motives}

In this subsection, we shall recall the notion of an Anderson $\mathbb{F}_q[z]$-motive, following \cite[Definition 2.3]{Gazda2024} and \cite[Definition 2.1]{Hartl2019}.
Readers familiar with these objects may skip this subsection, noting only our notation $\mathcal{M}$ and $\overline{\mathcal{M}}$, together with the conventions for defining matrices.

Given an $\mathbb{F}_q$-algebra $R$ together with an algebra homomorphism $\kappa:\mathbb{F}_q[z] \rightarrow R$, we define an endomorphism $\tau$ of $R[z]=R \otimes_{\mathbb{F}_q}\mathbb{F}_q[z]$ by $a \otimes b \mapsto a^q \otimes b$ for all $a \in R$ and $b \in \mathbb{F}_q[z]$. 
Here, we distinguish between the variable $z$ and its image $\kappa(z)$ in $R$. 
When $R$ is a subring of $\mathbb{C}_v$, the action of $\tau$ coincides with the first twist, where polynomials in $z$ are viewed as finite power series in $z$.
Given a module $M$ over $R[z]$, we write $\tau^* M$ for the scalar extension $M\otimes_{R[z],\,\tau} R[z]$ via the endomorphism $\tau$; then one has
\begin{equation}
fm\otimes g = m\otimes \tau(f)g 
\end{equation}
for $m\in M$ and $f,g\in R[z]$. 
An $R[z]$-linear mapping $u:M_1\rightarrow M_2$ gives rise to an $R[z]$-linear mapping $\tau^* u:\tau^* M_1 \rightarrow \tau^* M_2$ given by $\tau^*u(m\otimes g)=u(m)\otimes g$.

\begin{definition}[{\cite[Definition 2.3]{Gazda2024} and \cite[Definition 2.1]{Hartl2019}}] \label{Def:t-motive}
       For an $\mathbb{F}_q$-algebra $R$ together with an algebra homomorphism $\kappa:\mathbb{F}_q[z] \rightarrow R$, an \textit{Anderson $\mathbb{F}_q[z]$-motive over $R$} is a pair $M=(M,\,\tau_M)$ consisting of a projective $R[z]$-module $M$ of finite constant rank and an $R[z,\,1/(z-\kappa(z))]$-linear isomorphism
    \begin{equation}
        \tau_M: \tau^* M \otimes_{R[z]} R[z,\,1/(z-\kappa(z))] \rightarrow M\otimes_{R[z]} R[z,\,1/(z-\kappa(z))].
    \end{equation}
    We call the isomorphism $\tau_M$ the \textit{structure isomorphism} of $M$.
    A morphism of Anderson $\mathbb{F}_q[z]$-motives from $M_1$ to $M_2$ is an $R[z]$-linear mapping $u:M_1 \rightarrow M_2$ which is compatible with the structure isomorphisms; that is, we require the following diagram to commute:
    \begin{equation}
        \begin{tikzpicture}[auto] 
                   \node (22) at (-1, 2) {$\tau^* M_1 \otimes_{R[z]} R[z,\,1/(z-\kappa(z))]$}; \node (42) at (5, 2) {$M_1\otimes_{R[z]} R[z,\,1/(z-\kappa(z))]$}; 
                    \node (20) at (-1, 0) {$\tau^* M_2 \otimes_{R[z]} R[z,\,1/(z-\kappa(z))]$}; \node (40) at (5, 0) {$M_2\otimes_{R[z]} R[z,\,1/(z-\kappa(z))]$}; 
                    \draw[->] (42) to node {$u \otimes_{R[z]} R[z,\,1/(z-\kappa(z))] $}(40);
                    \draw[->] (22) to node {$\tau^*u \otimes_{R[z]} R[z,\,1/(z-\kappa(z))] $}(20);
                    
                    \draw[->] (22) to node {$\tau_{M_1}$}(42);
                    \draw[->] (20) to node {$\tau_{M_2}$}(40);

                \end{tikzpicture}.
    \end{equation}
    
\end{definition}

These are the objects called $A$-motives in \cite[Definition 2.3]{Gazda2024} and \cite[Definition 2.1]{Hartl2019}, with $A$ taken to be $\mathbb{F}_q[z]$ there, whereas in the present paper the symbol $A$ is reserved for the polynomial ring $\mathbb{F}_q[\theta]$. They are variants of the $t$-motives originally introduced by Anderson (\cite[\S 1.2]{Anderson1986}).

We define $\mathcal{M}_R=\mathcal{M}_{R,\kappa}$ to be the full subcategory consisting of Anderson $\mathbb{F}_q[z]$-motives that are free over $R[z]$ and, by abuse of notation, write $M \in \mathcal{M}_R$ when $M$ is an object of $\mathcal{M}_R$.
Further, we simply write $\mathcal{M}$ for $\mathcal{M}_{\Aloc,\,\kappa}$, where the homomorphism $\kappa:\mathbb{F}_q[z] \rightarrow \Aloc$ is the $\mathbb{F}_q$-algebra homomorphism given by $z\mapsto \zeta$, and write $\overline{\mathcal{M}}$ for $\mathcal{M}_{\quotfield,\,\overline{\kappa}}$.
Here, the homomorphism $\overline{\kappa}$ is obtained by composing $\kappa$ given above with the quotient by $v$.
We note that the action of $\tau$ on $\quotfield[z]$ is trivial since $\quotfield\simeq \mathbb{F}_q$ by the assumption $\deg v=1$.

For $M=(M,\tau_M)$ in $\mathcal{M}_R$, there exists an $R[z]$-basis
$
\mathbf{m}=(m_1,\dots,m_r)
$
and a matrix $\Phi_M\in \GL_r(R[z,\,1/(z-\kappa(z))])$ such that
$
\mathbf{m}\otimes 1:=(m_1\otimes 1,\dots,m_r\otimes 1)
$
satisfies
\begin{equation}
    \tau_M(\mathbf{m}\otimes 1)=\mathbf{m}\Phi_M.
\end{equation}
In this situation, we say that $\Phi_M$ is the defining matrix of $M$ with respect to the basis $\mathbf{m}$.

\begin{Example}
    The \textit{unit object} $\mathbf{1}$ of $\mathcal{M}$ is the rank-$1$ $\Aloc[z]$-module $\mathbf{1}$ with $\tau_{\mathbf{1}}$ given by $m\otimes 1 \mapsto m$ for a chosen generator $m$ of $\mathbf{1}$. In other words, $\mathbf{1} \in \mathcal{M}$ is the object defined by the identity matrix $I_1=(1)\in \GL_1(\Aloc[z,\,1/(z-\zeta)])$.
   Its reduction in $\overline{\mathcal{M}}$ is denoted by $\overline{\mathbf{1}}$ and this is also defined by the identity matrix $I_1\in \GL_1(\quotfield[z,\,1/z])$.
\end{Example}

\begin{Example}[{cf. \cite[Example 2.3.5]{Hartl2020}}]\label{Example:Carlitz}
    When the degree of $v$ is $1$, the \textit{Carlitz motive} $\mathcal{C}$ is the rank-$1$ $\Aloc[z]$-module with $\tau_{\mathcal{C}}$ given by $m\otimes 1 \mapsto (z-\zeta)m$ for a generating element $m\in \mathcal{C}$. That is, $\mathcal{C}$ is defined by the matrix $(z-\zeta)\in \GL_1(\Aloc[z,\,1/(z-\zeta)])$ of size $1$.
\end{Example}

Let us recall that $v\subset A=\mathbb{F}_q[\theta]$ is a degree one prime ideal with monic generator $\zeta$, and that we have put
$z:=\zeta|_{\theta=t}\in \mathbb{F}_q[t]$.
Reducing modulo $v=(\zeta)$, we obtain a surjective homomorphism
\begin{equation}
\Aloc[z,\,1/(z-\zeta)] 
%\rightarrow\Acomp[z,\,1/(z-\zeta)] 
\twoheadrightarrow \quotfield[z,\,1/z].
\end{equation}
Here, and in the rest of this paper, we write $\overline{a}$ for the image of
\begin{equation}a\in \Aloc[z,\,1/(z-\zeta)] 
\end{equation}
under the reduction map.
The reduction modulo $v$ can be restricted to
\begin{equation}
    \Aloc[z] \rightarrow \quotfield[z]
\end{equation}
and this is compatible with the aforementioned actions of $\tau$.

For any $M=(M,\,\tau_M) \in \mathcal{M} $, we can obtain its modulo $v$ reduction $\overline{M}=(\overline{M},\,\tau_{\overline{M}}) \in \overline{\mathcal{M}}$ as follows. The underlying module $\overline{M}$ is defined to be $M\otimes_{\Aloc[z]}\quotfield[z]$.
Let us consider the following diagram:
\begin{equation}
                \begin{tikzpicture}[auto] 
                   \node (22) at (-1, 2) {$\Aloc[z]$}; \node (42) at (3.5, 2) {$\Aloc[z]$}; \node (62) at (8, 2) {$\Aloc[z,\,1/(z-\zeta)]$};
                    \node (20) at (-1, 0) {$\quotfield[z]$}; \node (40) at (3.5, 0) {$\quotfield[z]$}; \node (60) at (8, 0) {$\quotfield[z,\,1/z]$};
    
                    \draw[->>] (62) to node[sloped] {}(60);
                    \draw[->>] (42) to node {}(40);
                    \draw[->>] (22) to node {}(20);
                    
                    \draw[->] (22) to node {$\tau$}(42);
                    \draw[->] (20) to node {$\tau$}(40);
    
                    \draw[->] (42) to node {}(62);
                    \draw[->] (40) to node {}(60);
                \end{tikzpicture}
            \end{equation}
The commutativity of the left square follows from the fact that the reduction $\Aloc[z]\twoheadrightarrow \quotfield [z]$ commutes with the actions of $\tau$; thus the whole diagram is commutative.
This commutative diagram yields natural isomorphisms
\begin{align}
&M\otimes_{\Aloc[z]}
\Aloc[z,1/(z-\zeta)]
\otimes_{\Aloc[z,1/(z-\zeta)]}
\quotfield[z,1/z]\\
&\simeq
M\otimes_{\Aloc[z]}\quotfield[z]
\otimes_{\quotfield[z]}
\quotfield[z,1/z]=
\overline M\otimes_{\quotfield[z]}\quotfield[z,1/z]
\end{align}
and
\begin{align}
    &\quad \ \
    \tau^* M \otimes_{\Aloc[z]}\Aloc[z,\,1/(z-\zeta)]\otimes_{\Aloc[z,\,1/(z-\zeta)]}\quotfield[z,\,1/z]\\
    &=M\otimes_{\Aloc[z],\tau}\Aloc[z]\otimes_{\Aloc[z]}\Aloc[z,\,1/(z-\zeta)]\otimes_{\Aloc[z,\,1/(z-\zeta)]}\quotfield[z,\,1/z]\\
   % &\simeq M\otimes_{\Aloc[z],\tau}\Aloc[z] \otimes_{\Aloc[z]}\quotfield[z,\,1/z]\\
    %&\simeq M\otimes_{\Aloc[z],\tau}\Aloc[z]\otimes_{\Aloc[z]}\quotfield[z]\otimes_{\quotfield[z]}\quotfield[z,\,1/z]\\
    &\simeq M\otimes_{\Aloc[z]}\quotfield[z]\otimes_{\quotfield[z],\tau}\quotfield[z] \otimes_{\quotfield[z]}\quotfield[z,\,1/z]\\
    &=\overline{M}\otimes_{\quotfield[z],\tau}\quotfield[z] \otimes_{\quotfield[z]}\quotfield[z,\,1/z]=\tau^*\overline{M}\otimes_{\quotfield[z]}\quotfield[z,\,1/z].
\end{align}
Therefore, we can deduce an isomorphism \begin{equation}\tau_{\overline{M}}:\tau^* \overline{M}\otimes_{\quotfield[z]}\quotfield[z,\,1/z]\rightarrow\overline{M}\otimes_{\quotfield[z]}\quotfield[z,\,1/z]
\end{equation}
from the isomorphism 
\begin{equation}
\tau_M:
 \tau^* M \otimes_{\Aloc[z]}\Aloc[z,\,1/(z-\zeta)]
\rightarrow 
M \otimes_{\Aloc[z]}\Aloc[z,\,1/(z-\zeta)]
\end{equation}
by tensoring with $\quotfield[z,\,1/z]$ over $\Aloc[z,\,1/(z-\zeta)]$.

Given an $\mathbb{F}_q$-algebra $R$ with a homomorphism $\kappa:\mathbb{F}_q[z] \rightarrow R$, and objects $M_1$ and $M_2$ in $\mathcal{M}_R$, we can consider the tensor product $M_1 \otimes M_2 =(M_1 \otimes_{R[z]} M_2,\,\tau_{M_1 \otimes M_2})$, where the linear isomorphism
\begin{align}
    \tau_{M_1 \otimes M_2}:&\tau^* ( M_1 \otimes_{R[z]}M_2)\otimes_{R[z]} R[z,\,1/(z-\kappa(z))] \\
    \ &\rightarrow ( M_1 \otimes_{R[z]}M_2)\otimes_{R[z]} R[z,\,1/(z-\kappa(z))]
\end{align}
is given by the following composition of the natural isomorphisms
\begin{align}
&\quad \ \tau^* ( M_1 \otimes_{R[z]}M_2)\otimes_{R[z]} R\left[z,\,\frac{1}{z-\kappa(z)}\right]\\
&=
\left( (M_1 \otimes_{R[z]} M_2)\otimes_{R[z],\,\tau}R[z]\right)\otimes_{R[z]} R\left[z,\,\frac{1}{z-\kappa(z)}\right]\\
&\simeq  \left( (M_1\otimes_{R[z],\,\tau}R[z]) \otimes_{R[z]} (M_2\otimes_{R[z],\,\tau}R[z])\right) \otimes_{R[z]} R\left[z,\,\frac{1}{z-\kappa(z)}\right]\\
&=\left( \tau^* M_1 \otimes_{R[z]} \tau^* M_2\right) \otimes_{R[z]} R\left[z,\,\frac{1}{z-\kappa(z)}\right]\\
&=\left( \tau^* M_1\otimes_{R[z]} R\left[z,\,\frac{1}{z-\kappa(z)}\right]\right) \otimes_{R[z,\,1/(z-\kappa(z))]} \left( \tau^* M_2\otimes_{R[z]} R\left[z,\,\frac{1}{z-\kappa(z)}\right] \right),
\end{align}
the tensor product $\tau_{M_1}\otimes \tau_{M_2}$ over $R[z,\,1/(z-\kappa(z))]$ of linear isomorphisms, and the natural isomorphism
\begin{align}
   &\quad \ \ ( M_1 \otimes_{R[z]}R[z,\,1/(z-\kappa(z))]) \otimes_{R[z,\,1/(z-\kappa(z))]} ( M_2 \otimes_{R[z]}R[z,\,1/(z-\kappa(z))])\\& \simeq ( M_1 \otimes_{R[z]}M_2)\otimes_{R[z]} R[z,\,1/(z-\kappa(z))].
\end{align}
In order to describe the matrices defining tensor products of objects of $\mathcal{M}_R$, we recall the notion of the \textit{Kronecker product} of matrices.
For any matrices $X=(X_{ij})$ and $Y=(Y_{i^\prime j^\prime})$, we define their Kronecker product $X\otimes Y$ to be the block matrix $(X_{ij} Y)$ if the products $X_{ij} Y_{i^\prime j^\prime}$ of entries are defined. 
We shall use the following standard properties of Kronecker products.
Given matrices $X$, $Y$, $X^\prime$, and $Y^\prime$ of suitable sizes, we have
\begin{equation}\label{eq:Kronecker_mixed}
    (X\otimes Y)(X^\prime\otimes Y^\prime)=(XX^\prime)\otimes(YY^\prime)
\end{equation}
and the associativity
\begin{equation}\label{eq:Kronecker_associativity}
    (X\otimes Y)\otimes Z=X\otimes (Y\otimes Z).
\end{equation}
Further, if $X$ and $Y$ are invertible square matrices of sizes $r$ and $s$, respectively, it holds that
\begin{equation}\label{eq:Kronecker_det}
    \det (X\otimes Y)=(\det X)^{s}(\det Y)^{r}.
\end{equation}
As the twists are applied entry-wise to matrices, we also have
\begin{equation}\label{eq:Kronecker_twist}
    (X\otimes Y)^{(n)}=X^{(n)}\otimes Y^{(n)}
\end{equation}
for $n \in \mathbb{Z}$ whenever the twisting is defined for the entries of $X$ and $Y$.
For $i=1,\,2$, let $M_i$ be an object of $\mathcal{M}_R$ which has the defining matrix $\Phi_{M_i} \in \GL_{r_i}(R[z,\,1/(z-\kappa(z))])$ with respect to an $R[z]$-basis $\mathbf{m}_i=(m_{i1}, \, \dots,$ $m_{ir_i})$. 
Then the defining matrix of the tensor product $M_1\otimes M_2$ is the Kronecker product $\Phi_{M_1}\otimes \Phi_{M_2} \in \GL_{r_1r_2}(R[z,\,1/(z-\kappa(z))])$ if we choose the basis $\mathbf{m}_1\otimes\mathbf{m}_2$ given by
\begin{align}
(m_{11}\otimes& m_{21}, \, m_{11}\otimes m_{22},\,\dots,\,m_{11}\otimes m_{2r_2},\,
m_{12}\otimes m_{21}, \, m_{12}\otimes m_{22},\,\dots\\&\dots,\,m_{12}\otimes m_{2r_2},
\,\cdots,\,
m_{1r_1}\otimes m_{21}, \, m_{1r_1}\otimes m_{22},\,\dots,\,m_{1r_1}\otimes m_{2r_2}).
\end{align}

\subsection{The comparison isomorphism}\label{subsec:comparisonthm}
We shall review \cref{Thm:comparison} and describe the isomorphism given in the theorem explicitly in several examples.
Before stating the comparison theorem, we explain how the defining matrices of an object $M\in\mathcal{M}$ and its reduction $\overline{M}$ can be viewed as matrices over a common ring.

For each $i\geq0$, the element $z-\zeta^{q^i}$ is invertible in $\Acomp[[z,1/z\}[\omegav]$.  Indeed, we have
\begin{align}
& \frac{1}{z-\zeta^{q^i}}\omegav^{-1}
=
\frac{1}{z}\left(1-\frac{\zeta^{q^i}}{z}\right)^{-1} \prod_{j=0}^{\infty}
    \left(
    1-\frac{\zeta^{q^j}}{z}
    \right)
=\frac{1}{z}\prod_{\substack{0\leq j < \infty\\ j \neq i}}
    \left(
    1-\frac{\zeta^{q^j}}{z}
    \right)
\\
&=\sum_{n=0}^{\infty}\left(
\sum_{\substack{j_1<\dots<j_n\\j_1\neq i,\cdots,\,j_n \neq i}}
\left(
-\zeta^{q^{j_1}}\right)\cdots \left(
-\zeta^{q^{j_n}}\right)
\right)z^{-n-1}
\in  \Acomp[[z,1/z\},
\end{align}
see \cite[\S 3.5]{Hartl2020a}.
Hence, there is a natural injective ring homomorphism
\begin{equation}
\Aloc[z,\,1/(z-\zeta)]
\longrightarrow
\Acomp[[z,1/z\}[\omegav]=\Acomp[[z,1/z\}\left[1/\Omegav^{(-1)}\right].
\end{equation}
The localization $\Acomp[[z,1/z\}\left[1/\Omegav^{(-1)}\right]$ is regarded as a topological localization of the topological ring $\Acomp[[z,1/z\}$ (\cite[Definition 2.1]{Requejo1994}); that is, if a continuous homomorphism $f:\Acomp[[z,1/z\} \rightarrow B$ to a topological ring $B$ is given and $f\left( \Omegav^{(-1)} \right)$ is invertible in $B$, then $f$ factors uniquely through $\Acomp[[z,1/z\}\left[1/\Omegav^{(-1)}\right]$ via a continuous homomorphism.

We next embed the ring $\quotfield[z,1/z]$ into the same ambient ring.
We regard $\quotfield$ as a subfield of $\Acomp$ via the Teichm\"uller section
$
\quotfield \longrightarrow \Acomp
$
of the reduction map $\Acomp \twoheadrightarrow \quotfield$.
Since $\deg v=1$, this section is given explicitly by
\begin{equation}
\overline{a}\longmapsto \lim_{m\to\infty}a^{q^m},
\end{equation}
where $a\in\Aloc$ is any preimage of $\overline{a}\in\quotfield$ under the
reduction $\Aloc\to\quotfield$.
Applying the Teichm\"uller section coefficientwise and keeping the variable
$z$ fixed, we obtain an embedding
\begin{equation}
\quotfield[z,1/z]\hookrightarrow \Acomp[z,\,1/z]
\subset \Acomp[[z,1/z\}[\omegav].
\end{equation}
Considering the reduction $\Acomp[[z,1/z\}[\omegav]\rightarrow \quotfield((z))$ modulo $v$, we observe that $\omegav$ is mapped to $1$ and that the reduction induces the identity mapping on $\quotfield[z,1/z]  \subset \Acomp[[z,1/z\}[\omegav]$.
In what follows, we regard matrices with entries in
$\Aloc[z,\,1/(z-\zeta)]$
and those with entries in
$\quotfield[z,1/z]$ as matrices over
$\Acomp[[z,1/z\}[\omegav]$ via these embeddings.

\begin{Theorem}[{\cite[Lemma 3.5.8]{Hartl2020a}}] \label{Thm:comparison}
    Let $v$ be a prime ideal of $A$ of degree $1$ with a monic generator $\zeta$ and let $z\in \mathbb{F}_q[t]$ be $\zeta|_{\theta=t}$.
    Take $M \in \mathcal{M}$ defined by $\Phi_M =(\Phi_{ij})\in \operatorname{GL}_r(\Aloc[z,\,1/(z-\zeta)])$ and its reduction $\overline{M}\in \overline{\mathcal{M}}$ defined by $\Xi_M=(\Xi_{ij}):=(\overline{\Phi}_{ij}) \in \GL_r(\quotfield[z,\,1/z])$. 
    If we put 
    \begin{equation}
        \Delta_{M,\,m} :=   \Phi_{M}\Phi_{M}^{(1)}\cdots%\Phi_{M}^{(dm-d)}
        \Phi_{M}^{(m)}(\Xi_{M}^{-1})^{(m)}%(\Xi_{M}^{-1})^{(dm-d)}
        \cdots(\Xi_{M}^{-1})^{(1)}\Xi_{M}^{-1}
    \end{equation}
    for each $m\geq0$, then the sequence $\{\Delta_{M,\,m}\}_{m\geq0}$ converges in $\GL_r\big(\Acomp[[z,\,1/z\}[\omegav]\big)$ as $m \to \infty$ and the resulting matrix $\Delta_{M}$ satisfies
    \begin{equation}\label{eq:Deltaproperty1}
        \Delta_{M} \cdot \Xi_{M}=\Phi_{M} \cdot \Delta_{M}^{(1)}
    \end{equation}
    and
    \begin{equation}\label{eq:Deltaproperty2}
        \Delta_{M} \equiv I_r \pmod{v}.
    \end{equation}
    Moreover, the matrix in $\GL_r\big(\Acomp[[z,\,1/z\}[\omegav]\big)$ satisfying \cref{eq:Deltaproperty1,eq:Deltaproperty2} is unique.
\end{Theorem}
It should be mentioned that \cref{Thm:comparison}, which goes back to \cite{Genestier2011}, has also been established for finite places $v$ of higher degree and over more general base rings; see \cite{Hartl2020a}.
The matrix $\Delta_M$ in the statement of \cref{Thm:comparison} is the inverse of the matrix constructed in \cite{Hartl2020a}.

\begin{Remark}\label{Rem:Topology_and_convergence}
The convergence of entries asserted in \cref{Thm:comparison} is convergence in the topological localization $\Acomp[[z,1/z\}[\omegav]=\Acomp[[z,1/z\}\left[1/\Omegav^{(-1)}\right]$.
Indeed, the estimates in the proof of \cite[Lemma 2.3.1]{Hartl2011} imply that, after multiplication by a sufficiently large fixed power of $\Omegav^{(-1)}$, the entries of $\Delta_{M,\,m}$ converge to those of $\Delta_M$ with respect to each Gauss norm $||\cdot||_r$ on $\Acomp[[z,1/z\}$.
The asserted convergence then follows from the continuity of the natural homomorphism to the localization and of multiplication by the inverse of this fixed power.
We claim that the convergence in this sense implies coefficient-wise convergence of the corresponding bidirectional power series. 
Let us consider the ring
\begin{equation}
    \Kcomp\left\langle \frac{z}{\zeta^{\frac{1}{2q}}},\, \frac{\zeta^{\frac{1}{q}}}{z}
    \right\rangle:=\left\{ \sum_{i=-\infty}^\infty a_i z^{i}\, \middle| a_i\in \Kcomp,\,\, \lim_{i \to \infty}|a_{i}|_v|\zeta^{i}|_v^{\frac{1}{2q}} =\lim_{i \to -\infty}|a_{i}|_v|\zeta^{i}|_v^{\frac{1}{q}}= 0\right\}
\end{equation}
whose norm is given by the maximum of the two Gauss norms as follows:
\begin{equation}
    ||f||_{\left[\frac{1}{2q},\,\frac{1}{q}\right]}:=\max\{||f||_{\frac{1}{2q}},\,||f||_{\frac{1}{q}}\},
\end{equation}
see \cite[\S 2.8]{Hartl2009}. 
Then the inclusion of $\Acomp[[z,1/z\}$ into this ring is continuous and $\Omegav^{(-1)}$ is invertible in the latter ring. 
Thus, by the universal property of topological localization,
this inclusion extends to a continuous injective homomorphism
from $\Acomp[[z,1/z\}[\omegav]$. 
Therefore, a given convergent sequence in $\Acomp[[z,1/z\}[\omegav]$ also converges in $\Kcomp\left\langle \frac{z}{\zeta^{\frac{1}{2q}}},\, \frac{\zeta^{\frac{1}{q}}}{z}
    \right\rangle$; and hence the claim follows.
Moreover, for every $f\in \Acomp[[z,1/z\}[\omegav]$, the series $f^{(1)}$ converges at $z=\zeta$ (\cite[Example 3.5.15]{Hartl2020a}) and belongs to $\Kcomp\left\langle \frac{z}{\zeta^{\frac{1}{2q}}},\, \frac{\zeta^{\frac{1}{q}}}{z}
    \right\rangle$;   
    further, since the coefficients of $f$ lie in $\Acomp$,
we have
\begin{equation}
    || f^{(1)}||_{\left[\frac{1}{2q},\,\frac{1}{q}\right]}
    \leq || f||_{\left[\frac{1}{2q},\,\frac{1}{q}\right]},
    \qquad
    |f^{(1)}(\zeta)|_v
    \leq || f||_{\left[\frac{1}{2q},\,\frac{1}{q}\right]}^q.
\end{equation}
Consequently, convergence $f_m\to f$ in $\Acomp[[z,1/z\}[\omegav]$ implies convergences $f_m^{(1)}\to f^{(1)}$ in the above norm and $f_m^{(1)}(\zeta)\to f^{(1)}(\zeta)$ in $\Kcomp$.
\end{Remark}

\begin{Example}[{see \cite[Proof of Lemma 3.5.8]{Hartl2020a}}]\label{Example:comparisonisomCarlitz}
    Let us consider the object $\mathcal{C}$ in \cref{Example:Carlitz}; then \cref{Thm:comparison} implies that the limit matrix $\Delta_{\mathcal{C}}=\lim_{m\to \infty}\Delta_{\mathcal{C},\,m}$ with
    \begin{equation}\label{eq:def_of_delta_Cm}
        \Delta_{\mathcal{C},\,m}=\begin{pmatrix}\delta_{\mathcal{C},\,m}\end{pmatrix},\quad 
        \delta_{\mathcal{C},\,m}:=\prod_{i=0}^{m}\left(1-\frac{\zeta^{q^i}}{z}\right)=\frac{\lambda_{m}}{z^{m+1}}
    \end{equation}
    satisfies the equation $\Delta_{\mathcal{C}}\cdot \Xi_{\mathcal{C}}=\Phi_{\mathcal{C}}\cdot\Delta_{\mathcal{C}}^{(1)}$ with $\Xi_{\mathcal{C}}=(z)$ and $\Phi_{\mathcal{C}}=(z-\zeta)$. 
   By the definition of $\omegav$ in \cref{eq:def_omegav}, the limit is given by $\Delta_{\mathcal{C}}=(\omegav^{-1})$; that is, we have
   \begin{equation}\label{eq:delta_and_omegav}
       \lim_{m\to \infty}\delta_{\mathcal{C},\,m}=\omegav^{-1}.
   \end{equation}
   In other words, the function $\omegav$, which generates the ambient ring $\Acomp[[z,\,1/z\}[\omegav]$, is nothing but the reciprocal of the entry of the comparison matrix of the Carlitz motive.
   We further note the equation
   \begin{equation}\label{eq:delta_Omega}
    \Delta_{\mathcal{C}}^{(1)}=\begin{pmatrix}\lim_{m\to \infty} \delta_{\mathcal{C},\,m}^{(1)}\end{pmatrix}=\begin{pmatrix}\Omegav\end{pmatrix},
    \end{equation}
   see \cref{eq:def_Omegav}.
\end{Example}

\begin{Remark}\label{Rem:comparison-realizations}
For $M$ in $\mathcal{M}$, the de Rham realization over $K$ is given by
\begin{equation}
\operatorname{H}^1_{\mathrm{dR}}(M,K)
:=
 \tau^* M
\otimes_{\Aloc[z]}K
\end{equation}
where the tensor product is taken via the homomorphism $\Aloc[z]\to K$ sending $z$ to $\zeta$; the crystalline realization of $M$ is given by 
\begin{equation}
    \operatorname{H}_{\mathrm{cris}}^{1}(M,\Kcomp)
:=
 \tau^* \overline{M} 
\otimes_{\quotfield[z]}\Kcomp
\end{equation}
where the tensor product is taken via the Teichm\"uller section $\quotfield \hookrightarrow \Kcomp$ and the substitution $z \mapsto \zeta$; see \cite[Definition 3.5.14]{Hartl2020a}.
With respect to the bases induced by the chosen basis of $M$, the evaluation $\Delta_M^{(1)}|_{z=\zeta}\in\GL_r(\Kcomp)$ represents an isomorphism from the crystalline realization to the de Rham realization after base change to $\Kcomp$ or to $\mathbb{C}_v$.
This is the inverse of the comparison isomorphism given in \cite[Theorem 3.5.18]{Hartl2020a}.
\end{Remark}

We let $M_i$ be an object of $\mathcal{M}$ which has the defining matrix $\Phi_{M_i} \in \GL_{r_i}(\Aloc[z,$ $1/(z-\zeta)])$ for $i=1,\,2$, and consider the tensor product $M_1\otimes M_2$, which is defined by the Kronecker product $\Phi_{M_1}\otimes \Phi_{M_2} \in \GL_{r_1r_2}(\Aloc[z,\,1/(z-\zeta)])$.
In this situation we have $\Delta_{M_1\otimes M_2}=\Delta_{M_1}\otimes \Delta_{M_2}$ by \cref{eq:Kronecker_mixed,eq:Kronecker_twist} and the uniqueness established in \cref{Thm:comparison}.
%\subsubsection{Examples: $v$-adic multiple zeta values}\label{Section:v-adicpolyloginterpretation}
\begin{Example}\label{Example:v-adicpolylogcomparison}
Let us observe that deformation series of special values of $v$-adic Carlitz multiple polylogarithms appear as entries of the matrix given in \cref{Thm:comparison}.
We assume $\deg v=1$, which gives us $\lambda_i=\left(z-\zeta^{(0)}\right)\cdots\left(z-\zeta^{(i)}\right)$ for $i\geq 0$ (see \cref{eq:def_lambda,eq:deg1consequence}); thus, if we take $\tilde{\mathbf{a}}=\big((s_1,\,a_1),\,\dots,\,(s_r,\,a_r)\big) \in \tilde{\mathbb{D}}_v^r$ and $m \geq 0$, then we have
\begin{equation}\label{eq:deftruncatedprepolylogv}
   \prelog_{\tilde{\mathbf{a}},\,m}= \sum_{m\geq i_1>\cdots>i_r\geq0}\frac{a_1^{(i_1)}\cdots a_r^{(i_r)}}{\big((z-\zeta^{(0)})\cdots(z-\zeta^{(i_1)})\big)^{s_1}\cdots\big((z-\zeta^{(0)})\cdots(z-\zeta^{(i_r)})\big)^{s_r}},
\end{equation}
see \cref{eq:deftruncatedprepolylog}.

For the chosen $\tilde{\mathbf{a}}\in \tilde{\mathbb{D}}_v^r$, 
the matrix $\Phi_{\tilde{\mathbf{a}}}$ given in \cref{eq:def matrix T} can be written as
\begin{equation}
\begin{pmatrix}
  1
  & 0 & \cdots & \cdots & 0\\%[25pt]
 \frac{a_1}{(z-\zeta)^{s_1}} & \frac{1}{(z-\zeta)^{s_1}} & \ddots & & \vdots\\%[25pt]
  0
  & \ddots
  & \frac{1}{(z-\zeta)^{s_1+s_2}}
  &\ddots& \vdots\\%[25pt]
  \vdots
  &\ddots&
  \ddots
  & \ddots
  &0\\%[25pt]
  0 & \cdots & 0&  \frac{a_r}{(z-\zeta)^{s_1+\cdots +s_r}}
  & \frac{1}{(z-\zeta)^{s_1+\cdots +s_r}}
\end{pmatrix}
\end{equation}
again by \cref{eq:deg1consequence}, and lies in $\GL_{r+1}(\Aloc[z,\,1/(z-\zeta)])$.
Let us consider the object $M_{\tilde{\mathbf{a}}}$ of $\mathcal{M}$ defined by this matrix and its modulo $v$ reduction $\overline{M}_{\tilde{\mathbf{a}}}$, which is defined by

\begin{equation}\renewcommand{\arraystretch}{1}
    \Xi_{\tilde{\mathbf{a}}}:=\begin{pmatrix}
  1
  & 0& \cdots & \cdots & 0\\
  \mathmakebox[\PowW][c]{\frac{\overline{a}_1}{(z-\overline{\zeta})^{s_1}}} & \frac{1}{(z-\overline{\zeta})^{s_1}} & \ddots &  & \vdots\\
  0
  & \frac{\overline{a}_2}{(z-\overline{\zeta})^{s_1+s_2}}
  & \frac{1}{(z-\overline{\zeta})^{s_1+s_2}}
  &\ddots& \vdots\\
  \vdots
  &\ddots&
  \ddots
  & \ddots
  & 0\\
  0 & \cdots & 0 & \frac{\overline{a}_r}{(z-\overline{\zeta})^{s_1+\cdots+s_r}}
  & \frac{1}{(z-\overline{\zeta})^{s_1+\cdots+s_r}}
\end{pmatrix}
\end{equation}

    Here, we have $a_i \in v\Aloc$ for $1 \leq i \leq r$, that is, $\overline{a}_i=0\in \mathbb{F}_v$ by assumption. Then the matrix $\Xi_{\tilde{\mathbf{a}}}$ is of the simplified form
    \begin{equation}\label{eq:Xi_s}
        \begin{pmatrix}
            1&0&\dots&0\\
            0& \frac{1}{(z-\overline{\zeta})^{s_1}}&\ddots&\vdots\\
            \vdots&\ddots&\ddots&0\\
            0&\cdots&0&\frac{1}{(z-\overline{\zeta})^{s_1+\cdots+s_r}}
        \end{pmatrix}=\begin{pmatrix}
            1&0&\dots&0\\
            0&z^{-s_1}&\ddots&\vdots\\
            \vdots&\ddots&\ddots&0\\
            0&\cdots&0&z^{-s_1-\cdots-s_r}
        \end{pmatrix}.
    \end{equation}
    Let us calculate the matrix $\Delta_{M_{\tilde{\mathbf{a}}}}$ in \cref{Thm:comparison} following \cite{Hartl2020a}. 
    To this end, we define $\prelog_{ji,\,m}$ to be $\prelog_{\tilde{\mathbf{a}}^\prime,\,m}$ with $\tilde{\mathbf{a}}^\prime:=\big((s_{i},\,a_{i}),\,\dots,\,(s_{j-1},\,a_{j-1})\big)$ for $1 \leq i<j \leq r+1$ (see \cref{eq:deftruncatedprepolylogv});
    then, for $m\geq 0$, we can show by induction that
    \begin{align}
        &\quad \ \ \Phi_{\tilde{\mathbf{a}}}\Phi_{\tilde{\mathbf{a}}}^{(1)}\cdots\Phi_{\tilde{\mathbf{a}}}^{(m)}(\Xi_{\tilde{\mathbf{a}}}^{-1})^{(m)}\cdots(\Xi_{\tilde{\mathbf{a}}}^{-1})^{(1)}\Xi_{\tilde{\mathbf{a}}}^{-1}\\
        &=
        \begin{pmatrix}
1 & 0 & \cdots & \cdots & 0\\
\prelog_{21,\,m} & \frac{1}{(\lambda_m)^{s_1}} & 0 & \cdots & 0\\
\vdots & \frac{\prelog_{32,\,m}}{(\lambda_m)^{s_1}} & \frac{1}{(\lambda_m)^{s_1+s_2}} & \ddots & \vdots\\[6pt]
\vdots & \vdots & \ddots & \ddots & 0\\[6pt]
\prelog_{(r+1)1,\,m} & \frac{\prelog_{(r+1)2,\,m}}{(\lambda_m)^{s_1}}& \cdots & \frac{\prelog_{(r+1)r,\,m}}{(\lambda_m)^{s_1+\cdots+s_{r-1}}} & \frac{1}{(\lambda_m)^{s_1+\cdots+s_r}}
\end{pmatrix}\\
        &\quad \ \ \cdot\begin{pmatrix}
            1&0&\dots&0\\
            0&z^{s_1}&\ddots&\vdots\\
            \vdots&\ddots&\ddots&0\\
            0&\cdots&0&z^{s_1+\cdots+s_r}
        \end{pmatrix}^{m+1}\\
        &=\begin{pmatrix}
1 & 0 & \cdots & \cdots & 0\\
\prelog_{21,\,m} & \frac{1}{\delta_{\mathcal{C},\,m}^{s_1}} & 0 & \cdots & 0\\
\vdots & \frac{\prelog_{32,\,m}}{\delta_{\mathcal{C},\,m}^{s_1}} & \frac{1}{\delta_{\mathcal{C},\,m}^{s_1+s_2}} & \ddots & \vdots\\[6pt]
\vdots & \vdots & \ddots & \ddots & 0\\[6pt]
\prelog_{(r+1)1,\,m} & \frac{\prelog_{(r+1)2,\,m}}{\delta_{\mathcal{C},\,m}^{s_1}}& \cdots & \frac{\prelog_{(r+1)r,\,m}}{\delta_{\mathcal{C},\,m}^{s_1+\cdots+s_{r-1}}} & \frac{1}{\delta_{\mathcal{C},\,m}^{s_1+\cdots+s_r}}
\end{pmatrix},
\end{align}
see also \cref{eq:def_of_delta_Cm}.
Considering the limit as $m \to \infty$ of the sequence 
\begin{equation}
\left\{
\Phi_{\tilde{\mathbf{a}}}\Phi_{\tilde{\mathbf{a}}}^{(1)}\cdots\Phi_{\tilde{\mathbf{a}}}^{(m)}(\Xi_{\tilde{\mathbf{a}}}^{-1})^{(m)}\cdots(\Xi_{\tilde{\mathbf{a}}}^{-1})^{(1)}\Xi_{\tilde{\mathbf{a}}}^{-1}
\right\}_m,
\end{equation}
we can deduce, by \cref{Thm:comparison}, that the sequence converges.
Since $\delta_{\mathcal{C},\,m}$ converges to the invertible element $\omegav^{-1}$ of $\Acomp[[z,\,1/z\}[\omegav]$ (\cref{eq:delta_and_omegav}), the limit
\begin{equation}
    \lim_{m\to \infty}\prelog_{\tilde{\mathbf{a}}^\prime,\,m} \in \Acomp[[z,\,1/z\}[\omegav]
\end{equation}
exists for each tuple $\tilde{\mathbf{a}}^\prime:=\big((s_{i},\,a_{i}),\,\dots,\,(s_{j-1},\,a_{j-1})\big)$ with $1 \leq i <j \leq r+1$.
Since the tuple $\tilde{\mathbf{a}}\in \tilde{\mathbb{D}}_v^r$ was taken arbitrarily, the limit
\begin{equation}\label{eq:def of t-motivic CMPLs}
    \prelog_{v,\,\tilde{\mathbf{b}}}:=\lim_{m\to \infty}\prelog_{\tilde{\mathbf{b}},\,m} \in \Acomp[[z,\,1/z\}[\omegav]
\end{equation}
exists for every $r^\prime \geq 1$ and $\tilde{\mathbf{b}} \in \tilde{\mathbb{D}}_v^{r^\prime}$; by \cref{Rem:Topology_and_convergence}, we can regard the limit $\prelog_{v,\,\tilde{\mathbf{b}}}$ as a bidirectional power series.

Keeping the tuple $\tilde{\mathbf{a}}$ fixed as above, we use the abbreviations
$\prelog_{v,\,ji}:=\prelog_{v,\,\tilde{\mathbf{a}}^\prime}$ with $\tilde{\mathbf{a}}^\prime:=\big((s_{i},\,a_{i}),\,\dots,\,(s_{j-1},\,a_{j-1})\big)$ for $1 \leq i <j \leq r+1$, and $\prelog_{v,\,ii}:=1$ for $1 \leq i \leq r+1$, in parallel with the symbols in \cref{Section:infty-side}; we note that the notations $\prelog_{ji,\,m}$ and $\prelog_{v,\,ji}$ depend on the tuple $\tilde{\mathbf{a}}$ fixed at the beginning.
Then, the matrix $\Delta_{M_{\tilde{\mathbf{a}}}}$ given by the theorem is expressed as
\begin{equation}
\begin{pmatrix}
1 & 0 & \cdots & \cdots & 0\\
\prelog_{v,\,21} & \omegav^{s_1} & 0 & \cdots & 0\\
\vdots & \omegav^{s_1}\prelog_{v,\,32} & \omegav^{s_1+s_2} & \ddots & \vdots\\
\prelog_{v,\,r1} & \vdots & \ddots & \ddots & 0\\
\prelog_{v,\,(r+1)1} & \omegav^{s_1}\prelog_{v,\,(r+1)2} & \cdots & \omegav^{s_1+\cdots+s_{r-1}}\prelog_{v,\,(r+1)r} & \omegav^{s_1+\cdots+s_r}
\end{pmatrix}.
\end{equation}
It holds that 
\begin{equation}\label{eq:prelogvdifferenceeq}
        \prelog_{v,\,ji}=\sum_{k=i}^{j}\frac{\prelog_{jk,\,m}\,\prelog_{v,\,ki}^{(m+1)}}{\lambda_m^{s_i+\cdots+s_{k-1}}}
    \end{equation}
    for any $m\geq 0$ if we put $\prelog_{jj,\,m}=1$ for $1\leq j \leq r+1$ by convention.
\end{Example}

Chen proved that these series are in the field of fractions of the ring $\Entirev$ of entire functions. More precisely,

\begin{Proposition}\label{Prop:v-adic polylog entire}
The series $\Omegav^{s_1+\cdots+s_r}$ $\prelog_{v,\,\tilde{\mathbf{a}}}^{(1)}$ is entire for $\tilde{\mathbf{a}}=\big((s_1,\,a_1),$ $\dots,$ $(s_r,\,a_r)\big)\in \tilde{\mathbb{D}}_v^r$ if $\deg v=1$.
\end{Proposition}

\begin{proof}
In \cite[Proposition 3.2.2]{Chen2025}, Chen proved that the series
    
    \begin{align}
        &\Omegav^{s_1+\cdots+s_r}\prelog_{v,\,\tilde{\mathbf{a}}}^{(1)}+a_r\Omegav^{s_1+\cdots+s_r}\prelog_{v,\,\tilde{\mathbf{a}}^\prime}^{(1)}\\
        &=\Omegav^{s_1+\cdots+s_r}\sum_{i_1>i_2>\cdots>i_r\geq 0}\frac{{a_1^{(i_1)}}\cdots a_r^{(i_r)}}{\big(\left(z-\zeta^{q}\right)\cdots\left(z-\zeta^{q^{i_1}}\right)\big)^{s_1}\cdots\big(\left(z-\zeta^{q}\right)\cdots\left(z-\zeta^{q^{i_r}}\right)\big)^{s_r}} 
    \end{align}
    is entire, where we put $\tilde{\mathbf{a}}^\prime=\big((s_1,\,a_1),\,\dots,\,(s_{r-1},\,a_{r-1})\big)$. By induction on $r$, we can show that the series $\Omega_v^{s_1+\cdots+s_r}\prelog_{v,\,\tilde{\mathbf{a}}}^{(1)}$ is also entire.
\end{proof}

\begin{Example}\label{Example:comparison-non-conv}
    Let us put $\tilde{\mathbf{a}}=((s,\,a))$ with $s \in \mathbb{N}$ and $a \in \Aloc\setminus \zeta\Aloc$, and consider the object $M_{\tilde{\mathbf{a}}}$ defined by the same matrix as in \cref{Example:v-adicpolylogcomparison}.
    Then the matrix $\Delta_{M_{\tilde{\mathbf{a}}}}$ is given as follows:
    \begin{align}
    &\quad \ \ \lim_{m \to \infty}\left(\begin{pmatrix}
        1&0\\
        \prelog_{\tilde{\mathbf{a}},\,m}& \left( \lambda_m \right)^{-s}
    \end{pmatrix}\cdot
    \begin{pmatrix}
        1&0\\
        -\overline{a} (1+z^{s}+z^{2s}+\cdots+z^{ms})&z^{(m+1)s}
    \end{pmatrix}
    \right)\\
    &=\lim_{m \to \infty}
    \begin{pmatrix}
        1&0\\
        \prelog_{\tilde{\mathbf{a}},\,m}-\overline{a} \left( \lambda_m \right)^{-s}(1+z^{s}+z^{2s}+\cdots+z^{ms}) &\delta_{\mathcal{C},\,m}^{-s}
    \end{pmatrix}.
    \end{align}
    We note here that the limit 
    \begin{align}
    &\quad \lim_{m \to \infty}\left( \prelog_{\tilde{\mathbf{a}},\,m}-\overline{a} \left( \lambda_m \right)^{-s}(1+z^{s}+z^{2s}+\cdots+z^{ms})\right)\\
    &=\lim_{m \to \infty}\left( \prelog_{\tilde{\mathbf{a}},\,m}-\overline{a} \left( \lambda_m \right)^{-s} z^{(m+1)s}( z^{-s}+z^{-2s}+\cdots+z^{-(m+1)s})\right)\\
   &= \lim_{m \to \infty}\left( \prelog_{\tilde{\mathbf{a}},\,m}-\overline{a} \delta_{\mathcal{C},\,m}^{-s}( z^{-s}+z^{-2s}+\cdots+z^{-(m+1)s} )\right)
    \end{align}
    exists as we have \cref{Thm:comparison} and \cref{Rem:Topology_and_convergence}; we emphasize that, in the setting of this example, the limit is taken for the combined expression as a whole, including the term  $-\overline{a}\,\delta_{\mathcal{C},\,m}^{-s}(z^{-s}+z^{-2s}+\cdots+z^{-(m+1)s})$, and we do not claim the convergence of each part separately, in contrast with the case $\tilde{\mathbf{a}}\in\tilde{\mathbb{D}}_v^r$ treated in \cref{Example:v-adicpolylogcomparison}, where the series $\prelog_{v,\,\tilde{\mathbf{a}}}$ itself converges.
\end{Example}

\begin{Example}\label{Example:hyperderivative_comparison}

In the study of $\infty$-adic periods, hyperderivatives of series and Maurischat's prolongations of $t$-motives enable us to realize Taylor coefficients of deformation series of the periods as periods of suitable $t$-motives; see \cite{Maurischat2018,Namoijam2024}. 
This approach has been applied, in particular, to Taylor coefficients of deformation series arising in the study of Thakur's multiple zeta values (\cite{Chang2021a,Matsuzuki2024a}). 
Motivated by these works, we give a simple example in which the first Taylor coefficient of a pre-period can again occur as a period; see also \cref{Example:hyperderivative_period}.

We define the $\Acomp$-linear operator $\partial_z$ on $ \Acomp[[z,1/z\}$ by
    \begin{equation}
    \partial_z:\sum_{i \in \mathbb{Z}}
    a_i z^{i}  \mapsto \sum_{i \in \mathbb{Z}}\binom{i}{1}
    a_i z^{i-1},
\end{equation}
as in \cite{Hasse1936}, and extend it to the localization $\Acomp[[z,1/z\}[\omegav]=\Acomp[[z,1/z\}[1/\Omegav^{(-1)}]$ by $\partial_z(\omegav)=-\omegav^2\partial_z(\omegav^{-1})$.
Taking $s \geq 1$ and $a \in \Aloc$ with $|a|_v<1$, we put $\tilde{\mathbf{a}}:=((s,\,a))$. 
If we consider the matrices given by
\begin{align}
    \rho_z\Phi_{\tilde{\mathbf{a}}}&:=\begin{pmatrix}
        1&0&0&0\\
        \frac{a}{(z-\zeta)^s}&\frac{1}{(z-\zeta)^s}&0&0\\
        0&0&1&0\\
        \partial_z\left( \frac{a}{(z-\zeta)^s}\right)&\partial_z\left(\frac{1}{(z-\zeta)^s} \right)&\frac{a}{(z-\zeta)^s}&\frac{1}{(z-\zeta)^s}
    \end{pmatrix}\\&=
    \begin{pmatrix}
        1&0&0&0\\
        \frac{a}{(z-\zeta)^s}&\frac{1}{(z-\zeta)^s}&0&0\\
        0&0&1&0\\
        -\frac{sa}{(z-\zeta)^{s+1}}&-\frac{s}{(z-\zeta)^{s+1}}&\frac{a}{(z-\zeta)^s}&\frac{1}{(z-\zeta)^s}
    \end{pmatrix},
\end{align}
\begin{equation}
    \rho_z\Xi_{\tilde{\mathbf{a}}}:=\overline{\rho_z\Phi_{\tilde{\mathbf{a}}}}=
    \begin{pmatrix}
        1&0&0&0\\
        0&\frac{1}{z^s}&0&0\\
        0&0&1&0\\
        0&-\frac{s}{z^{s+1}}&0&\frac{1}{z^s}
    \end{pmatrix} ,
\end{equation}
and
\begin{equation}
    \rho_z\Delta_{M_{\tilde{\mathbf{a}}}}:=
    \begin{pmatrix}
        1&0&0&0\\
        \prelog_{v,\,\tilde{\mathbf{a}}}&\omegav^s&0&0\\
        0&0&1&0\\
        \partial_z\prelog_{v,\,\tilde{\mathbf{a}}}&\partial_z \omegav^s&\prelog_{v,\,\tilde{\mathbf{a}}}&\omegav^s
    \end{pmatrix},
\end{equation}
then we have $\rho_z\Delta_{M_{\tilde{\mathbf{a}}}} \cdot \rho_z\Xi_{\tilde{\mathbf{a}}}= \rho_z\Phi_{\tilde{\mathbf{a}}}\cdot \rho_z\Delta_{M_{\tilde{\mathbf{a}}}}^{(1)}$.
Let $\rho_z M_{\tilde{\mathbf{a}}}$ be the object of $\mathcal{M}$ defined by $\rho_z\Phi_{\tilde{\mathbf{a}}}$; then its reduction $\rho_z\overline{M}_{\tilde{\mathbf{a}}}$ in $\overline{\mathcal{M}}$ is represented by $\rho_z\Xi_{\tilde{\mathbf{a}}}$. 
By the uniqueness established in \cref{Thm:comparison}, we see that the matrix $\Delta_{\rho_zM_{\tilde{\mathbf{a}}}}$ given by the theorem coincides with $\rho_z\Delta_{M_{\tilde{\mathbf{a}}}}$ since we have $\rho_z\Delta_{M_{\tilde{\mathbf{a}}}}\equiv I_4 \mod v$.
\end{Example}

\subsection{The formulation of $v$-adic periods}\label{subsec:v-adic_periods}
This subsection is devoted to the formulation of Andr\'{e}'s $v$-adic periods and a discussion of their examples.
The notion of Andr\'{e}'s $v$-adic periods for an object $M \in \mathcal{M}$ comes from the isomorphism in \cref{Thm:comparison} between $M$ and its reduction $\overline{M}$ tensored with $\Acomp[[z,\,1/z\}[\omegav]$. Before formulating Andr\'{e}'s $v$-adic periods themselves, we introduce the notion of pre-periods for convenience. We note that the term ``pre-periods'' is not standard, though the notion and similar ones frequently appear in the study of periods in the arithmetic of function fields.

\begin{definition}\label{Def:pre-period}
    Take $M \in \mathcal{M}$ defined by $\Phi_M =(\Phi_{ij})\in \operatorname{GL}_r(\Aloc[z,\,1/(z-\zeta)])$, and let $\mathbf{m}=(m_1,\,\dots,\,m_r)$ be the $\Aloc[z]$-basis corresponding to $\Phi_M$.
    Consider the reduction $\overline{M}$ of $M$ defined by $\Xi_M=(\overline{\Phi}_{ij}) \in \GL_r(\quotfield[z,\,1/z])$ and let $\overline{\mathbf{m}}=(\overline{m}_1,\,\dots,\,\overline{m}_r)$ be the row vector of $\quotfield[z]$-basis elements induced by $\mathbf{m}$ and corresponding to $\Xi_M$.
    Let $\Delta_M \in \GL_r\big(\Acomp[[z,\,1/z\}[\omegav]\big)$ be the matrix given by \Cref{Thm:comparison}.
    \begin{enumerate}
        \item Fix the $\quotfield[z]$-basis of $\overline{\mathbf{1}}$; for a morphism $\iota:\overline{\mathbf{1}}\rightarrow \overline{M}$ represented by $\nu_\iota \in \Mat_{r\times 1}(\quotfield[z])$ with respect to the basis $\overline{\mathbf{m}}$, the vector
    \begin{equation}
        \psi_{M,\, \iota}:= \Delta_M^{{(1)}}\cdot \nu_\iota^{{(1)}} \in\Mat_{r\times 1}\Big(\Acomp[[z,\,1/z\}[\omegav^{(1)}]\Big)
    \end{equation}
    is called the \textit{$v$-adic pre-period vector} of $M$ corresponding to $\iota$; note that the twisting defined in \cref{eq:def_v_twisting} preserves the ring $\Acomp[[z,\,1/z\}$.
    \item The entries of a $v$-adic pre-period vector are called \textit{$v$-adic pre-periods}. 
    \end{enumerate}
\end{definition}

%In \Cref{sec:difference galois}, we will attempt to  study the transcendence of pre-period via the theory of difference equations. 
Every element of $\Acomp[[z,\,1/z\}[\omegav^{(1)}]$ converges at $z=\zeta$.
Moreover, expansion at $z=\zeta$ defines an embedding of this ring into $\Kcomp[[z-\zeta]]$
(\cite[Theorem 3.5.18]{Hartl2020a}); indeed, the embedding of $\Acomp[[z,\,1/z\}$ is given by
\begin{equation}
    \sum_{i=-\infty}^\infty a_i z^{i}=\sum_{j=0}^{\infty} \zeta^{-j}\left(
    \sum_{i \in \mathbb{Z}}\binom{i}{j}a_i\zeta^i
    \right)
    (z-\zeta)^j
\end{equation}
and $\omegav^{(1)}=\Omegav^{-1}$ also has a convergent power-series expansion at $z=\zeta$, since $\Omegav$ does not vanish there; see \cite[Example 3.5.15]{Hartl2020a}.
Thus, pre-periods converge at $z=\zeta$ by \cref{Rem:Topology_and_convergence}.

\begin{definition}\label{def:period}
    Take $M \in \mathcal{M}$ which has a $v$-adic pre-period $P \in \Acomp[[z,\,1/z\}[\omegav^{(1)}]$; then the value $P|_{z=\zeta}$ is called \textit{Andr\'{e}'s $v$-adic period}.
    If $M$ has a pre-period vector $\psi$, then the $\Kcomp$-vector $\psi|_{z=\zeta}$ is called a \textit{$v$-adic period vector}.
\end{definition}

Although the individual entries of a $v$-adic period vector depend on the choice of basis, their $\Kalg$-linear span does not.

Let us indicate how the above definition realizes the idea described in \Cref{subsec:intro_idea}.
If we consider the situation of \cref{Rem:comparison-realizations} and \cref{Def:pre-period}, a morphism $\iota:\overline{\mathbf{1}} \rightarrow \overline{M}$ represented by $\nu_\iota=(\nu_1,\,\dots,\,\nu_r)^\top \in \Mat_{r\times 1}(\quotfield[z])$ gives rise to an element 
\begin{equation}
   ( \overline{m}_1\otimes \tau(\nu_1)|_{z=\zeta})+\cdots+( \overline{m}_r\otimes \tau(\nu_r)|_{z=\zeta})=(\overline{\mathbf{m}}\otimes 1)\cdot \big(\nu_\iota^{(1)}|_{z=\zeta}\big)
\end{equation}
of the crystalline realization of $M$; 
via the comparison isomorphism represented by $\Delta_{M}^{(1)}|_{z=\zeta}$, this element is transported to the element
\begin{equation}
   ( \mathbf{m} \otimes 1)\cdot\left(\Delta_{M}^{(1)}|_{z=\zeta} \right)\cdot \big(\nu_\iota^{(1)}|_{z=\zeta}\big)
\end{equation}
of the de Rham realization of $M$ after base change to $\Kcomp$. 
Thus, the $v$-adic period vector is the coordinate vector of this element with respect to the basis $\mathbf{m} \otimes 1$ of the de Rham realization induced by the chosen basis $\mathbf{m}$ of $M$. 
Therefore, the definition above may be viewed as a function-field counterpart of the idea of Ancona and Fr\u{a}\c{t}il\u{a}, where the morphisms $\iota:\overline{\mathbf{1}}\rightarrow \overline{M}$ play the role of the cycles of the reduction modulo $p$ of motives used in \cite{Ancona2025}.

\begin{Example}\label{Example:v-adicpolyloginterpretation}
    We continue to consider the case discussed in \cref{Example:v-adicpolylogcomparison}.
    Since the matrix $\Xi_{\tilde{\mathbf{a}}}$ representing the structure isomorphism $\tau_{\overline{M}_{\tilde{\mathbf{a}}}}$ of $\overline{M}_{\tilde{\mathbf{a}}}$ is given as in \cref{eq:Xi_s}, we have a morphism ${\iota_{\tilde{\mathbf{a}}}}:\overline{\mathbf{1}}\rightarrow \overline{M}_{\tilde{\mathbf{a}}}$, whose representing matrix $\nu_{{\iota_{\tilde{\mathbf{a}}}}}$ is given by
    $(1,\,0,\,\dots,\,0)^\top$.
    The pre-period vector $\psi_{\tilde{\mathbf{a}}}:=\Delta_{M_{\tilde{\mathbf{a}}}}^{(1)}\cdot \nu_{\iota_{\tilde{\mathbf a}}}^{(1)}=\Delta_{M_{\tilde{\mathbf{a}}}}^{(1)}\cdot \nu_{\iota_{\tilde{\mathbf a}}}=\Delta_{M_{\tilde{\mathbf{a}}}}^{(1)}\cdot \begin{pmatrix}
            1&0&\cdots&0
        \end{pmatrix}^\top$
        associated with $\iota_{\tilde{\mathbf{a}}}$ is as follows:
    \begin{equation}\label{eq:v-adic polylog pre-period vector}
        \psi_{\tilde{\mathbf{a}}}=\begin{pmatrix}               1&\prelog_{v,\,21}^{(1)}&\cdots&\prelog_{v,\,r1}^{(1)}&\prelog_{v,\,(r+1)1}^{(1)}
            \end{pmatrix}^\top.
    \end{equation}
    By the construction of the series 
    \begin{equation}
        \prelog_{v,\,\big((s_1,\,a_1),\,\dots,\,(s_{j},\,a_{j})\big)} \text{ and }\prelog_{v,\,\big((s_1,\,a_1),\,\dots,\,(s_{j-1},\,a_{j-1})\big)},
    \end{equation} 
    and \cref{eq:deg1consequence}, it holds for each $1 \leq j \leq r$ that 
    \begin{align}
        &\prelog_{v,\,(j+1)1}^{(1)}+a_{j}\prelog_{v,\,j1}^{(1)}
        =\prelog_{v,\,\big((s_1,\,a_1),\,\dots,\,(s_{j},\,a_{j})\big)}^{(1)}+a_{j}\prelog_{v,\,\big((s_1,\,a_1),\,\dots,\,(s_{j-1},\,a_{j-1})\big)}^{(1)}\\
        &=\sum_{ i_1>\cdots>i_{j}\geq0}\frac{a_1^{(i_1)}\cdots a_{j}^{(i_{j})}}{\big(\left(z-\zeta^{(1)}\right)\cdots\left(z-\zeta^{(i_1)}\right)\big)^{s_1}\cdots\big(\left(z-\zeta^{(1)}\right)\cdots\left(z-\zeta^{(i_{j})}\right)\big)^{s_{j}}}.
    \end{align}
    Here, we put $\left(z-\zeta^{(1)}\right)\cdots\left(z-\zeta^{(i_j)}\right):=1$ if $i_j=0$ and, in the case $j=1$, set $a_{j}\prelog_{v,\,\big((s_1,\,a_1),\,\dots,\,(s_{j-1},\,a_{j-1})\big)}^{(1)}$ to be  $a_{1}$ by convention. Thus, by \cref{Rem:Topology_and_convergence}, we have
    \begin{equation}\label{eq:v_prelog_and_log}
        \mathrm{Li}\big((s_1,\,a_1),\,\dots,\,(s_{j},\,a_{j})\big)_{v}=\prelog_{v,\,(j+1)1}^{(1)}\Big|_{z=\zeta}+a_{j}\prelog_{v,\,j1}^{(1)}\Big|_{z=\zeta},
    \end{equation}
    as the defining series of $ \mathrm{Li}\big((s_1,\,a_1),\,\dots,\,(s_{j},\,a_{j})\big)_{v}$ (see \cref{eq:defCMPL}) can be rewritten as
    \begin{equation}\label{eq:defvCMPL}
    \sum_{ i_1>\cdots>i_{j}\geq0}\frac{a_1^{(i_1)}\cdots a_{j}^{(i_{j})}}{\big(\left(\zeta-\zeta^{(1)}\right)\cdots\left(\zeta-\zeta^{(i_1)}\right)\big)^{s_1}\cdots\big(\left(\zeta-\zeta^{(1)}\right)\cdots\left(\zeta-\zeta^{(i_{j})}\right)\big)^{s_{j}}}
\end{equation}
by \cref{eq:deg1consequence}.
    Consequently, we can interpret special values of $v$-adic Carlitz multiple polylogarithms in terms of $v$-adic periods.
\end{Example}

\begin{Example}\label{Example:preperiods-non-conv}
    Considering the object $M_{\tilde{\mathbf{a}}}$ in \cref{Example:comparison-non-conv}, we notice that its reduction $\overline{M}_{\tilde{\mathbf{a}}}$ admits a morphism $\iota:\overline{\mathbf{1}}\rightarrow \overline{M}_{\tilde{\mathbf{a}}}$ represented by $\nu_{\iota}:=(z^s-1,\,\overline{a})^\top$ as we have
    \begin{align}
        \Xi_{M_{\tilde{\mathbf{a}}}} \nu_\iota^{(1)}&=
        \begin{pmatrix}
            1&0\\
            \overline{a}/z^s&1/z^s
        \end{pmatrix}
        \begin{pmatrix}
            z^s-1\\
            \overline{a}
        \end{pmatrix}^{(1)}=
        \begin{pmatrix}
            1&0\\
            \overline{a}/z^s&1/z^s
        \end{pmatrix}
        \begin{pmatrix}
            z^s-1\\
            \overline{a}
        \end{pmatrix}\\
        &=\begin{pmatrix}
            z^s-1\\
            \overline{a}
        \end{pmatrix}=\nu_\iota I_1.
    \end{align}
    Hence, we have a pre-period
    \begin{align}
        &\quad \ \lim_{m\to \infty}
        \left(
        \left(z^s-1\right)
        \left(\prelog_{\tilde{\mathbf{a}},\,m}-\overline{a} \delta_{\mathcal{C},\,m}^{-s}( z^{-s}+z^{-2s}+\cdots+z^{-(m+1)s})\right)+\overline{a}\delta_{\mathcal{C},\,m}^{-s}
        \right)^{(1)}\\
        &=\lim_{m\to \infty}
        \left(
        \left(z^s-1\right)
        \prelog_{\tilde{\mathbf{a}},\,m}-\overline{a}\delta_{\mathcal{C},\,m}^{-s}(1-z^{-(m+1)s})+\overline{a}\delta_{\mathcal{C},\,m}^{-s}
        \right)^{(1)}\\
        &=\lim_{m\to \infty}
        \left(
        \left(z^s-1\right)
        \prelog_{\tilde{\mathbf{a}},\,m}+\overline{a}\delta_{\mathcal{C},\,m}^{-s}z^{-(m+1)s}
        \right)^{(1)}.
    \end{align}
\end{Example}

\begin{Example}\label{Example:hyperderivative_period}
    Given the situation in \cref{Example:hyperderivative_comparison}, we observe that the reduction $\rho_z\overline{M}_{\tilde{\mathbf{a}}}$ admits a morphism $\iota$ from $\overline{\mathbf{1}}$ represented by $\nu_{\iota}:=(1,\,0,\,0,\,0)^\top$; thus we have the pre-period vector \begin{equation}
        \psi_{\iota}=\left(1,\,\prelog_{v,\,\tilde{\mathbf{a}}}^{(1)},\,0,\,(\partial_z\prelog_{v,\,\tilde{\mathbf{a}}})^{(1)}\right)^\top.
    \end{equation}
    Therefore, we obtain Andr\'{e}'s $v$-adic period $(\partial_z\prelog_{v,\,\tilde{\mathbf{a}}})^{(1)}|_{z=\zeta}$; this coincides with the first Taylor coefficient of the pre-period $\prelog_{v,\,\tilde{\mathbf{a}}}^{(1)}$ at $z=\zeta$, while $\prelog_{v,\,\tilde{\mathbf{a}}}^{(1)}|_{z=\zeta}$ is the $0$-th coefficient; we refer the reader to \cite[Corollary 2.7]{Uchino1998} for the relationship between hyperderivatives and Taylor coefficients of series.
\end{Example}

Let us take objects $M_i$ of $\mathcal{M}$ defined by matrices $\Phi_{M_i} \in \GL_{r_i}(\Aloc[z,\,1/(z-\zeta)])$ for $i=1,\,2$, and consider the tensor product $M_1\otimes M_2$. We recall that the tensor product is defined by the Kronecker product $\Phi_{M_1}\otimes \Phi_{M_2}$ and the matrix $\Delta_{M_1\otimes M_2}$ given by \cref{Thm:comparison} is also given by the Kronecker product $\Delta_{M_1}\otimes \Delta_{M_2}$.
If we take a morphism $\iota_i:\overline{\mathbf{1}}\rightarrow \overline{M}_i$ defined by a column matrix $\nu_{\iota_i}$ with the pre-period vector $\psi_{M_i,\,\iota_i}:= \Delta_{M_i}^{(1)}\cdot \nu_{\iota_i}^{(1)} $ for each $i=1,\,2$, then the tensor product $\iota_1\otimes \iota_2: \overline{\mathbf{1}}\rightarrow \overline{M}_1\otimes \overline{M}_2$ is defined by the Kronecker product $\nu_{\iota_1}\otimes \nu_{\iota_2}$ and the corresponding pre-period vector is also given by the Kronecker product $\psi_{M_1,\,\iota_1} \otimes \psi_{M_2,\,\iota_2}$ by \cref{eq:Kronecker_mixed}.

\begin{Remark}\label{Rem:CWY}
This remark is included to indicate the scope of the idea underlying our formulation of $v$-adic periods: we show that our idea of Andr\'e's $v$-adic periods is not specific to Carlitz multiple polylogarithms; readers who are only interested in the proof of the main theorem may skip this remark.

For the $t$-motives associated with the Carlitz $\mathbb{F}_{q^l}[t]$-modules $\mathbf{C}_l$ ($l\geq1$), Chang, Wei, and Yu \cite[\S 3.7]{Chang2025} transported the Frobenius action on the crystalline realization to the de Rham realization over $\Kcomp$ through the comparison isomorphism.
They computed the matrix of the resulting $\Kcomp$-linear automorphism with respect to a basis of the de Rham realization over $K$, extended to $\Kcomp$.
Their Chowla--Selberg-type formula relates the nonzero entries of this matrix to special values of Goss' $v$-adic gamma function (\cite[Theorem 3.7.7]{Chang2025}).
Here, by considering the case where $l=3$ and $\deg v=1$, we shall observe that the nonzero entries of the matrix representing the inverse of the induced Frobenius action considered in \cite{Chang2025} can also be treated using the basic idea underlying our formulation of Andr\'{e}'s $v$-adic periods.

We mention here that our formulation does not apply directly in this case, since the auxiliary objects $M_{\mathbf{C}_3}^{\prime}$ and $M_{\mathbf{C}_3}\otimes M_{\mathbf{C}_3}^\prime$ introduced in this remark are not in $\mathcal{M}$ (\cref{Def:t-motive}); their defining matrices have entries in $\Aloc[z,\,1/(z-\zeta^{q^i}) \mid  i \geq 0]$ rather than in $\Aloc[z,\,1/(z-\zeta)]$; thus \cref{Thm:comparison} is not directly applicable to them. 
However, the same basic idea still applies, and we shall use it to interpret the quantities considered in \cite{Chang2025}.
We hope to give a formulation that covers such objects directly in future work.

Let $M_{\mathbf{C}_3}$ be the object of $\mathcal{M}$ defined by the square matrix
\begin{equation}
\Phi_{M_{\mathbf{C}_3}}
    :=\begin{pmatrix}
        0&1&0\\
        0 &0&1\\
        z-\zeta&0&0
    \end{pmatrix}
\end{equation}
of size $3$.
Then, its reduction $\overline{M}_{\mathbf{C}_3}$ is defined by the matrix 
\begin{equation}
    \Xi_{M_{\mathbf{C}_3}}:=\overline{\Phi}_{M_{\mathbf{C}_3}}=\begin{pmatrix}
        0&1&0\\
        0 &0&1\\
        z&0&0
    \end{pmatrix}.
\end{equation}
Let us consider the infinite product $\omega_{v,\,3}$ given by
\begin{equation}\label{eq:def of omega v. l}
    \omega_{v,\,3}=\omega_{v,\,3}(z):=\prod_{i=0}^{\infty} \left(
    1-\frac{\zeta^{q^{3i}}}{z}
    \right)^{-1}
\end{equation}
(see \cite{Chang2025}); then we can see that $\omega_{v,\,3}^{-1} \in \Acomp[[z,\,1/z\}$ and the matrix

\begin{equation}
    \Delta_{M_{\mathbf{C}_3}}=
    \begin{pmatrix}
        \big(\omega_{v,\,3}^{-1}\big)^{(2)}&0&0\\
        0&\big(\omega_{v,\,3}^{-1}\big)^{(1)}&0\\
       
        0&0&\big(\omega_{v,\,3}^{-1}\big)^{(0)}
    \end{pmatrix}
\end{equation}
 satisfies $\Delta_{M_{\mathbf{C}_3}} \Xi_{M_{\mathbf{C}_3}}=\Phi_{M_{\mathbf{C}_3}} \Delta_{M_{\mathbf{C}_3}}^{(1)}$ since we have
\begin{equation}\label{eq:omegavl difference eq}
    \omega_{v,\,3}(z)^{(3)}=\left(1-\frac{\zeta}{z}\right)\omega_{v,\,3}(z)^{(0)},
\end{equation}
which can be confirmed directly by the definition. 
As we can show $\Delta_{M_{\mathbf{C}_3}}\equiv I_3 \pmod{v}$ from the definition and the uniqueness of such a matrix is established in \cref{Thm:comparison}, the matrix $\Delta_{M_{\mathbf{C}_3}}$ is the one given by the theorem.
We consider the $\Aloc[z]$-basis $(m_0,\,m_1,\,m_{2})$ of $M_{\mathbf{C}_3}$ corresponding to the matrix $\Phi_{M_{\mathbf{C}_3}}$, and let $(\overline{m}_0,\,\overline{m}_1,\,\overline{m}_{2})$ be the basis of $\overline{M}_{\mathbf{C}_3}$ coming from $(m_0,\,m_1,\,m_{2})$ via the reduction. 
    
We consider a free $\Aloc[z]$-module $M_{\mathbf{C}_3}^{\prime}$ of rank $3$ equipped with the isomorphism
\begin{align}
    \tau_{M_{\mathbf{C}_3}^{\prime}}: \tau^* M_{\mathbf{C}_3}^{\prime}&\otimes_{\Aloc[z]}\Aloc[z,\,1/(z-\zeta^{q^i}) \mid  i \geq 0]\\
    &\rightarrow M_{\mathbf{C}_3}^{\prime}\otimes_{\Aloc[z]}\Aloc[z,\,1/(z-\zeta^{q^i}) \mid  i \geq 0]
\end{align}
represented by a matrix
    \begin{equation}
    \Phi_{M_{\mathbf{C}_3}^\prime}
    :=\begin{pmatrix}
        0&1&0\\
        0 &0&1\\
        \frac{1}{z-\zeta^{q}}&0&0
    \end{pmatrix}
\end{equation}
of size $3$ with respect to an $\Aloc[z]$-basis $(n_0,\,n_1,\,n_{2})$ of $M_{\mathbf{C}_3}^{\prime}$ and the basis $(n_0\otimes 1,\,n_1\otimes 1,\,n_{2}\otimes 1)$ of $ \tau^* M_{\mathbf{C}_3}^{\prime}$. 
Then, its reduction $\overline{M}_{\mathbf{C}_3}^{\prime}$ is the free $\quotfield[z]$-module of rank $3$ equipped with the isomorphism 
\begin{equation}
    \tau_{\overline{M}_{\mathbf{C}_3}^{\prime}}: \tau^* \overline{M}_{\mathbf{C}_3}^{\prime} \otimes_{\quotfield[z]}\quotfield[z,\,1/z] \rightarrow \overline{M}_{\mathbf{C}_3}^{\prime}\otimes_{\quotfield[z]}\quotfield[z,\,1/z] 
\end{equation}
defined by the matrix 
\begin{equation}
    \Xi_{M_{\mathbf{C}_3}^{\prime}}:=\overline{\Phi}_{M_{\mathbf{C}_3}^{\prime}}=\begin{pmatrix}
        0&1&0\\
        0 &0&1\\
        \frac{1}{z}&0&0
    \end{pmatrix};
    \end{equation}  
    we let $(\overline{n}_0,\,\overline{n}_1,\,\overline{n}_{2})$ be the basis of $\overline{M}_{\mathbf{C}_3}^\prime$ which is induced by $(n_0,\,n_1,\,n_{2})$ under reduction and corresponds to the matrix $\Xi_{M_{\mathbf{C}_3}^{\prime}}$.
    Although we cannot directly apply \cref{Thm:comparison} to $M_{\mathbf{C}_3}^\prime$, we can obtain $\Delta_{M_{\mathbf{C}_3}^{\prime}}\in \GL_3\big(\Acomp[[z,\,1/z\}[\omegav]\big)$ such that $\Delta_{M_{\mathbf{C}_3}^\prime} \Xi_{M_{\mathbf{C}_3}^\prime}=\Phi_{M_{\mathbf{C}_3}^\prime} \Delta_{M_{\mathbf{C}_3}^\prime}^{(1)}$ and $\Delta_{M_{\mathbf{C}_3}^\prime} \equiv I_3 \pmod{v}$.
     Indeed, the matrix 
    \begin{equation}
    \Delta_{M_{\mathbf{C}_3}^\prime}=
    \begin{pmatrix}
        \omega_{v,\,3}^{(3)}&0&0\\
        0&\omega_{v,\,3}^{(2)}&0\\    
        0&0&\omega_{v,\,3}^{(1)}
    \end{pmatrix}
\end{equation}
satisfies the desired conditions.
Here, we note that
\begin{equation}
    \omega_{v,\,3}=\omega_{v} \prod_{i=0}^{\infty} \left(
    1-\frac{\zeta^{q^{3i+1}}}{z}
    \right)\prod_{i=0}^{\infty} \left(
    1-\frac{\zeta^{q^{3i+2}}}{z}
    \right) \in \Acomp[[z,\,1/z\}[\omega_{v}].
\end{equation}

Let us consider the tensor product 
$
M_{\mathbf{C}_3}\otimes M_{\mathbf{C}_3}^\prime :=M_{\mathbf{C}_3}\otimes_{\Aloc[z]} M_{\mathbf{C}_3}^\prime
$
and the isomorphism 
\begin{align}
    \tau_{M_{\mathbf{C}_3}\otimes M_{\mathbf{C}_3}^\prime}:&\tau^*\left(M_{\mathbf{C}_3}\otimes M_{\mathbf{C}_3}^\prime\right)\otimes_{\Aloc[z]}\Aloc[z,\,1/(z-\zeta^{q^i}) \mid  i \geq 0] \\
    &\quad \rightarrow M_{\mathbf{C}_3}\otimes M_{\mathbf{C}_3}^\prime\otimes_{\Aloc[z]}\Aloc[z,\,1/(z-\zeta^{q^i}) \mid  i \geq 0]
\end{align} 
induced by the tensor product of the two structure isomorphisms.
With respect to the $\Aloc[z]$-basis 
\begin{equation}
(m_{0}\otimes n_{0}, \,  m_{0}\otimes n_{1},\, m_{0}\otimes n_{2},\, m_{1}\otimes n_{0}, \,  m_{1}\otimes n_{1},\, m_{1}\otimes n_{2},\, m_{2}\otimes n_{0}, \,  m_{2}\otimes n_{1},\, m_{2}\otimes n_{2})
\end{equation}
of $M_{\mathbf{C}_3}\otimes M_{\mathbf{C}_3}^\prime$ and the corresponding basis of $\tau^*(M_{\mathbf{C}_3}\otimes M_{\mathbf{C}_3}^\prime)=\tau^*(M_{\mathbf{C}_3})\otimes_{\Aloc[z]} \tau^*(M_{\mathbf{C}_3}^\prime)$, the isomorphism $\tau_{M_{\mathbf{C}_3}\otimes M_{\mathbf{C}_3}^\prime}$ is given by the Kronecker product
\begin{align}
   & \Phi_{M_{\mathbf{C}_3}\otimes M^\prime_{\mathbf{C}_3}}:=\Phi_{M_{\mathbf{C}_3}}\otimes \Phi_{M^\prime_{\mathbf{C}_3}}\\
    &\ =\begin{pmatrix}
    0&0&0&0&1&0&0&0&0\\
    0&0&0&0&0&1&0&0&0\\
    0&0&0&\frac{1}{z-\zeta^{q}}&0&0&0&0&0\\
    0&0&0&0&0&0&0&1&0\\
    0&0&0&0&0&0&0&0&1\\
    0&0&0&0&0&0&\frac{1}{z-\zeta^{q}}&0&0\\
    0&z-\zeta&0&0&0&0&0&0&0\\
    0&0&z-\zeta&0&0&0&0&0&0\\
     \frac{z-\zeta}{z-\zeta^{q}}&0&0&0&0&0&0&0&0
\end{pmatrix}.
\end{align}

We further consider the reduction
$
\overline{M}_{\mathbf{C}_3}\otimes \overline{M}_{\mathbf{C}_3}^\prime:=\overline{M}_{\mathbf{C}_3}\otimes_{\quotfield[z]} \overline{M}_{\mathbf{C}_3}^\prime
$
and the isomorphism 
\begin{align}
\tau_{\overline{M}_{\mathbf{C}_3}\otimes \overline{M}_{\mathbf{C}_3}^\prime}:&\tau^* \left(\overline{M}_{\mathbf{C}_3}\otimes \overline{M}_{\mathbf{C}_3}^\prime \right)\otimes_{\quotfield[z]}\quotfield[z,\,1/z] \\*
&\quad \rightarrow \overline{M}_{\mathbf{C}_3}\otimes \overline{M}_{\mathbf{C}_3}^\prime\otimes_{\quotfield[z]}\quotfield[z,\,1/z] 
\end{align}
given by 
\begin{equation}\label{eq:Xi_MC3MC3prime}
\begin{aligned}
    &\Xi_{M_{\mathbf{C}_3}\otimes M^\prime_{\mathbf{C}_3}}:=\Xi_{M_{\mathbf{C}_3}}\otimes \Xi_ {M^\prime_{\mathbf{C}_3}}\\
    \quad &=\begin{pmatrix}
    0&0&0&0&1&0&0&0&0\\
    0&0&0&0&0&1&0&0&0\\
    0&0&0&\frac{1}{z}&0&0&0&0&0\\
    0&0&0&0&0&0&0&1&0\\
    0&0&0&0&0&0&0&0&1\\
    0&0&0&0&0&0&\frac{1}{z}&0&0\\
    0&z&0&0&0&0&0&0&0\\
    0&0&z&0&0&0&0&0&0\\
    z\cdot \frac{1}{z}&0&0&0&0&0&0&0&0
\end{pmatrix}=\begin{pmatrix}
    0&0&0&0&1&0&0&0&0\\
    0&0&0&0&0&1&0&0&0\\
    0&0&0&\frac{1}{z}&0&0&0&0&0\\
    0&0&0&0&0&0&0&1&0\\
    0&0&0&0&0&0&0&0&1\\
    0&0&0&0&0&0&\frac{1}{z}&0&0\\
    0&z&0&0&0&0&0&0&0\\
    0&0&z&0&0&0&0&0&0\\
    1&0&0&0&0&0&0&0&0
\end{pmatrix}
\end{aligned}
\end{equation}
with respect to the basis
\begin{equation}
(\overline{m}_{0}\otimes \overline{n}_{0}, \, \overline{m}_{0}\otimes \overline{n}_{1},\,\overline{m}_{0}\otimes \overline{n}_{2},\,
\overline{m}_{1}\otimes \overline{n}_{0}, \, \overline{m}_{1}\otimes \overline{n}_{1},\,\overline{m}_{1}\otimes \overline{n}_{2},
\,
\overline{m}_{2}\otimes \overline{n}_{0}, \, \overline{m}_{2}\otimes \overline{n}_{1},\overline{m}_{2}\otimes \overline{n}_{2}).
\end{equation}
Even though we cannot apply \cref{Thm:comparison} to $M_{\mathbf{C}_3}\otimes M_{\mathbf{C}_3}^\prime$, it is possible to show, using \cref{eq:Kronecker_mixed,eq:Kronecker_twist}, that the isomorphism given by the Kronecker product $\Delta_{M_{\mathbf{C}_3}}\otimes \Delta_{ M_{\mathbf{C}_3}^\prime}$, which is the diagonal matrix of size $9$ whose diagonal entries are
\begin{align}
     \Big(&\big(\omega_{v,\,3}^{-1}\big)^{(2)} \big(\omega_{v,\,3}\big)^{(3)},\,\big(\omega_{v,\,3}^{-1}\big)^{(2)} \big(\omega_{v,\,3}\big)^{(2)},\, \big(\omega_{v,\,3}^{-1}\big)^{(2)} \big(\omega_{v,\,3}\big)^{(1)},\\
     &\quad
     \big(\omega_{v,\,3}^{-1}\big)^{(1)} \big(\omega_{v,\,3}\big)^{(3)}, \,\big(\omega_{v,\,3}^{-1}\big)^{(1)} \big(\omega_{v,\,3}\big)^{(2)},\, \big(\omega_{v,\,3}^{-1}\big)^{(1)} \big(\omega_{v,\,3}\big)^{(1)},\\
     &\quad \quad
     \big(\omega_{v,\,3}^{-1}\big)^{(0)} \big(\omega_{v,\,3}\big)^{(3)},\,\big(\omega_{v,\,3}^{-1}\big)^{(0)} \big(\omega_{v,\,3}\big)^{(2)},\, \big(\omega_{v,\,3}^{-1}\big)^{(0)} \big(\omega_{v,\,3}\big)^{(1)} \Big),
\end{align}
has the properties stated in \cref{Thm:comparison}.
By \cref{eq:Xi_MC3MC3prime}, the $\quotfield[z]$-linear morphism $\iota:\overline{\mathbf{1}}\rightarrow \overline{M}_{\mathbf{C}_3}\otimes \overline{M}_{\mathbf{C}_3}^\prime$ represented by
\begin{equation}
    \nu_{\iota}:=(1,\,0,\,0,\,0,\,1,\,0,\,0,\,0,\,1)^\top
\end{equation}
is compatible with the isomorphism $\tau_{\overline{M}_{\mathbf{C}_3}\otimes \overline{M}_{\mathbf{C}_3}^\prime}$; that is, we have $\Xi_{M_{\mathbf{C}_3}\otimes M^\prime_{\mathbf{C}_3}}\cdot \nu_{\iota}^{(1)}=\Xi_{M_{\mathbf{C}_3}\otimes M^\prime_{\mathbf{C}_3}}\cdot \nu_{\iota}= \nu_{\iota} \cdot I_1$.
Hence, although $M_{\mathbf{C}_3}\otimes M_{\mathbf{C}_3}^\prime$ is not an object of $\mathcal{M}$, we obtain the vector
\begin{equation}
    \Delta_{M_{\mathbf{C}_3}\otimes M_{\mathbf{C}_3}^\prime}^{(1)}\cdot \nu_{\iota}=\left( \frac{\omega_{v,\,3}^{(4)}}{\omega_{v,\,3}^{(3)}},\,0,\,0,\,0,\,\frac{\omega_{v,\,3}^{(3)}}{\omega_{v,\,3}^{(2)}},\,0,\,0,\,0,\,\frac{\omega_{v,\,3}^{(2)}}{\omega_{v,\,3}^{(1)}}\right)^\top
\end{equation}
and the quantities 
\begin{equation}\label{eq:CWYperiods}
    \frac{\omega_{v,\,3}^{(4-i)}}{\omega_{v,\,3}^{(3-i)}}\Big|_{z=\zeta}, \quad (i=0,\,1,\,2)
\end{equation}
by a construction similar to those of pre-period vectors and Andr\'{e}'s $v$-adic periods.
The quantities in \eqref{eq:CWYperiods} agree, up to factors in $K^\times$, with the nonzero entries of the matrix representing the inverse of the Frobenius action considered in \cite[(3.7.1)--(3.7.3)]{Chang2025}.

\end{Remark}

\section{Chen's $v$-adic version of the ABP criterion}\label{Section:v-adic_ABP}

In this section, we apply Chen's $v$-adic version of the ABP criterion (\cref{Thm:YTcriterion}) to study the transcendence degrees of $v$-adic periods.
Roughly speaking, his criterion lifts $\Kalg$-linear relations among the special values of entire series satisfying a Frobenius difference equation to $\Kalg[1/z]$-linear relations among the series themselves; see \cref{Thm:YTcriterion} or \cite[Theorem 2.2.2]{Chen2025} for the precise statement.
After recalling it, we apply the criterion to suitable auxiliary systems, so as to detect not only the linear relations but also the algebraic ones.
This yields an equality between the transcendence degree of such series and that of their special values; the argument parallels those in \cite[\S 4.3]{Anderson2004} and \cite[\S 5.2]{Papanikolas2008}.
We further observe that the formulation of pre-period vectors (\cref{Def:pre-period}) yields the difference equation required to apply the $v$-adic version of the ABP criterion.
Applying these results to the objects introduced in \cref{Example:v-adicpolyloginterpretation}, we relate the transcendence degree of the special values of $v$-adic Carlitz multiple polylogarithms to that of the corresponding pre-periods.

\subsection{The criterion and consequences for algebraic independence}\label{subsection:Criterion_algebraic}
In this subsection, we recall the $v$-adic version of the ABP criterion established by Chen and its consequence for algebraic independence.
Throughout this subsection, $v$ may be a finite place of any degree and $z$ is regarded as an independent variable over
$\Kalg$; specializing $z$ to $\zeta$ means applying the $\Kalg$-algebra homomorphism $\Kalg[1/z]\to\Kalg$,
$1/z\mapsto1/\zeta$, and its extensions to the rings of convergent series
considered here.

The formulation below is adapted to the notation and the choice of variable used in the present paper. More precisely, it is obtained from \cite[Theorem 2.2.2]{Chen2025} by using $1/z$ in place of the variable $t$ appearing there; thus the following statement is a special case of Chen’s theorem and does not differ from it in any essential way.

\begin{Theorem}[{\cite[Theorem 2.2.2]{Chen2025}}]\label{Thm:YTcriterion}
Take a square matrix $\Phi \in \Mat_{r\times r}(\Kalg[1/z])$, an element $\gamma \in (\Kalg)^\times\setminus (\mathbb{F}_q^{\textrm{alg}})^\times$, and a column vector $\psi \in \Mat_{r\times 1}(\Entirev)$ such that $\psi^{(-1)}=\Phi \psi$. 
If $\det \Phi|_{1/z=\gamma^{(-i)}}\neq 0$ for any $i \geq 1$, then, for any $\ell \in \Mat_{1\times r}(\Kalg)$ such that $\ell \cdot(\psi|_{1/z=\gamma})=0$, there exists $L \in \Mat_{1\times r}(\Kalg[1/z]) $ such that $L|_{1/z=\gamma}=\ell$ and $L \cdot \psi=0$.
\end{Theorem}

We note that Chen established the $v$-adic ABP criterion for finite places $v$ of any degree.
Recalling here that $\Omegav =\big(\omegav^{(1)}\big)^{-1}\in \Entirev$ and $\Omegav^{(-1)}=\left(1-\zeta/z\right)\Omegav$,
we use \cref{Thm:YTcriterion} in the following form:

\begin{Corollary}\label{Cor:ABPapplication}
Take a square matrix $\Phi \in \Mat_{r\times r}(\Kalg[1/z,\,1/(z-\zeta)])$, an element $h$ of $\mathbb{F}_q(z)^\times$, and a vector $\psi$ in $\Mat_{r\times 1}\big(\Acomp[[z,\,1/z\}[\omegav^{(1)}]\big)$ such that $h \psi $ is in $\Mat_{r\times 1}(\Entirev[1/\Omegav])$ and $\psi^{(-1)}=\Phi \psi$.
Assume that $\det \Phi  |_{z=\zeta^{(-i)}}\neq0$ for any $i\geq 1$; then, for any row vector $\ell \in \Mat_{1\times r}(\Kalg)$ such that $\ell \cdot (\psi|_{z=\zeta})=0$, we have $L \in \Mat_{1\times r}(\Kalg[1/z])$ satisfying $L|_{z=\zeta}=\ell$ and $L \cdot \psi=0$.
\end{Corollary}
\begin{proof}
    As we have
    \begin{equation}
        \frac{1}{z-\zeta}\left(1-\frac{\zeta}{z}\right)=\frac{1}{z},
    \end{equation}
    for sufficiently large $m \geq 0$, it holds that \begin{equation}
    \left(1-\frac{\zeta}{z}\right)^m\Phi \in \operatorname{Mat}_{r\times r}(\Kalg[1/z]) \text{ and }
     %\end{equation} and \begin{equation}
    \Omegav^m h\psi \in \operatorname{Mat}_{r\times 1}(\Entirev).
    \end{equation}
    We note that $h \in\mathbb{F}_q(z)^\times$ is regular at $z=\zeta$ as $\zeta$ is transcendental over $\mathbb{F}_q$, and $\Omegav$ is regular and non-zero at $z=\zeta$; now, we have 
    \begin{equation}     
     \ell \cdot \big(\left(\Omegav^m\cdot h \psi\right)|_{z=\zeta}\big)=0, \quad  \det \left(\left(1-\frac{\zeta}{z}\right)^m\Phi\middle)\right|_{z=\zeta^{(-i)}}\neq 0 \text{ for $i \geq 1$},
    \end{equation}
     and 
    \begin{equation}
     \left(\Omegav^m h \psi\right)^{(-1)}=\left(\left(1-\frac{\zeta}{z}\right)^m\Phi\right)\cdot \left(\Omegav^m\cdot h\psi \right)
    \end{equation}
    since we have $\Omegav^{(-1)}=\left(1-\zeta/z\right)\Omegav$ and $h^{(-1)}=h$.
    Hence, \cref{Thm:YTcriterion} yields a vector $L\in \Mat_{1\times r}(\Kalg[1/z])$ such that $L|_{z=\zeta}=\ell$ and $L \cdot \big(\Omegav^m h \psi \big)=0$, from which it follows that $L \cdot \psi=0$.
\end{proof}

%By this proposition, the dimensions of 
\begin{Corollary}\label{cor:dimensioncoincide}
    In the situation of \cref{Cor:ABPapplication}, we have
    \begin{equation}
        \dim_{\Kalg\left(\frac{1}{z}\right)}\left\langle\psi_1,\,\dots,\,\psi_r\right\rangle_{\Kalg\left(\frac{1}{z}\right)}=\dim_{\Kalg}\left\langle \psi_1|_{z=\zeta},\,\dots,\,\psi_r|_{z=\zeta}\right\rangle_{\Kalg}.
    \end{equation}
    Here, the vector spaces appearing on the left- and right-hand sides are subspaces of $\Frac\mathbb{T}_v$ and $\mathbb{C}_v$, respectively.
\end{Corollary}

\begin{proof}
    We can obtain this corollary by the same arguments as in \cite[Proposition 4.4.3]{Anderson2004}; we reproduce the argument here to make the present paper self-contained.
    Letting $H_1$ be the torsion-free $\Kalg[1/z]$-module $\{L \in \Mat_{1\times r}(\Kalg[1/z])\mid L \cdot \psi=0\}$, we have an exact sequence
    \begin{equation}
        0\hookrightarrow H_1 \rightarrow \Mat_{1\times r}(\Kalg[1/z])\twoheadrightarrow H_0 \rightarrow 0,
    \end{equation}
    of $\Kalg[1/z]$-modules, where $H_0$ is the submodule of $\Frac\mathbb{T}_v$ spanned by the set $\{\psi_1,\dots,\psi_{r}\}$ and the surjection is given by $L \mapsto L\cdot \psi$. 
    As $\Kalg[1/z]$ is a PID and $H_0 \subset \Frac\mathbb{T}_v $ is also torsion-free, both $H_0$ and $H_1$ are free and we have $H_0\oplus H_1 \simeq \Mat_{1\times r}(\Kalg[1/z])$. 
    Since $H_0$ is free, the above exact sequence splits and therefore remains exact after tensoring with $\Kalg$ via the specialization $1/z\mapsto1/\zeta$.
    The dimension over $\Kalg$ of the base change $H_1\otimes_{\Kalg[1/z]}\Kalg$ with respect to the $\Kalg$-algebra homomorphism $\Kalg[1/z] \rightarrow \Kalg$ given by $1/z\mapsto1/\zeta$ is equal to the rank of $H_1$ over $\Kalg[1/z]$ as $H_1$ is free. 
    If we regard $H_1\otimes_{\Kalg[1/z]}\Kalg$ as a $\Kalg$-vector subspace of 
    \begin{equation}
        \Mat_{1\times r}(\Kalg)=\Mat_{1\times r}(\Kalg[1/z])\otimes_{\Kalg[1/z]} \Kalg,
    \end{equation}
    it is contained in $V_1:=\{\ell \in \Mat_{1\times r}(\Kalg)\mid \ell \cdot(\psi|_{z=\zeta})=0\}$; and \cref{Cor:ABPapplication} gives the opposite inclusion $V_1 \subset H_1\otimes_{\Kalg[1/z]}\Kalg$. Therefore, we have
    \begin{align}
        &\dim_{\Kalg\left(\frac{1}{z}\right)}\left\langle\psi_1,\,\dots,\,\psi_r\right\rangle_{\Kalg\left(\frac{1}{z}\right)}=\rank_{\Kalg[1/z]}H_0\\=&r-\rank_{\Kalg[1/z]}H_1=r-\dim_{\Kalg}V_1=\dim_{\Kalg}V_0,
    \end{align}
    where $V_0$ stands for the $\Kalg$-vector subspace of $\mathbb{C}_v$ spanned by $\psi_1|_{z=\zeta},\,\dots,\,\psi_r|_{z=\zeta}$ since we have an exact sequence $0 \rightarrow V_1 \rightarrow \Mat_{1\times r}(\Kalg) \rightarrow V_0\rightarrow0$.
\end{proof}

Using this corollary, we obtain the following theorem on the algebraic independence of $v$-adic periods:
\begin{Theorem}\label{Thm:period_vector_algebraic_independence}
    Take a square matrix $\Phi \in \Mat_{r\times r}(\Kalg[1/z,\,1/(z-\zeta)])$, an element $h$ of $\mathbb{F}_q(z)^\times$, and a vector $\psi=(\psi_1,\,\dots,\,\psi_r)^{\top}$ in $\Mat_{r\times 1}\big( \Acomp[[z,\,1/z\}[\omegav^{(1)}] \big)$ such that $h \psi \in \Mat_{r\times 1}(\Entirev[1/\Omegav])$.
If $\det \Phi  |_{z=\zeta^{(-i)}}\neq0$ for any $i\geq 1$ and $\psi^{(-1)}=\Phi \psi$, then we have
    \begin{equation}
        \trdeg_{\Kalg\left(\frac{1}{z}\right)}\Kalg\left(\frac{1}{z},\,\psi_1,\,\dots,\,\psi_r\right)=\trdeg_{\Kalg}\Kalg\left( \psi_1|_{z=\zeta},\,\dots,\,\psi_r|_{z=\zeta}\right)
    \end{equation}
    where the field on the left-hand side is a subfield of $\Frac\mathbb{T}_v$ and the field on the right-hand side is a subfield of $\mathbb{C}_v$. 
\end{Theorem}

%We use the notion of Kronecker products of matrices in the proof.
\begin{proof}
    Although we can obtain the theorem by an argument similar to that in \cite[Theorem 5.2.2]{Papanikolas2008}, we reproduce the proof here to make the present paper self-contained.
    We let $S$ be the ring $\Kalg(1/z)[\psi_1,\,\dots,\,\psi_r]$ and consider its filtration $F_0S\subset F_1S \subset\cdots$ given by 
    \begin{equation}
        F_nS:=\left\langle 
            \psi_{i_1}\cdots \psi_{i_m}
            \ \middle|\ 
            0\leq m \leq n, 1\leq i_1,\,\dots,\,i_m\leq r
        \right\rangle_{\Kalg(1/z)},
    \end{equation}
    that is, $F_nS$ is the $\Kalg(1/z)$-linear subspace of $S$ spanned by products of at most $n$ entries of the vector $\psi$.
    Further, we write $Q$ for the ring $\Kalg[\psi_1|_{z=\zeta},\,\dots,\,\psi_r|_{z=\zeta}]$ and consider its filtration $F_0Q\subset F_1Q \subset\cdots$ given by 
    \begin{equation}
    F_nQ:=\left\langle 
            \psi_{i_1}|_{z=\zeta}\cdots \psi_{i_m}|_{z=\zeta}
            \ \middle|\ 
            0\leq m \leq n, 1\leq i_1,\,\dots,\,i_m\leq r
        \right\rangle_{\Kalg}.
    \end{equation}
   For each $n\geq 0$, let us consider the matrix $\Phi_n^\prime:=(1)\oplus \Phi\oplus \Phi^{\otimes2}\oplus\cdots\oplus \Phi^{\otimes n}$ and the vector $\psi_n^\prime:=(1)\oplus (h \psi)\oplus (h \psi)^{\otimes2}\oplus\cdots\oplus (h \psi)^{\otimes n}$; then we have $\left(\psi_n^\prime\right)^{(-1)}=\Phi_n^\prime \psi_n^\prime $ by Equations \eqref{eq:Kronecker_mixed}, \eqref{eq:Kronecker_twist}, and the fact that $h \in \mathbb{F}_q(z)$; further, the determinant of $ \Phi_n^\prime$ is a positive power of $\det \Phi$ by \cref{eq:Kronecker_det}; thus we have $\det \Phi_n^\prime|_{z=\zeta^{(-i)}} \neq 0$ for $i \geq 1 $ and $n \geq 1$. 
    By definition, the entries of $\psi_n^\prime$ coincide with the generators spanning $F_nS$ up to multiplication by powers of $h \in \mathbb{F}_q(z)^\times$ and hence we have $\dim_{\Kalg(1/z)}F_nS=\dim_{\Kalg}F_nQ$ by \cref{cor:dimensioncoincide}; we can apply the corollary since the entries of $\psi_n^\prime$ are in $\Entirev[1/\Omegav]$.

    To deduce the desired equality of transcendence degrees from these dimension equalities, we work directly with degree filtrations without passing to the associated graded rings, following the conventions of \cite[\S 11.1]{Kemper2011}, although the underlying argument is essentially the same as the standard one from dimension theory.
    Now let us write $Q \simeq \Kalg[x_1,\,\dots,\,x_r]/\mathfrak{a}$, where $\mathfrak{a}$ is an ideal of the polynomial ring and $x_i$ corresponds to $\psi_i|_{z=\zeta}$, and write
    $S\simeq \Kalg(1/z)[x_1,\,\dots,\,x_r]/\mathfrak{b}$ similarly; then $\dim_{\Kalg(1/z)}F_nS=\dim_{\Kalg}F_nQ$ implies that the Hilbert series of the ideals $\mathfrak{a}$ and $\mathfrak{b}$ coincide.
    Therefore, the Krull dimension of $Q$ is equal to that of $S$ (\cite[Theorem 11.13]{Kemper2011}), and hence the equality of the transcendence degrees follows.
\end{proof}

\subsection{Conditions for applying the criterion}
In this subsection, we consider the conditions needed for applying \cref{Thm:YTcriterion} and its corollaries.
After observing that the difference equation required in \cref{Thm:YTcriterion,Thm:period_vector_algebraic_independence} automatically comes from the formulation of pre-period vectors (\cref{Def:pre-period}), we apply \cref{Thm:period_vector_algebraic_independence} to pre-periods in \cref{Example:v-adicpolyloginterpretation} (\cref{cor:vprelogpolylog}); this is a crucial step in the proof of \cref{Thm:main_CMPLs}.
Then, we will deal with a $v$-adic counterpart of Proposition 3.1.3 in \cite{Anderson2004}, which makes it easy to check a condition required to apply \cref{Thm:YTcriterion,Thm:period_vector_algebraic_independence}.

We shall check that the formulation of pre-period vectors (\cref{Def:pre-period}) naturally implies that the vectors satisfy the difference equation required in \cref{Thm:YTcriterion,Thm:period_vector_algebraic_independence}.

\begin{Lemma}\label{Lem:psi difference equation}
Assume $\deg v=1$, take an object $M$ of $\mathcal{M}$ defined by a matrix $\Phi_M \in \GL_r(\Aloc[z,$ $1/(z-\zeta)])$ and let $\overline{M}\in \overline{\mathcal{M}}$ be its reduction defined by $\Xi_M \in \GL_r(\quotfield[z,\,1/z])$.
Take a morphism $\iota:\overline{\mathbf{1}}\rightarrow \overline{M}$, which is defined by a vector $\nu_\iota \in \Mat_{r\times 1}(\quotfield[z])$, and let $\Delta_M$ be the matrix given in \Cref{Thm:comparison}.
Then, the pre-period vector $\psi=\psi_{M,\,\iota}:= \Delta_M^{(1)} \cdot \nu_{\iota}^{(1)}$ satisfies the difference equation 
        \begin{equation}
            \psi^{(-1)}=\Phi_M \cdot \psi \in \Mat_{r\times 1} (\Acomp[[z,\,1/z\}[\omegav]).
        \end{equation}
\end{Lemma}
 
\begin{proof}
    The definition of $\iota$ shows that the right square of the diagram
    \begin{equation}\label{eq:iotacommutativediagram}
        \begin{tikzpicture}[auto]
        \node (0-3) at (-4.1, 2) {$\overline{\mathbf{1}}\otimes_{\quotfield[z]}\quotfield[z,\,1/z]$}; \node (03) at (0, 2) {$\tau^{*}\overline{\mathbf{1}}\otimes_{\quotfield[z]}\quotfield[z,\,1/z] $}; \node (33) at (4, 2) {$\overline{\mathbf{1}}\otimes_{\quotfield[z]}\quotfield[z,\,1/z]$}; 
        \node (0-0) at (-4.1, 0) {$\overline{M}\otimes_{\quotfield[z]}\quotfield[z,\,1/z]$}; \node (00) at (0, 0) {$\tau^{*}\overline{M}\otimes_{\quotfield[z]}\quotfield[z,\,1/z]$}; \node (30) at (4, 0) {$\overline{M}\otimes_{\quotfield[z]}\quotfield[z,\,1/z]$}; 
        \draw[->] (03) to node {$\tau$}(33);
        \draw[->] (00) to node {$\tau$}(30);
        \draw[->] (03) to node {$\tau^{*}\iota\otimes_{\quotfield[z]}\quotfield[z,\,1/z]$} (00);
        \draw[->] (33) to node {$\iota\otimes_{\quotfield[z]}\quotfield[z,\,1/z]$} (30);
         \draw[->] (0-3) to node {$\iota\otimes_{\quotfield[z]}\quotfield[z,\,1/z]$} (0-0);
         \draw[-] (0-3) to node {$=$} (03);
         \draw[-] (0-0) to node {$=$} (00);
        \end{tikzpicture}
    \end{equation}
    is commutative; since $\tau$ acts trivially on $\mathbb{F}_v[z]$, the base change along $\tau$, denoted by $\tau^*$, has no effect, and hence the left square is also commutative.
    In terms of the representing matrix $\nu_\iota$, it holds that
    $
           \Xi_M \cdot \nu_\iota^{(1)}=\nu_\iota\cdot I_1=\nu_\iota
    $. Thus, we have
    \begin{equation}
        \psi^{(-1)}
        =\left(\Delta_M^{(1)}\cdot \nu_{\iota}^{(1)} \right)^{(-1)}
        =\Delta_M\cdot \nu_{\iota}
        =\Delta_M\cdot\Xi_M\cdot \nu_{\iota}^{(1)}
        =\Phi_M \cdot \Delta_M^{(1)}\cdot  \nu_{\iota}^{(1)}
        =\Phi_M \cdot \psi
    \end{equation}
    by \cref{eq:Deltaproperty1}; hence the result follows.
\end{proof}

\begin{Example}
    We consider the case discussed in \cref{Example:v-adicpolyloginterpretation,Example:v-adicpolylogcomparison}.
    Then we have a difference equation $\psi_{\tilde{\mathbf{a}}}^{(-1)}=\Phi_{\tilde{\mathbf{a}}}\psi_{\tilde{\mathbf{a}}}$.
\end{Example}

\begin{Example}\label{Example:divergent_points}
    In the situation of \cref{Example:preperiods-non-conv}, we have a pre-period vector
    \begin{align}
        \psi_{\iota}&:=
        \begin{pmatrix}
            z^s-1\\
           \displaystyle{ \lim_{m\to \infty}
        \left(
        \left(z^s-1\right)
        \prelog_{\tilde{\mathbf{a}},\,m}+\overline{a}\delta_{\mathcal{C},\,m}^{-s}z^{-(m+1)s}
        \right)}
        \end{pmatrix}^{(1)}\\
        &=
        \begin{pmatrix}
            z^s-1\\
           \displaystyle{ \lim_{m\to \infty}
        \left(
        \left(z^s-1\right)
        \prelog_{\tilde{\mathbf{a}},\,m}^{(1)}+\frac{\overline{a}}{z^{(m+1)s}}\left(\delta_{\mathcal{C},\,m}^{(1)}\right)^{-s}
        \right)}
        \end{pmatrix};
    \end{align}
    \cref{Lem:psi difference equation} yields the following difference equation:
    \begin{align}
        &\quad 
           \lim_{m\to \infty}
        \left(
        \left(z^s-1\right)
        \prelog_{\tilde{\mathbf{a}},\,m}^{(1)}+\frac{\overline{a}}{z^{(m+1)s}}\left(\delta_{\mathcal{C},\,m}^{(1)}\right)^{-s}
        \right)^{(-1)}\\
        &= \frac{a(z^s-1)}{(z-\zeta)^s}+\frac{1}{(z-\zeta)^s}\lim_{m\to \infty}
        \left(
        \left(z^s-1\right)
        \prelog_{\tilde{\mathbf{a}},\,m}^{(1)}+\frac{\overline{a}}{z^{(m+1)s}}\left(\delta_{\mathcal{C},\,m}^{(1)}\right)^{-s}
        \right).
    \end{align}
    
\end{Example}

\begin{Corollary}\label{Cor:v-adic_periods_to_ABP}
Let $v$ be a place of degree $1$, $M$ be an object of $\mathcal{M}$ defined by $\Phi_M \in \Mat_r(\Aloc[1/(z-\zeta)])\cap \GL_r(\Aloc[z,\,1/(z-\zeta)]) $, and $\overline{M}\in \overline{\mathcal{M}}$ be its reduction defined by $\Xi_M $. %\in \Mat_r(\quotfield[1/z])$.
Take a morphism $\iota:\overline{\mathbf{1}}\rightarrow \overline{M}$, which is defined by $\nu_\iota \in \Mat_{r\times 1}(\quotfield[z])$, and let $\Delta_M$ be the matrix given in \Cref{Thm:comparison}.
Consider the pre-period vector $\psi=\psi_{M,\,\iota}:= \Delta_M^{(1)} \cdot \nu_{\iota}^{(1)}$; if there exists $h\in\mathbb{F}_q(z)^\times$ such that $h \psi \in \Mat_{r\times 1}(\Entirev[1/\Omegav])$, then we have 
\begin{equation}
        \trdeg_{\Kalg\left(\frac{1}{z}\right)}\Kalg\left(\frac{1}{z},\,\psi_1,\,\dots,\,\psi_r\right)=\trdeg_{\Kalg}\Kalg\left( \psi_1|_{z=\zeta},\,\dots,\,\psi_r|_{z=\zeta}\right)
    \end{equation}
    where we write $\psi=(\psi_1,\,\dots,\,\psi_r)^\top$.
\end{Corollary}

\begin{proof} 
    As $\Phi_M$ is in $\GL_r(\Aloc[z,\,1/(z-\zeta)])$, its determinant is a unit of $\Aloc[z,\,1/(z-\zeta)]$; 
    thus the determinant is of the form $c(z-\zeta)^n$ with $c \in \Aloc^\times$ and $n \in \mathbb{Z}$. Hence we have $\det \Phi_M|_{z=\zeta^{(-i)}}=c\,(\zeta^{q^{-i}}-\zeta)^n\neq0$ for any $i\geq 1$. 
    Moreover, the entries of $\Phi_M$ lie in $\Aloc[1/(z-\zeta)]\subset \Kalg[1/z,\,1/(z-\zeta)]$ by assumption. 
    Since \cref{Lem:psi difference equation} shows that we have $\psi^{(-1)}=\Phi_M \psi$, we can use \cref{Thm:period_vector_algebraic_independence} to obtain the desired result.
\end{proof}

\begin{Corollary}\label{cor:vprelogpolylog}
    Take a finite place $v$ of degree $1$ and a finite subset $I$ of $\bigcup_{r=1}^\infty\tilde{\mathbb{D}}_v^r$.
    Suppose that $\tilde{\mathbf{a}}=\big((s_1,\,a_1),$ $\dots,\,(s_r,\,a_r)\big)\in I$ implies $\big((s_1,\,a_1),$ $\dots,$ $(s_{r-1},\,a_{r-1})\big)\in I$ if $r \geq 2$.
    Then it follows that
    \begin{equation}
        \trdeg_{\Kalg(z)}\Kalg \left(z,\,\prelog_{v,\,\tilde{\mathbf{a}}}^{(1)}\,\middle|\,\tilde{\mathbf{a}} \in I \right)=\trdeg_{\Kalg}\Kalg \left(\mathrm{Li}(\tilde{\mathbf{a}})_{v} \,\middle|\,\tilde{\mathbf{a}} \in I \right).
    \end{equation}
\end{Corollary}

\begin{proof}
    We use the notation in \cref{Example:v-adicpolylogcomparison,Example:v-adicpolyloginterpretation} in this proof. Let us consider the object 
    \begin{equation}
        M:= \bigoplus_{\tilde{\mathbf{a}}\in I}M_{\tilde{\mathbf{a}}}
    \end{equation}
    of $\mathcal{M}$ whose structure isomorphism is represented by the matrix
    \begin{equation}
        \Phi= \bigoplus_{\tilde{\mathbf{a}}\in I}\Phi_{\tilde{\mathbf{a}}}
    \end{equation}
    where the components $\Phi_{\tilde{\mathbf{a}}}$ are those defined in \cref{eq:def matrix T}. 
    Here, since the degree of $v$ is $1$, we have $t-\theta=z-\zeta$ and hence the entries of $\Phi$ are in $\Aloc[1/(z-\zeta)]$.
    The modulo $v$ reduction 
    \begin{equation}
        \overline{M}:=\bigoplus_{\tilde{\mathbf{a}}\in I}\overline{M}_{\tilde{\mathbf{a}}}
    \end{equation}
    admits a morphism 
    $
        \iota:= (\iota_{\tilde{\mathbf{a}}})_{\tilde{\mathbf{a}}\in I}.
    $
   This is represented by the column matrix
    \begin{equation}
        \nu:= \bigoplus_{\tilde{\mathbf{a}}\in I}\nu_{\iota_{\tilde{\mathbf{a}}}},
    \end{equation}
    and the corresponding pre-period vector $\Delta_M^{(1)}\cdot \nu^{(1)}$ is given by
    \begin{equation}
        \psi:= \bigoplus_{\tilde{\mathbf{a}}\in I}\psi_{\tilde{\mathbf{a}}}
    \end{equation}
    (see \cref{eq:v-adic polylog pre-period vector}), where $\Delta_M$ is the square matrix given by \cref{Thm:comparison}.
    Here, for column matrices $\mathbf{a}=(a_1,\,\dots,\,a_r)^\top$ and $\mathbf{b}=(b_1,\,\dots,$ $b_s)^\top$, we denote the column matrix $(a_1,\,\dots,\,a_r,\,b_1,\,\dots,\,b_s)^\top$ by $\mathbf{a}\oplus \mathbf{b}$.    
    By virtue of \cref{Prop:v-adic polylog entire},
    all entries of $\psi$ are in $\Entirev[1/\Omegav]$; thus we can use \cref{Cor:v-adic_periods_to_ABP} and \cref{eq:v_prelog_and_log} to obtain 
    \begin{align}
        \trdeg_{\Kalg(z)}\Kalg(z) \left(\prelog_{v,\,\tilde{\mathbf{a}}}^{(1)}\,\middle|\,\tilde{\mathbf{a}} \in I \right)&=
        \trdeg_{\Kalg(1/z)}\Kalg \left(\frac{1}{z}\right)\left(\prelog_{v,\,\tilde{\mathbf{a}}}^{(1)}\,\middle|\,\tilde{\mathbf{a}} \in I \right)\\
        &=\trdeg_{\Kalg}\Kalg \left(\mathrm{Li}(\tilde{\mathbf{a}})_{v} \,\middle|\,\tilde{\mathbf{a}} \in I \right),
    \end{align}
    where we use
    \begin{equation}
        \Kalg(z,\,\prelog_{v,\,\tilde{\mathbf{a}}}^{(1)}\mid\tilde{\mathbf{a}} \in I) =\Kalg(z) (\text{entries of $\psi$})
    \end{equation}
    which follows from the assumption on the set $I$.
\end{proof}

The rest of this section is devoted to \cref{Prop:ABP_condition}, a $v$-adic counterpart of \cite[Proposition 3.1.3]{Anderson2004}.
This proposition is not needed for the proof of our main result \cref{Thm:main_CMPLs}.
For the rest of this subsection, the place $v$ can again be of any degree, and $z$ is regarded as an independent variable over $\Kalg$. 
In addition to moving from the $\infty$-adic setting to the $v$-adic setting, we shall weaken a hypothesis for applying the proposition. To weaken the hypothesis, the following lemma is needed:

\begin{Lemma}\label{Lemma:ABPCondition}
    Given $\mathbf{d}\in \Mat_{r\times 1}(\Kalg)$ and $B \in \Mat_{r\times r}(\Kalg)$, the set
    \begin{equation}
        A_{\mathbf{d},\,B}:=\{\mathbf{x}=(x_1,\,\dots,\,x_r)^\top \in \Mat_{r\times 1}(\mathbb{C}_v) \mid \mathbf{x}-B \cdot \mathbf{x}^{(1)}=\mathbf{d} \}
    \end{equation}
    is contained in $\Mat_{r\times 1}(\Kalg)$.
\end{Lemma}
\begin{proof}
    We may assume that the set $A_{\mathbf{d},\,B}$ is non-empty.
    In the case where $\mathbf{d}$ is the zero vector $\mathbf{0}$, the set $A_{\mathbf{0},\,B}$ is an $\mathbb{F}_q$-linear subspace of $\Mat_{r\times 1}(\mathbb{C}_v)$.
    We claim that any set of $\mathbb{C}_v$-linearly dependent elements of $A_{\mathbf{0},\,B}$ is $\mathbb{F}_q$-linearly dependent; this follows from a well-known argument; see, for example, \cite[Lemma 3.3.7]{Papanikolas2008}.
    Suppose on the contrary that there exist $l \geq 1$ and a $\mathbb{C}_v$-linear relation
    \begin{equation}\label{eq:c_linear}
        \mathbf{x}_0+c_1 \mathbf{x}_1+\cdots+c_l \mathbf{x}_l=0, \quad c_1,\,\dots,\,c_l \in \mathbb{C}_v,
    \end{equation}
    among elements $\mathbf{x}_0,\,\dots,\,\mathbf{x}_{l}$ of $A_{\mathbf{0},\,B}$ which are $\mathbb{F}_q$-linearly independent; we can further assume that $l$ is taken minimally.
    Then, \cref{eq:c_linear} implies 
    $B \mathbf{x}_0^{(1)}+c_1^{(1)} B\mathbf{x}_1^{(1)}+\cdots+c_l^{(1)} B\mathbf{x}_l^{(1)}=B(\mathbf{x}_0^{(1)}+c_1^{(1)} \mathbf{x}_1^{(1)}+\cdots+c_l^{(1)} \mathbf{x}_l^{(1)})=0$; as $\mathbf{x}_0,\,\dots,\,\mathbf{x}_l$ are elements of $A_{\mathbf{0},\,B}$, we further have $\mathbf{x}_0+c_1^{(1)} \mathbf{x}_1+\cdots+c_l^{(1)} \mathbf{x}_l=0$. 
    Subtracting this equation from \cref{eq:c_linear}, we have $(c_1-c_1^{(1)}) \mathbf{x}_1+\cdots+(c_l-c_l^{(1)} )\mathbf{x}_l=0$; this linear relation must be trivial by the minimality of $l$, that is, we have $c_1-c_1^{(1)}=\cdots=c_l-c_l^{(1)}=0$. 
    Thus, $c_1,\dots,\,c_l$ are elements of $\mathbb{F}_q$, the fixed part of $\mathbb{C}_v$ under twist; this contradicts the $\mathbb{F}_q$-linear independence of $\mathbf{x}_0,\,\dots,\,\mathbf{x}_{l}$. 
    Hence the claim follows, and it shows that the rank of $A_{\mathbf{0},\,B}$ over $\mathbb{F}_q$ is at most $r$, the rank of $\Mat_{r\times 1}(\mathbb{C}_v)$ over $\mathbb{C}_v$; thus $A_{\mathbf{0},\,B}$ is a finite set.
    
   For arbitrary $\mathbf{d} \in \Mat_{r\times 1}(\Kalg)$ and any $\mathbf{y} \in A_{\mathbf{d},\,B}$, we have a bijection from $A_{\mathbf{0},\,B}$ to $ A_{\mathbf{d},\,B}$ given by $\mathbf{x}\mapsto \mathbf{x}+\mathbf{y}$. 
   Consequently, $A_{\mathbf{d},\,B}$ is a finite algebraic set defined over $\Kalg$. 
   Since $\mathbb{C}_v$ is algebraically closed, the coordinate ring of this set over $\Kalg$ is Artinian, and all its residue fields are algebraic over $\Kalg$ by Hilbert's Nullstellensatz; thus they are equal to $\Kalg$ as $\Kalg$ is algebraically closed. 
   Since every point of $A_{\mathbf{d},\,B}$ factors through one of these residue fields, we have $A_{\mathbf{d},\,B}\subset\Mat_{r\times1}(\Kalg)$.
\end{proof}

The following proposition is a $v$-adic counterpart of Proposition 3.1.3 in \cite{Anderson2004}, which makes it easy to check the condition for applying the ABP criterion. 
\begin{Proposition}\label{Prop:ABP_condition}
    Given $\Phi \in \Mat_{r\times r} (\Kalg[1/z])$ and $  \psi \in \Mat_{r\times 1} (\mathbb{T}_v)$
    such that $\psi^{(-1)}=\Phi \psi$, we have $\psi \in \Mat_{r\times 1} (\Entirev)$.
\end{Proposition}
Regarding \cref{Prop:ABP_condition}, if we further assume $\det \Phi|_{z^{-1}=0}\neq 0$, then we can prove the proposition by the same argument as that in \cite{Anderson2004}.
In order to deal with the series treated in \cref{Example:divergent_points}, we have to remove the assumption $\det \Phi|_{z^{-1}=0}\neq 0$.
Though we need a somewhat complicated argument to avoid the assumption, the basic idea of the proof is still the same as the one in the previous work.
We use the following notation for the estimates in the proof; for any matrix $X=(X_{ij}) \in \Mat_{m\times n} (\mathbb{C}_v)$ of any size, we write $||X||_v$ for the quantity 
\begin{equation}
\max_{1\leq i\leq m,\,1 \leq j \leq n}(|X_{ij}|_v),
\end{equation}
the Gauss norm induced by the $v$-adic absolute value on $\mathbb{C}_v$.

\begin{proof}[Proof of \cref{Prop:ABP_condition}]
Let us write
\begin{equation}
    \Phi=\sum_{i=0}^N b_{(i)} z^{-i} \text{ and } \psi =\sum_{i=0}^\infty a_{(i)} z^{-i}
\end{equation}
with $ b_{(i)} \in\Mat_{r\times r} (\Kalg) $ and $a_{(i)} \in \Mat_{r\times 1} (\mathbb{C}_v)$; we allow $b_{(i)}$ to be $0$ and assume that $N\geq 1$.
Then, the equation $\psi^{(-1)}=\Phi \psi$ implies
\begin{equation}
    a_{(n)}^{(-1)}=\sum_{i=0}^{\min(n,\,N)}b_{(i)}a_{(n-i)}
\end{equation}
for any $n\geq 0$. 
Thus, writing $B:=b_{(0)}$, we have
\begin{equation}\label{eq:a_n_recursion}
    a_{(n)}-B^{(1)}a_{(n)}^{(1)}=\sum_{i=1}^{\min(n,\,N)}b_{(i)}^{(1)}a_{(n-i)}^{(1)}=:c_{(n)}.
\end{equation}
    Then \cref{Lemma:ABPCondition} and induction on $n$ imply $a_{(n)} \in \Mat_{r\times 1}(\Kalg)$.

    By \cref{eq:a_n_recursion}, we have
    \begin{equation}
        a_{(n)}=c_{(n)}+\sum_{j=2}^m B^{(1)}\cdots B^{(j-1)} c_{(n)}^{(j-1)}+B^{(1)}\cdots B^{(m)}a_{(n)}^{(m)}
    \end{equation}
    for any $m \geq 1$.
    As we have assumed that the entries of $\psi$ are in $\mathbb{T}_v$, we have $\lim_{n\to \infty}||a_{(n)}||_v=0$; thus, for any $\lambda \geq 1$, we can take $n_0$ such that $||a_{(n)}||_v<\lambda^{-\frac{q^2}{q-1}}$ holds for each $n \geq n_0$. 
    If we choose $\lambda>\max (1,\,||B||_v )$, then we have
    \begin{equation}
        ||B^{(1)}\cdots B^{(m)}a_{(n)}^{(m)}||_v\leq
        ||B^{(1)}\cdots B^{(m)}||_v\cdot ||a_{(n)}^{(m)}||_v\leq \lambda^{q\frac{q^m-1}{q-1}}\cdot \lambda^{-q^m\frac{q^2}{q-1}}
    \end{equation}
    for any $n \geq n_0$.
    Since the quantity $\lambda^{\frac{q^{m+1}-q}{q-1}}\cdot \lambda^{-q^m\frac{q^2}{q-1}}=\lambda^{\frac{-q^{m+2}+q^{m+1}-q}{q-1}}$ goes to $0$ as $m$ goes to infinity, we have
    \begin{equation}\label{eq:a_n_closed}
        a_{(n)}=c_{(n)}+\sum_{j=2}^{\infty} B^{(1)}\cdots B^{(j-1)} c_{(n)}^{(j-1)}
    \end{equation}
    for $n \geq n_0$; we mention that the sum on the right-hand side corresponds to $\tilde{a}_{(n)}$ in the proof of Proposition 3.1.3 of \cite{Anderson2004}.
    Let us consider the field $\mathbf{F}$ given by adjoining all entries of $B$, $b_{(1)},\,\dots,\,b_{(N)}$, and $a_{(0)},\,\dots,\, a_{(n_0)}$ to $\Kcomp$; we then claim, by induction on $n$, that entries of $a_{(n)}$ are in $\mathbf{F}$ for any $n$.
    Indeed, if the claim holds for $a_{(0)},\,\dots,\,a_{(n-1)}$, then the definition of $c_{(0)},\,\cdots,\,c_{(n)}$ in \cref{eq:a_n_recursion} and \cref{eq:a_n_closed} show that the entries of $a_{(n)}$ are limits of sequences of elements of $\mathbf{F}$; since $\mathbf{F}$, a finite extension of $\Kcomp$, is topologically closed, the claim follows for $a_{(n)}$.

    By the definition of $c_{(n)}$ in \cref{eq:a_n_recursion}, we have $\lim_{n \to \infty}||c_{(n)}||_v=0$; thus, by replacing $n_0$ by a larger one if necessary, we further have $||c_{(n)}||_v<\lambda^{-\frac{q}{q-1}}$ for $n \geq n_0$.
    Then we have 
    \begin{align}
        &\left|\middle| B^{(1)}\cdots B^{(j-1)} c_{(n)}^{(j-1)}\middle|\right|_v
        \leq
        \lambda^{q\frac{q^{j-1}-1}{q-1}}||c_{(n)}||_v^{q^{j-1}-1}\cdot||c_{(n)}||_v\\
        &\leq
       \lambda^{\frac{q^{j}-q}{q-1}}\lambda^{-(q^{j-1}-1)\frac{q}{q-1}}||c_{(n)}||_v
        =%\lambda^{\frac{-q^{2}+q}{q-1}}||c_{(n)}||_v
        %\leq 
        ||c_{(n)}||_v
    \end{align}
    for any $j \geq 2$.
   Let us fix $C \geq 1$; then, for large $n$, we have
    \begin{align}
        &C^n ||a_{(n)}||_v
        \leq  
        C^n ||c_{(n)}||_v
        \leq 
        \max_{i=1}^N \left(C^n ||b_{(i)}^{(1)}||_v\cdot ||a_{(n-i)}^{(1)}||_v\right)\\
        &\leq
        \max_{i=1}^N \left(C^i ||b_{(i)}^{(1)}||_v\right)\max_{i=1}^N \left(  ||a_{(n-i)}||_v\right)^{q-1}
        \max_{i=1}^N \left(C^{n-i} ||a_{(n-i)}||_v\right)\\
        &\leq \max_{i=1}^N \left(C^{n-i}||a_{(n-i)}||_v\right);
    \end{align}
    thus we have
    \begin{equation}
        \sup_{n=0}^\infty C^n ||a_{(n)}||_v<\infty \text{ and } \limsup_{n\to \infty} ||a_{(n)}||_v^{1/n}\leq C^{-1}.
    \end{equation}
    Hence, we have
    \begin{equation}
        \lim_{n\to \infty}||a_{(n)}||_v^{1/n}=0
    \end{equation}
    as $C$ is arbitrary, and we can conclude that entries of $\psi$ are entire.
\end{proof}

\cref{Prop:ABP_condition} gives us another justification of \cref{Prop:v-adic polylog entire}.
In addition, we have the following example.
\begin{Example} 
    Let us consider the situation in \cref{Example:divergent_points}. 
    For each $s^\prime$ larger than $s$, the second entry of the column vector $\Omegav^{s^\prime} z^{-s}\psi_\iota $ is
    \begin{align}
         
         &
        \left(\lim_{m\to \infty}\left(\delta_{\mathcal{C},\,m}^{(1)}\right)^{s^\prime}  z^{-s} \right)\cdot\left( \lim_{m\to \infty}
        \left(
        \left(z^s-1\right) 
        \prelog_{\tilde{\mathbf{a}},\,m}^{(1)}+\frac{\overline{a}}{z^{(m+1)s}}\left( \delta_{\mathcal{C},\,m}^{(1)}\right)^{-s}
        \right) \right)\\
        &=\lim_{m\to \infty}
        \left(
        
        \left(\delta_{\mathcal{C},\,m}^{(1)}\right)^{s^\prime}(1-z^{-s})\prelog_{\tilde{\mathbf{a}},\,m}^{(1)}+\left(\delta_{\mathcal{C},\,m}^{(1)}\right)^{s^\prime-s}\frac{\overline{a}}{z^{(m+2)s}}
        \right);\label{eq:non-conv_insidelimit}
    \end{align}
    \cref{Thm:comparison} shows that this limit belongs to $\Acomp[[z,$ $1/z\}[\omegav^{(1)}]=\Acomp[[z,\,1/z\}[\Omegav^{-1}]$; taking $s^\prime$ sufficiently large, we may assume that the limit is in $\Acomp[[z,$ $1/z\}$.
    For each $m$, we note that the expression inside the limit in \cref{eq:non-conv_insidelimit} belongs to $\Acomp[1/z]$; thus the limit is in $\Acomp[[1/z]]$ 
    %since the subring $\Acomp[[1/z]] \cap \Acomp[[z,\,1/z\}$ is topologically closed in $\Acomp[[z,\,1/z\}$, see also \cref{Rem:Topology_and_convergence};
    since \cref{Rem:Topology_and_convergence} yields coefficientwise convergence;
    thus, we have 
    $
    \Omegav^{s^\prime} z^{-s}\psi_\iota \in \Mat_{2 \times 1}( \Acomp[[1/z]]\cap \Acomp[[z,\,1/z\}) \subset \Mat_{2 \times 1}( \mathbb{T}_v). 
    $
    The difference equation in \cref{Example:divergent_points}, together with $\Omegav^{(-1)}=(1-\zeta/z)\Omegav$, gives
\begin{equation}
    \left(\Omegav^{s^\prime}z^{-s}\psi_\iota\right)^{(-1)}
    =
    \begin{pmatrix}
        (1-\zeta/z)^{s^\prime} & 0 \\
        a z^{-s}(1-\zeta/z)^{s^\prime-s}
        & z^{-s}(1-\zeta/z)^{s^\prime-s}
    \end{pmatrix}
    \left(\Omegav^{s^\prime}z^{-s}\psi_\iota\right).
\end{equation}
Since $s^\prime>s$, the coefficient matrix belongs to $\Mat_{2\times2}(\Kalg[1/z])$.
Therefore, \cref{Prop:ABP_condition} shows that
$\Omegav^{s^\prime}z^{-s}\psi_\iota
\in \Mat_{2\times1}(\Entirev)$, that is,
$z^{-s}\psi_\iota
\in \Mat_{2\times1}(\Entirev[1/\Omegav])$.
    Therefore, we can apply \cref{Thm:period_vector_algebraic_independence} to obtain
    \begin{align}
        & \trdeg_{\Kalg\left(\frac{1}{z}\right)}\Kalg\left(\frac{1}{z},\, \lim_{m\to \infty}
        \left(
        \left(z^s-1\right) 
        \prelog_{\tilde{\mathbf{a}},\,m}^{(1)}+\frac{\overline{a}}{z^{(m+1)s}}\left( \delta_{\mathcal{C},\,m}^{(1)}\right)^{-s}
        \right)  \right)\\
        &=\trdeg_{\Kalg}\Kalg\left(  \lim_{m\to \infty}
        \left(
        \left(z^s-1\right) 
        \prelog_{\tilde{\mathbf{a}},\,m}^{(1)}+\frac{\overline{a}}{z^{(m+1)s}}\left( \delta_{\mathcal{C},\,m}^{(1)}\right)^{-s}
        \middle)\right|_{z=\zeta}\right).
    \end{align}
\end{Example}

\section{A difference-Galois bridge between the $\infty$-adic and $v$-adic sides}\label{sec:difference galois}

The aim of this section is to relate the transcendence properties of deformations of $\infty$-adic periods to those of $v$-adic pre-periods.
To this end, we adapt techniques from difference Galois theory.
The discussion in this section culminates in \cref{Prop:PVTmain} and \cref{Cor:PVTmain}, which play a key role in the transcendence theory developed in the present paper.
The purpose of this section is not to develop new results in difference Galois theory itself, but rather to use its methods to compare the transcendence properties of deformation series arising from the $\infty$-adic and $v$-adic settings.
The arguments in this section are inspired in part by recent studies of Picard-Vessiot theory, especially those in \cite{Maurischat2017}.

Let us recall some notation of difference Galois theory; we refer the reader to \cite[Chapter 1]{Put1997}.
A \textit{difference ring} is a pair $(R,\,\phi)$ of a ring $R$ and its automorphism $\phi$. 
If the automorphism $\phi$ is clear from the context, we simply refer to the difference ring as $R$.
The set of elements in $R$ fixed by $\phi$ is denoted by $R^\phi$ and is a subring of $R$.
If $R$ is an integral domain or a field, we call the pair a \textit{difference domain} or a \textit{difference field}, respectively.
For a given difference ring $(R,\,\phi)$ and an ideal $\mathfrak{m}$ of $R$, we call $\mathfrak{m}$ a \textit{difference ideal} if $\mathfrak{m}$ is stable under $\phi$, that is, $\phi(\mathfrak{m})\subset \mathfrak{m}$.
Difference rings which have no difference ideals other than the trivial ones $(0)$ and $(1)$ are said to be \textit{simple}.
For a matrix $\Psi \in \Mat_{r\times r^\prime}(R)$, we write ${\phi^m}\Psi$ or ${\phi^m}(\Psi)$ for the matrix given by applying the $m$-th iterate $\phi^m$ entrywise.

\begin{Example}\label{Example:Tate algebras as difference rings}
    For $n\geq 1$, we can consider the difference domain $(\Frac\mathbb{T}_\infty, \,\phi=\tau^{-n})$,
    where $\tau$ is the automorphism of $\mathbb{T}_\infty$ given by the first twist $f\mapsto f^{(1)}$, extended to $\Frac\mathbb{T}_\infty$ by $\tau(f/g)=f^{(1)}/g^{(1)}$ for $f,g\in\mathbb{T}_\infty$ with $g\neq0$; this extension is unique.
    The fixed ring $(\Frac\mathbb{T}_\infty)^\phi$ is $\mathbb{F}_{q^n}(t)$ as shown by the same arguments as those in \cite[Lemma 3.3.2]{Papanikolas2008}.
    On the other hand, for any finite place $v$, the twist $\tau:f\mapsto f^{(1)}$ defined in \cref{eq:def_v_twisting} is an automorphism of $\mathbb{T}_v$ and extends uniquely to an automorphism of $\Frac\mathbb{T}_v$ in the same way.
    Thus, with $\phi=\tau^{-n}$, the pair $(\mathbb{T}_v,\phi)$ is a difference domain and $(\Frac\mathbb{T}_v,\phi)$ is a difference field for any $n \geq 1$.
\end{Example}

\begin{Lemma}\label{Lemma:fixed_field_of_TT_v}
    For each $n\geq1$, the ring of fixed elements in $\mathbb{T}_v$ under the automorphism $\phi:=\tau^{-n}$ is $\mathbb{F}_{q^n}[1/z]$ and the difference domain $(\Frac\mathbb{T}_v, \,\phi=\tau^{-n} )$ has the fixed field $\mathbb{F}_{q^{n}}(z)$.
\end{Lemma}
\begin{proof}
    Although this lemma can be justified by an argument similar to that in the proof of \cite[Lemma 3.3.2]{Papanikolas2008}, we include a proof to keep the paper self-contained.
    For $a \in \mathbb{C}_v$, we have $a=a^{q^{-n}}$ if and only if $a \in \mathbb{F}_{q^n}$; thus an element $f \in\mathbb{T}_v $ fixed by $\phi=\tau^{-n}$ is in $\mathbb{F}_{q^n}[[1/z]]$.
    
    As elements of $\mathbb{F}_{q^n}^{\times}$ have $v$-adic absolute value $1$ and the coefficients of $f$ tend to $0$, the fixed element $f$ is in $\mathbb{F}_{q^n}[1/z]$.
    As $\mathbb{C}_v$ is algebraically closed, for $f \in (\Frac\mathbb{T}_v)^\times$, we have a unique factorization
    \begin{equation}
        f=
        \left(
        z^{-1}-a_1
        \right)^{m_1}
        \cdots
        \left(
        z^{-1}-a_r
        \right)^{m_r}
        f^\prime
    \end{equation}
    with pairwise distinct $a_1,\,\dots,\,a_r$ such that $|a_i|_v\leq1$ for $i=1,\,\dots,\,r$, $m_1,\,\dots,\,m_r \in \mathbb{Z}\setminus \{0\}$, and $f^\prime \in \mathbb{T}_v^\times$ by the Weierstrass preparation theorem (see \cite[\S 2.2, Corollary 9]{Bosch2014} for example). If $\phi f =f $, we have
    \begin{equation}
        \left(
        z^{-1}-a_1
        \right)^{m_1}
        \cdots
        \left(
        z^{-1}-a_r
        \right)^{m_r}
        =
        \left(
        z^{-1}-a_1^{q^{-n}}
        \right)^{m_1}
        \cdots
        \left(
        z^{-1}-a_r^{q^{-n}}
        \right)^{m_r}
    \end{equation}
    by the uniqueness of the factorization; this means that $a \mapsto a^{q^{-n}}$ is a permutation of the set $\{a_1,\,\dots,\,a_r\}$;
    thus, if $n_0$ is the order of this permutation, then we have $a_i=a_i^{q^{-nn_0}}$ and hence $a_i \in \mathbb{F}_{q^{nn_0}}$ for all $1 \leq i \leq r$.
    We can take a polynomial $P(z^{-1}) \in \mathbb{F}_{q}[z^{-1}]$ in the variable $z^{-1}$ whose order of vanishing at $z^{-1}=a_i$ is equal to or greater than $-m_i$ for each $i$; then the product $P(z^{-1})f$ has no pole in the closed unit disk; thus it is in $\mathbb{T}_v$. As $P(z^{-1})f$ is also fixed by $\phi$, it is an element of the fixed part $\mathbb{F}_{q^n}[1/z]$ of $\mathbb{T}_v$ and thus $(\Frac \mathbb{T}_v)^\phi \subset \mathbb{F}_{q^n}(z)$; the opposite inclusion is clear.
\end{proof}

\begin{comment}
\begin{Proposition}
    If $(R,\,\phi)$ is simple difference domain, then $R^{\phi}=\Frac(R^\phi)=(\Frac R)^\phi$.
\end{Proposition}
\begin{proof}
    
\end{proof}
\end{comment}
Let us take a difference field $(\mathbf{K},\,\phi)$ and a matrix $\Phi \in \GL(\mathbf{K})$; then we simply refer to the first-order linear difference system $\phi Y=\Phi Y$ in the sense of \cite[\S 1]{Put1997} as the difference equation given by the pair $\big(\mathbf{K},\,\Phi\big)=\big((\mathbf{K},\,\phi),\,\Phi\big)$ for short.
If $\Psi   \in \GL(\mathbf{K}^\prime)$ with $\mathbf{K}^\prime$ being a difference ring containing $\mathbf{K}$ satisfies $\phi(\Psi)=\Phi \Psi$, then $\Psi$ is said to be a \textit{fundamental matrix} of the difference equation.

    \begin{Lemma}[{\cite[Corollary 4.2.5]{Papanikolas2008}}]\label{lem:PVsimple}
    Take a positive integer $n$; if a difference equation $\big((\Kalg(t),\,\tau^{-n}),\,\Phi\big)$ with $\Phi \in \GL_r(\Kalg(t))$ has a fundamental matrix $\Psi \in \GL_r(\Frac \mathbb{T}_\infty)$, then the ring $\Kalg(t)[\Psi_{ij},\,1/\det \Psi\mid 1\leq i,\,j \leq r]$ is a simple difference domain. 
    \end{Lemma}
    To justify this lemma, we consider the $\tau^{-n}$-admissible triple $(\mathbb{F}_{q^{n}}(t),$ $ \Kalg(t),$ $\Frac\mathbb{T}_{\infty})$ and the fundamental matrix $\Psi$ of the difference equation defined by $\Big( (\Kalg(t),\,\tau^{-n}),\,\Phi\Big)$ in the sense of \cite[\S4]{Papanikolas2008}. Then, \cite[Corollary 4.2.5]{Papanikolas2008} shows that the ring $\Kalg(t)[\Psi_{ij},\,1/\det \Psi\mid 1\leq i,\,j \leq r]$ is isomorphic to a quotient of a difference ring by a maximal proper difference ideal. Thus it is simple and is a domain as it is a subring of $\Frac\mathbb{T}_\infty$.

\begin{Example}\label{Example:inftypolylogperiodinterpretation}
    Let us consider the matrices $\Phi_{\tilde{\mathbf{a}}}$ defined in \cref{eq:def matrix T} and $\Psi_{\tilde{\mathbf{a}}}$ defined in \cref{eq:def matrix Psi} under the assumption $\tilde{\mathbf{a}}=\big((s_1,\,a_1),\,\dots,\,(s_r,\,a_r)\big)\in \tilde{\mathbb{D}}_{\infty}^r$.
    Let us consider the automorphism $\phi:=\tau^{-1}$ on $\Frac \mathbb{T}_\infty$. 
    Then \cref{eq:inftypolylogmatrixdifferenceeq} can be rewritten as $\phi\Psi_{\tilde{\mathbf{a}}}=\Phi_{\tilde{\mathbf{a}}}\Psi_{\tilde{\mathbf{a}}}$.
    Therefore, if we consider the ring $\mathbf{R}=\Kalg(t)[ \Psi_{\tilde{\mathbf{a}}\,ij},\,1/\det \Psi_{\tilde{\mathbf{a}}}\mid 1\leq i,\,j\leq r+1 ]$, then  the difference ring $\big(\mathbf{R},\,\tau^{-1}\big)$ is a simple difference domain and the matrix $\Psi_{\tilde{\mathbf{a}}}$ is a fundamental matrix of the difference equation given by the pair $(\Kalg(t),\, \Phi_{\tilde{\mathbf{a}}})$.
    
    By induction, we can show that $\Psi_{\tilde{\mathbf{a}}}$ is a fundamental matrix of the difference equation given by $\big((\Kalg(t),\,\tau^{-n}),\, \tau^{-n+1}(\Phi_{\tilde{\mathbf a}})\cdots\tau^{-1}(\Phi_{\tilde{\mathbf a}})\Phi_{\tilde{\mathbf a}}\big)$ for each $n \geq 1$, and we can also prove that
    $(\mathbf{R},\,\tau^{-n}\big)$ is also a simple difference domain.
\end{Example}

We will consider difference fields $\mathbf{P}/\mathbf{K}$ such that $\mathbf{P}^{\phi^m}=\mathbf{K}^{\phi^m}$ and $\mathbf{K}^{\phi^m}/\mathbf{K}^{\phi}$ is a finite extension for any $m\geq 1$.
Then $\mathbf{P}$ and $\mathbf{K}$ are simple as they are fields.
The typical case we consider is where $\mathbf{K}$ equals $\Kalg(t)$ and a difference field $\mathbf{P}$ is given by adjoining to $\mathbf{K}$ a finite number of elements of $\Frac\mathbb{T}_v$ for a degree-one $v$ such that the automorphism $\phi=\tau^{-1}$ of $\Frac\mathbb{T}_v$ (\cref{Example:Tate algebras as difference rings}) can be restricted to an automorphism $\phi=\phi|_{\mathbf{P}}$ on $\mathbf{P}$; 
since $\deg v=1$, we have $z=t-\alpha$ for some $\alpha\in\mathbb{F}_q$, and hence $\Kalg(t)=\Kalg(z)$.
Then we have $\mathbf{K}^{\phi^m}=\mathbf{P}^{\phi^m}=\mathbb{F}_{q^{m}}(z)$ for each $m\geq1$ since we have $\mathbb{F}_{q^{m}}(z)\subset \mathbf{K}^{\phi^m}\subset \mathbf{P}^{\phi^m} \subset (\Frac \mathbb{T}_v)^{\tau^{-m}}= \mathbb{F}_{q^{m}}(z) $.

\begin{Example}\label{Example:mathbfPv-adicpolylog}
    We continue to use the setting and notation in \cref{Example:v-adicpolyloginterpretation}.
   If we assume $\deg v=1$, put $\mathbf{K}:=\Kalg(z)$, and define  
   \begin{equation}
       \mathbf{P}=\mathbf{K}(\prelog_{v,\,ji}^{(1)}\mid 1\leq i <j\leq r+1) \subset \Frac \mathbb{T}_v,
   \end{equation}
   then $\mathbf{P}$ is stable under the action of $\phi=\tau^{-1}$ on $\Frac\mathbb{T}_v$. 
   Indeed, taking $m=0$ in \cref{eq:prelogvdifferenceeq}, we obtain
\begin{equation}\label{eq:tau-stability}
(z-\zeta)^{s_i+\cdots+s_{j-1}}\,\prelog_{v,\,ji}=\prelog_{v,\,ji}^{(1)}+a_{j-1}\,\prelog_{v,\,(j-1)i}^{(1)},
\end{equation}
which implies $\tau^{-1}(\mathbf{P})\subseteq \mathbf{P}$.
Conversely, induction on $j-i$ shows the opposite inclusion as $\prelog_{v,\,ii}=1$; hence $\tau^{-1}(\mathbf{P})=\mathbf{P}$.
\end{Example}

Our next goal is to establish \cref{Prop:PVTmain} and its corollaries.
This proposition plays a key role when we relate the algebraic independence of $v$-adic pre-periods to that of the corresponding $\infty$-adic periods,  as in \cref{Cor:inftyprelogvprelog}.
In order to prove \cref{Prop:PVTmain}, we will use the following two lemmas:

\begin{Lemma}\label{lemma:alg. ext. of const. field}
    Let $(\mathbf{P}, \phi)$ be a difference field and $(\mathbf{S},\,\phi)$ be a simple difference ring containing $\mathbf{P}$.
 If $\mathbf{S}$ is finitely generated over $\mathbf{P}$, then the field extension $\mathbf{S}^{\phi}/\mathbf{P}^\phi$ is algebraic.
\end{Lemma}
This lemma is a consequence of a well-known argument (see, for example, \cite[Lemma 1.8]{Put1997}), and we refer the reader to \cite[Theorem 1.2.6]{Maier2011}, where this fact is stated explicitly.
Although the lemma below is also proven by standard arguments in Picard-Vessiot theory (see \cite[Corollary 1.16]{Put1997} for instance), 
we reproduce the detailed argument to make the present paper self-contained.

\begin{Lemma}\label{Lemma:difference domain}
    Take a simple difference ring  $(\mathbf{R},\,\phi)$ with $\mathbf{R}$ Noetherian. 
    If the number of minimal prime ideals of $\mathbf{R}$ is $m \geq 1$ and $\mathfrak{n}$ is one of them, then the iterate $\phi^{m}$ stabilizes $\mathfrak{n}$ and the difference ring $(\mathbf{R}/\mathfrak{n},\,\phi^m)$ is a simple difference domain.
\end{Lemma}

\begin{proof}
    Since the ring $\mathbf{R}$ has a finite number $m$ of minimal prime ideals, we have $1\leq e\leq m$ such that $\phi^e \mathfrak{n}=\mathfrak{n}$ as the automorphism $\phi$ permutes the  minimal prime ideals of $\mathbf{R}$; we choose the smallest such $e$.
    The intersection $\bigcap_{j=0}^{e-1}\phi^j\mathfrak{n}$ is stable under the action of $\phi$ and hence is equal to $(0)$ as $\mathbf{R}$ is assumed to be simple; thus any minimal prime ideal $\mathfrak{n}^\prime$ contains the intersection $\bigcap_{j=0}^{e-1}\phi^j\mathfrak{n}$ and hence we have $0 \leq j \leq e-1$ such that $\phi^j\mathfrak{n} \subset \mathfrak{n}^\prime$ as a prime ideal containing a finite intersection of ideals contains one of them; then we have $\phi^j\mathfrak{n} = \mathfrak{n}^\prime$ as $\mathfrak{n}^\prime$ is a minimal prime.
    Hence we have $m=e$, that is, $ \mathfrak{n},\,\phi \mathfrak{n},\,\dots,\,\phi^{m-1}\mathfrak{n}$ are distinct, and $\mathfrak{n}$ is stable under the action of $\phi^m$.
    %Therefore, we have $e=m$ by the uniqueness of the primary decomposition of $(0)$;if we let $\mathfrak{n}=\mathfrak{n}_0,\,\mathfrak{n}_1,\,\dots,\,\mathfrak{n}_{m-1}$ be distinct minimal prime ideals of $\mathbf{R}$, then the intersection $\mathfrak{n}_0\cap\cdots\cap\mathfrak{n}_{m-1}$ is equal to the trivial ideal as $\mathbf{R}$ is simple (nilpotent radicals of difference rings must be a difference ideal).That is, $\mathfrak{n}$ is stable under the action of $\phi^m$.
    
    In order to prove that the difference domain $(\mathbf{R}/\mathfrak{n},\,\phi^m)$ is simple, we show that $\mathbf{R}$ has no proper $\phi^m$-stable ideals containing $\mathfrak{n}$ other than $\mathfrak{n}$ itself.
    Let us take a proper ideal $\mathfrak{n}^{\prime\prime}$ which is $\phi^m$-stable and contains $\mathfrak{n}$. 
    Then the intersection $\bigcap_{j=0}^{m-1}\phi^j \mathfrak{n}^{\prime\prime}$ is a proper $\phi$-stable ideal as $\phi^m \mathfrak{n}^{\prime\prime}\subset\mathfrak{n}^{\prime\prime}$ and it is contained in $\mathfrak{n}^{\prime\prime}$, and hence this intersection is equal to $(0)$, which is a subset of $\mathfrak{n}$. Since $\mathfrak{n}$ is prime, we have $0\leq j \leq m-1$ such that $\phi^j\mathfrak{n}^{\prime\prime}\subset \mathfrak{n}$.
    This inclusion implies $\phi^j\mathfrak{n} \subset \phi^j\mathfrak{n}^{\prime\prime}\subset \mathfrak{n}$, which gives $\phi^j\mathfrak{n} =\mathfrak{n}$ as $\mathfrak{n}$ is minimal.
    As $ \mathfrak{n},\,\phi \mathfrak{n},\,\dots,\,\phi^{m-1}\mathfrak{n}$ are distinct minimal prime ideals, we have $j=0$, so we conclude that  $\mathfrak{n}^{\prime\prime} =\phi^j \mathfrak{n}^{\prime\prime}\subset \mathfrak{n}$.
    Thus, we have $\mathfrak{n}^{\prime\prime}=\mathfrak{n}$, which shows that the difference domain $(\mathbf{R}/\mathfrak{n},\,\phi^m)$ is simple.
\end{proof}

In the applications we have in mind, we consider two difference rings; 
one, denoted by $\mathbf{R}$, is finitely generated over $\mathbf{K}=\Kalg(t)$ by elements of $\Frac\mathbb{T}_\infty$, whereas the other, $\mathbf{P}\subset\Frac\mathbb{T}_v$, is obtained by adjoining $v$-adic pre-periods to $\mathbf{K}$.
The following proposition gives us a simple difference domain $\mathbf{S}$; in the situations arising in our applications, its first corollary (\cref{Cor:PVTmain_injective}) then shows that both $\mathbf{R}$ and $\mathbf{P}$ are difference subrings of $\mathbf{S}$ after replacing the automorphism $\phi$ of $\mathbf{R}$ and $\mathbf{P}$ by a suitable iterate $\phi^m$.

\begin{Proposition}\label{Prop:PVTmain}
Let $(\mathbf{P},\phi)$ be a difference field and $(\mathbf{R},\phi)$ a difference ring, both containing the same difference subfield $(\mathbf{K},\phi)$.
If $\mathbf{R}$ is a finitely generated $\mathbf{K}$-algebra,
then there exist an integer $m\ge 1$ and a simple difference domain $(\mathbf{S},\varphi)$ such that
\begin{equation}
(\mathbf{R}\otimes_{\mathbf K}\mathbf P,\phi^m)\twoheadrightarrow (\mathbf S,\varphi)
\end{equation}
is a surjective difference homomorphism.
Moreover, the extension of fields of constants
$
\mathbf S^\varphi / \mathbf P^{\phi^m}
$
is algebraic.
\end{Proposition}

\begin{proof}

Let us consider the tensor product $\mathbf{R}_1:=\mathbf{R} \otimes_{\mathbf{K}} \mathbf{P}$ on which $\phi$ acts diagonally.
Since $\mathbf{P}$ is a field extension of $\mathbf{K}$ and $\mathbf{R}$ is a non-zero $\mathbf{K}$-algebra, the tensor product is non-zero.
Hence, by the same argument as that for the existence of a maximal ideal in a commutative ring using Zorn's lemma, we can obtain a maximal difference ideal $\mathfrak{m}$ of $\mathbf{R}_1$; as $\mathbf{R}_1$ is finitely generated over the field $\mathbf{P}$, it is Noetherian and thus we have $\phi(\mathfrak{m})=\mathfrak{m}$.
We let $\mathbf{R}_2:=\big(\mathbf{R} \otimes_{\mathbf{K}} \mathbf{P}\big)/\mathfrak{m}$; this quotient ring has a difference ring structure induced by the automorphism $\phi$ of $\mathbf{R}_1$ as $\mathfrak{m}$ is a difference ideal.
Since the difference ring $\mathbf{R}_2$ is finitely generated over the field $\mathbf{P}$ by construction, the ring $\mathbf{R}_2$ is Noetherian as well; by the maximality of $\mathfrak{m}$, the difference ring $(\mathbf{R}_2,\,\phi)$ is simple.
If $m$ is the number of minimal prime ideals of $\mathbf{R}_2$, then every minimal prime ideal $\mathfrak{n}$ is $\phi^{m}$-invariant and we put $\mathbf{S}:=(\mathbf{R}_2/\mathfrak{n},\,\varphi)$ where $\varphi$ is the automorphism of $\mathbf{R}_2/\mathfrak{n}$ induced by $\phi^m$ on $\mathbf{R}_2$. 
Clearly, the difference ring $(\mathbf{S},\,\varphi)$ admits a surjective difference homomorphism $(\mathbf{R}\otimes_{\mathbf{K}}\mathbf{P},\,\phi^m)\twoheadrightarrow(\mathbf{S},\,\varphi)$ and \cref{Lemma:difference domain} implies that it is a simple difference domain.
Since $\mathbf{P}$ is a field, the composition $\mathbf{P}\rightarrow\mathbf{R}_1 \rightarrow\mathbf{R}_2\rightarrow \mathbf{S} $ is an injection; as the domain $\mathbf{S}$ is finitely generated over $\mathbf{P}$,
the field extension $\mathbf{S}^{\varphi}/\mathbf{P}^{\phi^m}$ is algebraic by \cref{lemma:alg. ext. of const. field}.
\end{proof}

Since the kernel of a difference homomorphism is a difference ideal, every (unital) difference homomorphism from a simple difference ring is injective.
In particular, this applies to difference fields, as we have already seen for $\mathbf{P}$ in the proof of \cref{Prop:PVTmain}. 
Thus we obtain the following corollary.

\begin{Corollary}\label{Cor:PVTmain_injective}
    Consider the situation of \cref{Prop:PVTmain}; then the difference homomorphism $(\mathbf{P},\,\phi^m)\rightarrow (\mathbf{S},\,\varphi)$ given by the difference homomorphism $(\mathbf{R}\otimes_{\mathbf K}\mathbf P,\phi^m)\twoheadrightarrow (\mathbf S,\varphi)$ is injective. If $(\mathbf{R},\,\phi^m)$ is simple, the homomorphism $(\mathbf{R},\,\phi^m)\rightarrow (\mathbf{S},\,\varphi)$ is also an injection.
\end{Corollary}

The following corollary is the form of \cref{Prop:PVTmain} used in our application to transcendence questions. 
It will be applied to a difference equation admitting two fundamental matrices: one fundamental matrix has entries built from the deformation series of $\infty$-adic periods, while the other is obtained by replacing some of those series with the corresponding $v$-adic pre-periods.
We use this result to carry out the comparison between the $\infty$-adic and $v$-adic sides in \cref{Cor:inftyprelogvprelog}.

\begin{Corollary}\label{Cor:PVTmain}
Let $(\mathbf{P},\phi)$ be a difference field and $(\mathbf{R},\phi)$ a difference ring, both containing the same difference subfield $(\mathbf{K},\phi)$.
Let $m$ and $(\mathbf{S},\varphi)$ be as in \cref{Prop:PVTmain}.
Assume that $(\mathbf{R},\phi^m)$ is simple, that $\mathbf{R}$ is a finitely generated $\mathbf{K}$-algebra, and that $\mathbf{P}^{\phi^m}/\mathbf{K}^{\phi^m}$ is finite.
Suppose that a matrix $\Phi\in\GL_r(\mathbf{K})$ admits two fundamental matrices $\Theta=(\Theta_{ij})$ and $\Psi=(\Psi_{ij})$ in $\GL_r(\mathbf{S})$; that is, we have $\varphi(\Theta)=\Phi \Theta$ and $\varphi(\Psi)=\Phi \Psi$.
If $\mathbf{P}\subset \mathbf{K}^{\mathrm{alg}}(\Theta_{ij})$, then
\begin{equation}
\trdeg_{\mathbf{K}}\mathbf{P}
=
\trdeg_{\mathbf{K}}\mathbf{K}(\Psi_{ij})
-\trdeg_{\mathbf{P}}\mathbf{P}(\Theta_{ij}).
\end{equation}
\end{Corollary}

\begin{proof}
By \cref{Cor:PVTmain_injective}, we regard both $\mathbf{R}$ and $\mathbf{P}$ as difference subrings of $\mathbf{S}$.
\cref{Prop:PVTmain} implies that the extension $\mathbf{S}^{\varphi}/\mathbf{P}^{\phi^m}$ is algebraic;
then we can see that the extension $\mathbf{S}^{\varphi}/\mathbf{K}^{\phi^m}$ is also algebraic as we further assumed that $\mathbf{P}^{\phi^m}/\mathbf{K}^{\phi^m}$ is finite.

Since both $\Psi$ and $\Theta$ are fundamental matrices of the difference equation given by the matrix $\Phi$, we obtain
\begin{equation}
    \varphi(\Psi^{-1}\Theta)=\Psi^{-1}\Phi^{-1} \Phi \Theta=\Psi^{-1}\Theta,
\end{equation}
    which implies that $\Psi^{-1} \Theta \in \operatorname{GL}(\mathbf{S}^{\varphi})$; thus we can obtain $\mathbf{K}( \Psi_{ij})\cdot\mathbf{S}^{\varphi}=\mathbf{K}( \Theta_{ij})\cdot\mathbf{S}^{\varphi} \subset \Frac \mathbf{S}$, which yields
    $\trdeg_{\mathbf{K}}\mathbf{K}( \Psi_{ij})=\trdeg_{\mathbf{K}}\mathbf{K}( \Theta_{ij})$ since the extension $\mathbf{S}^{\varphi}/\mathbf{K}^{\phi^{m}}$ is algebraic.
As we assumed that $\mathbf{P}$ is a subfield of $\mathbf{K}^{\textrm{alg}}( \Theta_{ij}) \subset (\Frac \mathbf{S})^{\mathrm{alg}}$, we have
\begin{align}
    &\quad \ \ \trdeg_{\mathbf{K}}\mathbf{K}( \Psi_{ij})=\trdeg_{\mathbf{K}}\mathbf{K}( \Theta_{ij})
    =\trdeg_{\mathbf{K}}\mathbf{K}^{\textrm{alg}}( \Theta_{ij})\\
    &=\trdeg_{\mathbf{P}(\Theta_{ij})}\mathbf{K}^{\textrm{alg}}( \Theta_{ij})+\trdeg_{\mathbf{P}}\mathbf{P}(\Theta_{ij})+\trdeg_{\mathbf{K}}\mathbf{P}\\
    &=\trdeg_{\mathbf{P}}\mathbf{P}(\Theta_{ij})+\trdeg_{\mathbf{K}}\mathbf{P},%\\
    %(&\leq \trdeg_{\mathbf{K}}\mathbf{P}+1 \text{ for the case of MZV}).
\end{align}
hence the result.
\end{proof}

As an application, we obtain a description of the transcendence degree of the deformation of $v$-adic Carlitz multiple polylogarithms in terms of the transcendence degree of the deformation of $\infty$-adic Carlitz multiple polylogarithms.

    \begin{Example}\label{Example:Theta_s}
     Consider the situation in \cref{Example:v-adicpolyloginterpretation,Example:mathbfPv-adicpolylog,Example:inftypolylogperiodinterpretation}, and further assume that we have $\tilde{\mathbf{a}}=\big((s_1,\,a_1),\,\dots,\,(s_r,\,a_r)\big)\in \tilde{\mathbb{D}}_{v}^r$.
     As we have assumed $\deg v=1$, we have $z=t-\alpha$ for some $\alpha\in\mathbb{F}_q$, and hence $\Kalg(t)=\Kalg(z)$.
    If we define $\mathbf{P}$ as in \cref{Example:mathbfPv-adicpolylog}, then \cref{Prop:PVTmain} yields $m \geq 1$ and a simple difference domain $(\mathbf{S},\,\varphi)$; since $(\mathbf{R},\,\phi^m)$ is simple, \cref{Cor:PVTmain_injective} allows us to regard both $(\mathbf{R},\,\phi^m)$ and $(\mathbf{P},\,\phi^m)$ as difference subrings of $(\mathbf{S},\,\varphi)$.
    We put
    \begin{equation}
        \Theta_{\tilde{\mathbf{a}}}:=
        \begin{pmatrix}
1 & 0 & \cdots & \cdots & 0\\
\prelog_{v,\,21}^{(1)} &\frac{1}{\Omega_{\infty}^{s_1}} & 0 & \cdots & 0\\
\vdots & \frac{\prelog_{v,\,32}^{(1)}}{\Omega_{\infty}^{s_1}} & \frac{1}{\Omega_{\infty}^{s_1+s_2}} & \ddots & \vdots\\
\prelog_{v,\,r1}^{(1)} & \vdots & \ddots & \ddots & 0\\
\prelog_{v,\,(r+1)1}^{(1)} &\frac{\prelog_{v,\,(r+1)2}^{(1)}}{\Omega_{\infty}^{s_1}} & \cdots & \frac{\prelog_{v,\,(r+1)r}^{(1)}}{\Omega_{\infty}^{s_1+\cdots+s_{r-1}}} & \frac{1}{\Omega_{\infty}^{s_1+\cdots+s_r}}
\end{pmatrix}
    \end{equation}
     where the entries are considered to be elements of $\mathbf{S}$ via the embeddings given by \cref{Cor:PVTmain_injective}. Then it is a fundamental matrix of the difference equation given by $\big((\mathbf{K},\,\phi^{m}),\,\phi^{m-1}(\Phi_{\tilde{\mathbf{a}}})\cdots\phi(\Phi_{\tilde{\mathbf{a}}}) \cdot \Phi_{\tilde{\mathbf{a}}}\big)$, that is, 
     \begin{equation}\label{eq:Thetadiffeq}
         \varphi (\Theta_{\tilde{\mathbf{a}}})= \phi^{m-1}(\Phi_{\tilde{\mathbf{a}}})\cdots\phi(\Phi_{\tilde{\mathbf{a}}}) \cdot \Phi_{\tilde{\mathbf{a}}}\cdot \Theta_{\tilde{\mathbf{a}}}.
         \end{equation}
     Indeed, the matrix $\varphi(\Theta_{\tilde{\mathbf{a}}})$ is given by
     \begin{equation}
          \begin{pmatrix}
1 & 0 & \cdots & \cdots & 0\\
\prelog_{v,\,21}^{(-m+1)} &\frac{1}{\left(\Omega_{\infty}^{(-m)}\right)^{s_1}} & 0 & \cdots & 0\\
\vdots & \frac{\prelog_{v,\,32}^{(-m+1)}}{\left(\Omega_{\infty}^{(-m)}\right)^{s_1}} & \frac{1}{\left(\Omega_{\infty}^{(-m)}\right)^{s_1+s_2}} & \ddots & \vdots\\
\prelog_{v,\,r1}^{(-m+1)} & \vdots & \ddots & \ddots & 0\\
\prelog_{v,\,(r+1)1}^{(-m+1)} &\frac{\prelog_{v,\,(r+1)2}^{(-m+1)}}{\left(\Omega_{\infty}^{(-m)}\right)^{s_1}} & \cdots & \frac{\prelog_{v,\,(r+1)r}^{(-m+1)}}{\left(\Omega_{\infty}^{(-m)}\right)^{s_1+\cdots+s_{r-1}}} & \frac{1}{\left(\Omega_{\infty}^{(-m)}\right)^{s_1+\cdots+s_r}}
\end{pmatrix}
     \end{equation}
     and the matrix 
     \begin{equation}
         \phi^{m-1}(\Phi_{\tilde{\mathbf{a}}})\cdots\phi(\Phi_{\tilde{\mathbf{a}}}) \cdot \Phi_{\tilde{\mathbf{a}}}=\phi^{m-1}
         \left(
            \Phi_{\tilde{\mathbf{a}}}\Phi_{\tilde{\mathbf{a}}}^{(1)}\cdots\Phi_{\tilde{\mathbf{a}}}^{(m-1)}
         \right)
    \end{equation}
    is given by
    \begin{equation}
        \begin{pmatrix}
1 & 0 & \cdots & \cdots & 0\\
\prelog_{21,\,m-1} & \frac{1}{(\lambda_{m-1})^{s_1}} & 0 & \cdots & 0\\
\vdots & \frac{\prelog_{32,\,m-1}}{(\lambda_{m-1})^{s_1}} & \frac{1}{(\lambda_{m-1})^{s_1+s_2}} & \ddots & \vdots\\
\vdots & \vdots & \ddots & \ddots & 0\\
\prelog_{(r+1)1,\,m-1} & \frac{\prelog_{(r+1)2,\,{m-1}}}{(\lambda_{m-1})^{s_1}}& \cdots & \frac{\prelog_{(r+1)r,\,m-1}}{(\lambda_{m-1})^{s_1+\cdots+s_{r-1}}} & \frac{1}{(\lambda_{m-1})^{s_1+\cdots+s_r}}
\end{pmatrix}^{(-m+1)},
    \end{equation}
    see \cref{Example:v-adicpolylogcomparison}; furthermore, we have \cref{eq:prelogvdifferenceeq}, which implies
    \begin{equation}
        \prelog_{v,\,ji}^{(-m+1)}=\sum_{k=i}^{j}\left(\frac{\prelog_{jk,\,m-1}\,}{\lambda_{m-1}^{s_i+\cdots+s_{k-1}}}\right)^{(-m+1)}\prelog_{v,\,ki}^{(1)};
    \end{equation}
    therefore, we have 
    \begin{equation}
        \frac{\prelog_{v,\,ji}^{(-m+1)}}{\left(\Omega_{\infty}^{(-m)}\right)^{s_1+\cdots+s_{i-1}}}=\sum_{k=i}^{j}\left(\frac{\prelog_{jk,\,m-1}\,}{\lambda_{m-1}^{s_1+\cdots+s_{k-1}}}\right)^{(-m+1)}\frac{\prelog_{v,\,ki}^{(1)}}{\Omega_{\infty}^{s_1+\cdots+s_{i-1}}}
    \end{equation}
     since \cref{eq:Omega_inftyFrobeniuseq} implies $\Omega_\infty^{(-m)}=\lambda_{m-1}^{(-m+1)}\Omega_\infty$.
     As the matrix $\Psi_{\tilde{\mathbf{a}}}$ also has the same property (\cref{Example:inftypolylogperiodinterpretation}), 
      \cref{Cor:PVTmain} gives us
     \begin{align}
         &\quad \ \ \trdeg_{\Kalg (t)}\Kalg(t,\,\prelog_{v,\,ji}^{(1)}\mid1\leq i< j \leq r+1)\\
         & \geq \trdeg_{\Kalg (t)}\Kalg(t,\,\Omega_\infty,\, \prelog_{\infty,\,ji}^{(1)}\mid1\leq i< j \leq r+1)-1.
     \end{align}
\end{Example}

Recalling that $\Kalg(t)=\Kalg(z)$ when $\deg v=1$,
and using the inequality
$\trdeg_{\mathbf{P}}\mathbf{P}(\Omega_\infty)\leq 1$,
we obtain the following more general corollary:
\begin{Corollary}\label{Cor:inftyprelogvprelog}
    Take a prime ideal $v\subset A$ of degree $1$ and a finite subset $I$ of $\bigcup_{r=1}^\infty(\tilde{\mathbb{D}}_\infty \cap \tilde{\mathbb{D}}_v)^r$. Assume that $\tilde{\mathbf{a}}=\big((s_1,\,a_1),\,\dots,\,(s_r,\,a_r)\big)\in I$ implies $\big((s_i,\,a_i),$ $\dots,$ $(s_j,\,a_j)\big)\in I$ for any $1 \leq i \leq j \leq r$.
    Then it follows that
    \begin{equation}
        \trdeg_{\Kalg(t)}\Kalg \left(t,\,\prelog_{v,\,\tilde{\mathbf{a}} }^{(1)}\,\middle|\,\tilde{\mathbf{a}} \in I \right)\geq \trdeg_{\Kalg(t)}\Kalg \left(t,\,\Omega_\infty,\,\prelog_{\infty,\,\tilde{\mathbf{a}}}^{(1)} \,\middle|\,\tilde{\mathbf{a}} \in I \right)-1.
    \end{equation}
\end{Corollary}
\begin{proof}
    We put $\mathbf{K}:=\Kalg(t)$, $\phi:=\tau^{-1}$, and let $\Phi$ be the block diagonal matrix
    $(t-\theta)\oplus \bigoplus_{\tilde{\mathbf{a}}\in I} \Phi_{\tilde{\mathbf{a}}}$
    where the block $\Phi_{\tilde{\mathbf{a}}}$ is the matrix in \cref{eq:def matrix T}. 
    Then the difference equation given by the pair $\big((\mathbf{K},\,\phi),\,\Phi \big)$ has a fundamental matrix
    $\Psi=(\Psi_{ij}):=(\Omega_\infty)\oplus\bigoplus_{\tilde{\mathbf{a}}\in I} \Psi_{\tilde{\mathbf{a}}} \in \GL(\Frac \mathbb{T}_\infty)$,
    where the blocks $\Psi_{\tilde{\mathbf{a}}}$ are as in \cref{Example:inftypolylogperiodinterpretation}, and
    the ring $\mathbf{R}:=\mathbf{K}[\Psi_{ij},\, 1/\det\Psi]$ is a simple difference domain by \cref{lem:PVsimple} and the discussion in \cref{Example:inftypolylogperiodinterpretation}.
    If we let $\mathbf{P}$ be $\mathbf{K} \left(\prelog_{v,\,
    \tilde{\mathbf{a}}}^{(1)} \,\middle|\,\tilde{\mathbf{a}} \in I \right)\subset \Frac\mathbb{T}_v$, then $\mathbf{P}$ is closed under $\phi:=\tau^{-1}$
    by \cref{eq:prelogvdifferenceeq} (see also \cref{Example:mathbfPv-adicpolylog}) and the assumption on the set $I$.
    Let $(\mathbf{S},\,\varphi)$ and $m$ be those yielded by \cref{Prop:PVTmain}.
    We consider the matrix $\Theta=(\Theta_{ij}):=(\Omega_\infty)\oplus\bigoplus_{\tilde{\mathbf{a}}\in I} \Theta_{\tilde{\mathbf{a}}}$, whose entries are regarded as elements of $\mathbf{S}$ (see \cref{Example:Theta_s}); then this matrix is a fundamental matrix of $\left((\mathbf{K},\,\phi^{m}),\,\phi^{m-1}(\Phi)\cdot\phi^{m-2}(\Phi)\cdots \phi(\Phi) \cdot \Phi \right)$
 since we have \cref{eq:Thetadiffeq} for each $\tilde{\mathbf{a}}\in I$.
    Therefore, \cref{Cor:PVTmain} implies 
    \begin{align}
        &\quad \ \trdeg_{\Kalg(t)}\Kalg \left(t,\,\prelog_{v,\,\tilde{\mathbf{a}} }^{(1)}\,\middle|\,\tilde{\mathbf{a}} \in I \right)=\trdeg_{\mathbf{K}}\mathbf{P}\\
        &=\trdeg_{\mathbf{K}}\mathbf{K}\left(\,\Omega_\infty,\,\prelog_{\infty,\,\tilde{\mathbf{a}}}^{(1)} \,\middle|\,\tilde{\mathbf{a}} \in I \right)-\trdeg_{\mathbf{P}}\mathbf{P}\left(\Omega_\infty,\,\prelog_{v,\,\tilde{\mathbf{a}} }^{(1)}\,\middle|\,\tilde{\mathbf{a}} \in I \right)\\
        &=\trdeg_{\mathbf{K}}\mathbf{K} \left(\Omega_\infty,\,\prelog_{\infty,\,\tilde{\mathbf{a}}}^{(1)} \,\middle|\,\tilde{\mathbf{a}} \in I \right)-\trdeg_{\mathbf{P}}\mathbf{P}(\Omega_\infty)\\
        &\geq\trdeg_{\mathbf{K}}\mathbf{K} \left(\Omega_\infty,\,\prelog_{\infty,\,\tilde{\mathbf{a}}}^{(1)} \,\middle|\,\tilde{\mathbf{a}} \in I \right)-1\\
        &=\trdeg_{\Kalg(t)}\Kalg \left(t,\,\Omega_\infty,\,\prelog_{\infty,\,\tilde{\mathbf{a}}}^{(1)} \,\middle|\,\tilde{\mathbf{a}} \in I \right)-1
    \end{align}
    as we have
    \begin{equation}
        \mathbf{K}( \Omega_\infty,\,\prelog_{v,\,\tilde{\mathbf{a}} }^{(1)}\,|\,\tilde{\mathbf{a}} \in I )=\mathbf{K}( \Theta_{ij} ) \text{ and } \mathbf{K}(\Omega_\infty,\,\prelog_{\infty,\,\tilde{\mathbf{a}}}^{(1)} \,|\,\tilde{\mathbf{a}} \in I)=\mathbf{K}(\Psi_{ij} )
    \end{equation}
    by the assumption on $I$.
\end{proof}

\section{Conclusions}\label{Section:conclusiongeneral}

In this section, we prove \cref{Thm:main_CMPLs}.
We first review the results of Chang, Chen, and Mishiba (\cite{Chang2022}, \cite{Chang2019a}, and \cite{Chang2021}) and show that they imply the inequality
\begin{equation}\label{eq:oppositeineq}
    \trdeg_{\Kalg}\Kalg \left(
        \mathrm{Li}(\tilde{\mathbf{a}})_{v}
        \,\middle|\,\tilde{\mathbf{a}} \in I
    \right)
    \leq
    \trdeg_{\Kalg}\Kalg \left(
        \tilde{\pi},\,\mathrm{Li}(\tilde{\mathbf{a}})_{\infty}
        \,\middle|\,\tilde{\mathbf{a}} \in I
    \right)-1.
\end{equation}
The reverse inequality follows by combining
\cref{Cor:inftyprepolylogpolylog,cor:vprelogpolylog,Cor:inftyprelogvprelog}.
We shall close the paper with some remarks on future work.

\subsection{The proof of the main result and algebraic independence consequences}\label{SS:opposite_ineq}
Before proving \cref{Thm:main_CMPLs}, we will confirm that Inequality \eqref{eq:oppositeineq} is already known from previous works.
We summarize the results of Chang, Chen, and Mishiba that we shall use, referring the reader to \cite{Chang2022}, \cite{Chang2019a}, and \cite{Chang2021}  for the full theory;
they focused on $\infty$-adic and $v$-adic \textit{Carlitz multiple star polylogarithms} rather than on the non-star versions, since the star versions are better adapted to the theory of $t$-modules and their logarithms.
Carlitz multiple star polylogarithms are defined in the same way as the Carlitz multiple polylogarithms, except that the range $i_1>i_2>\cdots>i_r\geq0$ in \eqref{eq:defCMPL} is replaced by $i_1\geq i_2\geq\cdots\geq i_r\geq0$; thus, the Carlitz multiple star polylogarithm $\mathrm{Li}_{\mathbf{s}}^{\star}(x_1,\,\dots,\,x_r)$ and the $v$-adic Carlitz multiple star polylogarithm $\mathrm{Li}_{\mathbf{s}}^{\star}(x_1,\,\dots,\,x_r)_v$ are given by
\begin{equation}\label{eq:defCMSPL}
    \sum_{i_1\geq i_2 \geq \cdots\geq i_r\geq 0}\frac{x_1^{(i_1)}\cdots x_r^{(i_r)}}{\Big(\left(\theta-\theta^{(1)}\right)\cdots\left(\theta-\theta^{(i_1)}\right)\Big)^{s_1}\cdots\Big(\left(\theta-\theta^{(1)}\right)\cdots\left(\theta-\theta^{(i_r)}\right)\Big)^{s_r}}
\end{equation}
for each tuple $\mathbf{s}:=(s_1,\,\dots,\,s_r)$ of positive integers; the star versions converge on the same domains as the non-star ones, namely $\mathbb{D}_{\mathbf{s},\,\infty}$ and $\mathbb{D}_{\mathbf{s},\,v}^{\mathrm{Conv}}$ (see \cref{eq:Dinfty,eq:Dv}).
Building on the idea in \cite{Anderson1990}, Chang and Mishiba (\cite[\S 6.1.2]{Chang2021}) extend the domain of the function $\mathrm{Li}_{\mathbf{s}}^{\star}(x_1,\,\dots,\,x_r)_v$ to
\begin{equation}
    \mathbb{D}_{\mathbf{s},\,v}^{\textrm{Def}}:=\left\{(x_1,\,\dots,\,x_r) \in \mathbb{C}_v^r\,\,\middle| \, |x_i|_v \leq 1 \text{ for } 1 \leq i \leq r \right\}.
\end{equation}

\begin{Remark}\label{Rem:starnonstar}
Carlitz multiple polylogarithms can be written in terms of Carlitz multiple star polylogarithms (\cite[Proposition 5.2.3]{Chang2021}); for example, we can observe, by definition, that
\begin{equation}
    \mathrm{Li}_{s_1,\,s_2}(x_1,\,x_2)=\mathrm{Li}_{s_1,\,s_2}^{\star}(x_1,\,x_2)-\mathrm{Li}_{s_1+s_2}^{\star}(x_1x_2) 
\end{equation}
and
\begin{align}
    \mathrm{Li}_{s_1,\,s_2,\,s_3}(x_1,\,x_2,\,x_3)&=\mathrm{Li}_{s_1,\,s_2,\,s_3}^{\star}(x_1,\,x_2,\,x_3)-\mathrm{Li}_{s_1+s_2,\,s_3}^{\star}(x_1x_2,\,x_3)\\
    &\qquad-\mathrm{Li}_{s_1,\,s_2+s_3}^{\star}(x_1,\,x_2x_3)+\mathrm{Li}_{s_1+s_2+s_3}^{\star}(x_1x_2x_3)
\end{align}
hold on the domains of definition; analogous identities for $r\geq 4$ follow in the same way.
We further notice that $v$-adic Carlitz multiple polylogarithms and $v$-adic Carlitz multiple star polylogarithms satisfy the same relations on their domain of convergence.
\end{Remark}
Chang, Chen, and Mishiba proved that the $\Kalg$-span
\begin{equation}
    \overline{\mathcal{L}}_\infty^{\textrm{Conv}}:= \left\langle 1,\,\mathrm{Li}_{\mathbf{s}}^{\star}(\mathbf{x}) \,\,\middle|\,\,r\geq 1,\,s_1,\,\dots,\,s_r \geq 1,\,\mathbf{x}\in \mathbb{D}_{\mathbf{s},\,\Kalg}^{\textrm{Conv}}
   \right\rangle_{\Kalg} \subset \mathbb{C}_\infty
\end{equation}
is a $\Kalg$-subalgebra (\cite[Proposition 2.2.3]{Chang2022}) and coincides with
\begin{equation}
    \overline{\mathcal{L}}_\infty^{\textrm{Def}}:= \left\langle 1,\,\mathrm{Li}_{\mathbf{s}}^{\star}(\mathbf{x}) \,\,\middle|\,\,r\geq 1,\,s_1,\,\dots,\,s_r \geq 1,\,\mathbf{x}\in \mathbb{D}_{\mathbf{s},\,\Kalg}^{\textrm{Def}}
   \right\rangle_{\Kalg} \subset \mathbb{C}_\infty
\end{equation}
(\cite[Theorem 4.2.1]{Chang2022}), where we put
\begin{align}
    \mathbb{D}_{\mathbf{s},\,\Kalg}^{\textrm{Conv}}&:=\left( \mathbb{D}_{\mathbf{s},\,\infty} \cap (\Kalg)^r \right) \cap\left( \mathbb{D}_{\mathbf{s},\,v}^{\textrm{Conv}} \cap (\Kalg)^r \right) \text{ and}\\
    \mathbb{D}_{\mathbf{s},\,\Kalg}^{\textrm{Def}}&:=\left( \mathbb{D}_{\mathbf{s},\,\infty} \cap (\Kalg)^r \right) \cap\left( \mathbb{D}_{\mathbf{s},\,v}^{\textrm{Def}} \cap (\Kalg)^r \right) 
\end{align}
for each $r$-tuple $\mathbf{s}$ of positive integers; 
they further constructed a well-defined $\Kalg$-algebra homomorphism 
\begin{equation}\phi_v:\overline{\mathcal{L}}_\infty^{\textrm{Def}} \twoheadrightarrow \overline{\mathcal{L}}_v^{\textrm{Def}}:= \left\langle 1,\,\mathrm{Li}_{\mathbf{s}}^{\star}(\mathbf{x})_v \,\,\middle|\,\,r\geq 1,\,s_1,\,\dots,\,s_r \geq 1,\,\mathbf{x}\in \mathbb{D}_{\mathbf{s},\,\Kalg}^{\textrm{Def}}
   \right\rangle_{\Kalg} \subset \mathbb{C}_v 
\end{equation}
given by $\mathrm{Li}_{\mathbf{s}}^{\star}(\mathbf{x}) \mapsto\mathrm{Li}_{\mathbf{s}}^{\star}(\mathbf{x})_v$ (\cite[Theorem 5.1.5 and \S 5.2]{Chang2022}); this homomorphism $\phi_v$ is distinct from the difference automorphism $\phi$ of \Cref{sec:difference galois}; as a consequence, we have an equality
\begin{equation}
   \overline{\mathcal{L}}_v^{\textrm{Def}}=\overline{\mathcal{L}}_v^{\textrm{Conv}}:= \left\langle 1,\,\mathrm{Li}_{\mathbf{s}}^{\star}(\mathbf{x})_v \,\,\middle|\,\,r\geq 1,\,s_1,\,\dots,\,s_r \geq 1,\,\mathbf{x}\in \mathbb{D}_{\mathbf{s},\,\Kalg}^{\textrm{Conv}}
   \right\rangle_{\Kalg} \subset \mathbb{C}_v,
\end{equation}
see \cite[\S 5.2]{Chang2022}.
By \cref{Rem:starnonstar}, we have $\mathrm{Li}_{\mathbf{s}}(\mathbf{x}) \in \overline{\mathcal{L}}_\infty^{\textrm{Conv}}=\overline{\mathcal{L}}_\infty^{\textrm{Def}} $ and $\phi_v(\mathrm{Li}_{\mathbf{s}}(\mathbf{x}))=\mathrm{Li}_{\mathbf{s}}(\mathbf{x})_v$ for any tuple $\mathbf{s}$ and $\mathbf{x}\in  \mathbb{D}_{\mathbf{s},\,\Kalg}^{\textrm{Conv}}$; at a point $\mathbf{x}\in  \mathbb{D}_{\mathbf{s},\,\Kalg}^{\textrm{Def}}$ that is not necessarily in $\mathbb{D}_{\mathbf{s},\,\Kalg}^{\textrm{Conv}}$, the value of the $v$-adic Carlitz multiple polylogarithm is given by $\mathrm{Li}_{\mathbf{s}}(\mathbf{x})_v:=\phi_v(\mathrm{Li}_{\mathbf{s}}(\mathbf{x}))$ (\cite[Definition 4.2.3]{Chang2019a}).

The \textit{Carlitz zeta values}
\begin{equation}
    \zeta_A(s)
    :=
    \sum_{\substack{a\in A\\ a\text{ monic}}}\frac{1}{a^{s}} \in K_\infty \quad (s \in \mathbb{Z},\,s \geq1),
\end{equation}
which were introduced by Carlitz (\cite{Carlitz1935}), can be written as $K$-linear combinations of special values of $\mathrm{Li}_{s}(x)$ at elements of $A$ (\cite[\S 3.7, \S 3.8, and \S 3.9]{Anderson1990}); thus these values are elements of  $\overline{\mathcal{L}}_\infty^{\textrm{Def}}=\overline{\mathcal{L}}_\infty^{\textrm{Conv}}$.
Furthermore, the value $\zeta_A(s)$ is a product of $\tilde{\pi}^s$ and an explicit element of $K$ if $s$ is a multiple of $q-1$, that is,
\begin{equation}\label{eq:Euler--Carlitz}
    \zeta_A(s)/\tilde{\pi}^s\in K
\end{equation} 
for each multiple $s$ of $q-1$; this is known as the Euler--Carlitz relation (\cite{Carlitz1935}, \cite[Theorem 5.2.1]{ThakurBook}).
For instance, one has
\begin{equation}\label{eq:eulerian-q-1}
    \mathrm{Li}_{(q-1)}(1)
    =
    \zeta_A(q-1)
    =
    \frac{\tilde{\pi}^{q-1}}{\theta-\theta^q}.
\end{equation}
Since this value is in $\overline{\mathcal{L}}_\infty^{\textrm{Conv}}$, it holds that
\begin{equation}\label{eq:ring-identification-q-1}
    \begin{aligned}\Kalg\bigl[\mathrm{Li}(\tilde{\mathbf{a}})_\infty \mid \tilde{\mathbf{a}}\in I\cup\{((q-1,\,1))\}\bigr]
    &= 
     \Kalg\bigl[\zeta_A(q-1),\,\mathrm{Li}(\tilde{\mathbf{a}})_\infty \mid \tilde{\mathbf{a}}\in I\bigr]\\
   &= \Kalg\bigl[\tilde{\pi}^{q-1},\,\mathrm{Li}(\tilde{\mathbf{a}})_\infty \mid \tilde{\mathbf{a}}\in I\bigr]
    \end{aligned}
\end{equation}
for each subset $I$ of $\bigcup_{r=1}^\infty(\tilde{\mathbb{D}}_\infty \cap \tilde{\mathbb{D}}_v)^r$ with the notation in \cref{eq:polylogabb}.
By \cref{Rem:starnonstar}, this $\Kalg$-algebra is a subalgebra of $\overline{\mathcal{L}}_\infty^{\textrm{Def}}=\overline{\mathcal{L}}_\infty^{\textrm{Conv}}$; hence we have a restriction 
\begin{equation}
\phi_v^\prime:\Kalg\Bigl[\zeta_A(q-1),\,\mathrm{Li}(\tilde{\mathbf{a}})_\infty \mid \tilde{\mathbf{a}}\in I\Bigr] \twoheadrightarrow \Kalg\Bigl[\phi_v\big(\zeta_A(q-1)\big),\,\mathrm{Li}(\tilde{\mathbf{a}})_v \mid \tilde{\mathbf{a}}\in I\Bigr].
\end{equation}
Since $\phi_v\big(\zeta_A(q-1)\big)=0$ (see \cite[Theorem 1.2.3]{Chang2022} and \cite[Theorem 5.3]{Goss1979} for example), the kernel of $\phi_v^\prime$ contains the non-zero element $\zeta_A(q-1)$; thus, if the set $I$ is finite, we have the inequality
\begin{equation}\label{eq:opposite_ineq}
    \begin{aligned}
    \dim \Kalg[\tilde{\pi}^{q-1},\,\mathrm{Li}(\tilde{\mathbf{a}})_\infty \mid \tilde{\mathbf{a}} \in I ]
    &=\dim \Kalg[\zeta_A(q-1),\,\mathrm{Li}(\tilde{\mathbf{a}})_\infty \mid \tilde{\mathbf{a}} \in I ]\\
    &\geq \dim \Kalg[\mathrm{Li}(\tilde{\mathbf{a}})_v \mid \tilde{\mathbf{a}} \in I ]+1
    \end{aligned}
\end{equation}
between the Krull dimensions of finitely generated integral algebras over a field, by which we obtain the inequality \eqref{eq:oppositeineq} as $\tilde{\pi}$ is algebraic over $\Kalg(\tilde{\pi}^{q-1})$.
\begin{proof}[Proof of \cref{Thm:main_CMPLs}]
    Let us take a prime ideal $v\subset A$ of degree $1$ and a finite subset $I$ of $\bigcup_{r=1}^\infty(\tilde{\mathbb{D}}_\infty \cap \tilde{\mathbb{D}}_v)^r$. 
    Suppose that $\tilde{\mathbf{a}}=\big((s_1,\,a_1),\,\dots,\,(s_r,\,a_r)\big)\in I$ implies $\big((s_i,\,a_i),$ $\dots,$ $(s_j,\,a_j)\big)\in I$ for any $1 \leq i \leq j \leq r$.
    Then, \cref{Cor:inftyprepolylogpolylog,Cor:inftyprelogvprelog,cor:vprelogpolylog} imply the inequality
    \begin{equation}
        \trdeg_{\Kalg}\Kalg \left(\mathrm{Li}(\tilde{\mathbf{a}})_{v}\,\middle|\,\tilde{\mathbf{a}} \in I \right)\geq\trdeg_{\Kalg}\Kalg \left(\tilde{\pi},\,\mathrm{Li}(\tilde{\mathbf{a}})_{\infty}\,\middle|\,\tilde{\mathbf{a}} \in I \right)-1.
    \end{equation}
    The opposite inequality has already been confirmed.
\end{proof}

Combining \cref{Thm:main_CMPLs} with existing algebraic independence results for special values of $\infty$-adic Carlitz multiple polylogarithms, we can obtain the following corollaries.

\begin{Corollary}
    Take a prime ideal $v \subset A$ of degree $1$ and natural numbers $s_1,\,\dots,$ $s_d$ not divisible by $q-1$, and assume that the ratio $s_i/s_j$ is not an integer power of $p$ for each $1 \leq i < j \leq d$. 
    For any $a_1,\,\dots,\,a_d \in A \setminus 0$ such that $|a_i|_\infty<|\theta|_\infty^{\frac{s_iq}{q-1}}$ and $|a_i|_v<1$ for each $1 \leq i \leq d$, we have
    \begin{equation}
        \trdeg_{\Kalg}\Kalg(\mathrm{Li}(\tilde{\mathbf{a}})_{v} \mid \tilde{\mathbf{a}}\in I)=\frac{d(d+1)}{2}
    \end{equation}
    where $I:=\{ \tilde{\mathbf{a}}=((s_i,\,a_i),\,\dots,\,(s_j,\,a_j))\mid 1\leq i \leq j \leq d\}$.
\end{Corollary}
\begin{proof}
    \cref{Thm:main_CMPLs} yields the equality
    \begin{equation}
        \trdeg_{\Kalg}\Kalg(\mathrm{Li}(\tilde{\mathbf{a}})_{v} \mid \tilde{\mathbf{a}}\in I)=\trdeg_{\Kalg}\Kalg(\tilde{\pi},\,\mathrm{Li}(\tilde{\mathbf{a}})_{\infty} \mid \tilde{\mathbf{a}}\in I)-1.
    \end{equation}
    Because of the assumptions on the natural numbers $s_1,\,\dots,\,s_d$, Theorem 1.5 in \cite{Mishiba2017} shows that the transcendence degree on the right-hand side is $1+d(d+1)/2$; hence the result follows.
\end{proof}

\begin{Corollary}
    Fix a prime ideal $v \subset A$ of degree $1$ and natural numbers $s_1,\,\dots,\,s_d$ such that the ratio $s_i/s_j$ is not an integer power of $p$ for each $1 \leq i < j \leq d$. 
    For each $1 \leq i \leq d$, further take any $a_{i1},\,\dots,\,a_{ir_i} \in A \setminus 0$ such that $|a_{ij}|_\infty<|\theta|_\infty^{\frac{s_iq}{q-1}}$ and $|a_{ij}|_v<1$ for each $1 \leq j \leq r_i$.
    If, for each $1 \leq i \leq d$, the special values
    \begin{equation}        
        \mathrm{Li}((s_i,\,a_{i1}))_{v},\,\dots,\,\mathrm{Li}((s_i,\,a_{ir_i}))_{v}         
    \end{equation}
    are linearly independent over $K$, then we have
    \begin{equation}
        \trdeg_{\Kalg}\Kalg (\mathrm{Li}((s_i,\,a_{ij}))_{v} \mid 1\leq i \leq d,\,1 \leq j \leq r_i)=\sum_{i=1}^d r_i.
    \end{equation}
\end{Corollary}
\begin{proof}
    Fix $1 \leq i \leq d$ and first consider the case where $(q-1) \nmid s_i $, which forces $q > 2$; then we notice that $\tilde{\pi}^{s_i} \notin K_\infty$.
    Indeed, it is immediate from the definition of $\tilde{\pi}$ (\cref{eq:def_pitilde}) that $v_\infty(\tilde{\pi}) = -q/(q-1)$; thus $\tilde{\pi}^{s} \in K_\infty$ only if $(q-1) \mid s$, since $v_\infty(K_\infty^\times) = \mathbb{Z}$; the same observation appears in \cite[proof of Cor.~2.3.3]{Chang2014}.
    On the other hand, the special values $\mathrm{Li}((s_i,\,a_{i1}))_{\infty},\,\dots,\,\mathrm{Li}((s_i,\,a_{ir_i}))_{\infty} $  are elements of $K_\infty$. Since we obtain, by restricting $\phi_v$ considered in \cite[\S 5.2]{Chang2022}, a $\Kalg$-linear map
    \begin{equation}
        \langle \mathrm{Li}((s_i,\,a_{i1}))_{\infty},\,\dots,\,\mathrm{Li}((s_i,\,a_{ir_i}))_{\infty} \rangle_{\Kalg} \rightarrow \langle  \mathrm{Li}((s_i,\,a_{i1}))_{v},\,\dots,\,\mathrm{Li}((s_i,\,a_{ir_i}))_{v}  \rangle_{\Kalg}
    \end{equation}
    given by $\mathrm{Li}((s_i,\,a_{ij}))_{\infty} \mapsto \mathrm{Li}((s_i,\,a_{ij}))_{v}$ for $1 \leq j \leq r_i$, we can deduce the $K$-linear independence of 
    \begin{equation}
        \mathrm{Li}((s_i,\,a_{i1}))_{\infty},\,\dots,\,\mathrm{Li}((s_i,\,a_{ir_i}))_{\infty}
    \end{equation}
    from that of $\mathrm{Li}((s_i,\,a_{i1}))_{v},\,\dots,\,\mathrm{Li}((s_i,\,a_{ir_i}))_{v} $, and hence the elements
    \begin{equation}
        \tilde{\pi}^{s_i},\,\mathrm{Li}((s_i,\,a_{i1}))_{\infty},\,\dots,\,\mathrm{Li}((s_i,\,a_{ir_i}))_{\infty} 
    \end{equation}
    are linearly independent over $K$ as $\tilde{\pi}^{s_i} \notin K_\infty$.
    Second, we consider the case $(q-1)|s_i$; works of Chang, Chen, and Mishiba (\cite[Theorem 1.2.3]{Chang2022}) also yield the $\Kalg$-linear map
    \begin{align}
       & \langle \zeta_A(s_i),\,\mathrm{Li}((s_i,\,a_{i1}))_{\infty},\,\dots,\,\mathrm{Li}((s_i,\,a_{ir_i}))_{\infty} \rangle_{\Kalg} \\
        &\quad\rightarrow \langle  \mathrm{Li}((s_i,\,a_{i1}))_{v},\,\dots,\,\mathrm{Li}((s_i,\,a_{ir_i}))_{v}  \rangle_{\Kalg}
    \end{align}
    given by $\zeta_A(s_i)\mapsto 0$ and $\mathrm{Li}((s_i,\,a_{ij}))_{\infty} \mapsto \mathrm{Li}((s_i,\,a_{ij}))_{v}$ for $1 \leq j \leq r_i$, which enables us to establish the linear independence of 
    \begin{equation}
        \tilde{\pi}^{s_i},\,\mathrm{Li}((s_i,\,a_{i1}))_{\infty},\,\dots,\,\mathrm{Li}((s_i,\,a_{ir_i}))_{\infty} 
    \end{equation}
    as \cref{eq:Euler--Carlitz} holds and we have assumed that \begin{equation}        
        \mathrm{Li}((s_i,\,a_{i1}))_{v},\,\dots,\,\mathrm{Li}((s_i,\,a_{ir_i}))_{v}         
    \end{equation}
    are linearly independent over $K$.
    Now, Theorem 4.2 in \cite{Mishiba2017} implies
    \begin{equation}
        \trdeg_{\Kalg}\Kalg (\tilde{\pi},\,\mathrm{Li}((s_i,\,a_{ij}))_{\infty} \mid 1\leq i \leq d,\,1 \leq j \leq r_i)=1+\sum_{i=1}^d r_i;
    \end{equation}
    combining this equality with  
    \begin{align}
        &\trdeg_{\Kalg}\Kalg (\mathrm{Li}((s_i,\,a_{ij}))_{v} \mid 1\leq i \leq d,\,1 \leq j \leq r_i)\\
        &\qquad =\trdeg_{\Kalg}\Kalg (\tilde{\pi},\,\mathrm{Li}((s_i,\,a_{ij}))_{\infty} \mid 1\leq i \leq d,\,1 \leq j \leq r_i)-1 
    \end{align}
    given by \cref{Thm:main_CMPLs}, we obtain the desired result.
\end{proof}

\subsection{Future research}
We conclude the present paper by discussing directions for future research.
It is natural to attempt to remove the conditions on the degree of the place $v$ and on the $v$-adic valuations of arguments imposed in \cref{Thm:main_CMPLs}. 
Regarding the restriction on the $v$-adic valuations of the arguments, it would be interesting to investigate the relationship between the values dealt with in \cref{Example:comparison-non-conv,Example:preperiods-non-conv} and the corresponding special values of the analytic continuations of $v$-adic Carlitz (multiple) polylogarithms that appear in \cite[Definitions 4.1.2 and 4.2.3]{Chang2019a}.

A further question concerns the relationship between the formulation in the present paper, or a suitable generalization of it, and the approach of Chang, Wei, and Yu \cite[\S 3.7]{Chang2025}, which focuses on the matrix representing the Frobenius action induced on the de Rham realization.
We have discussed a special case in \cref{Rem:CWY}.
For related questions in the setting of $p$-adic periods, see \cite[Remarks 5.6 and 5.17]{Ancona2025} and \cite{DanCohen2025}.

It would also be interesting to describe the transcendence properties of $v$-adic periods in terms of Tannakian duality and algebraic groups.
Ancona and Fr\u{a}\c{t}il\u{a} formulated a conjecture relating the transcendence properties of $p$-adic periods to dimensions of motivic Galois groups, which are certain algebraic groups associated with the chosen motives via Tannakian duality; see \cite[Theorem 6.3 and \S 8]{Ancona2025} for details.
We hope to pursue an analogue in positive characteristic of the theory of Ancona and Fr\u{a}\c{t}il\u{a} in a future project.

%%%%%%%%%%%%%%%
%\bibliographystyle{abbrv}    
%\bibliography{BibTexLibrary1.bib} 

\end{document}